\newcommand{\R}{I\!\!R}
\newcommand{\Z}{\mathbb Z}
\newcommand{\Q}{\mathbb Q}
\newcommand{\N}{I\!\!N}
\newcommand{\C}{\mathbb C}
\newcommand{\Id}{\mathrm{Id}}
\newcommand{\Dim}{\mathrm{dim}}
\newcommand{\Codim}{\mathrm{codim}}
\newcommand{\Ker}{\mathrm{Ker}}
\newcommand{\Bsym}{\mathrm B_{\mathrm{sym}}}
\newcommand{\Lsym}{\mathcal L_{\mathrm{sym}}}
\newcommand{\sgn}{\sigma}
\newcommand{\ord}{\mathrm{ord}}
\newcommand{\Imm}{\mathrm{Im}}
\newcommand{\mul}{\mathrm{mul}}
\newcommand{\Gr}{\mathrm{Gr}}
\newcommand{\LLL}{\Sigma_{L_0}}

\newcommand{\noo}[1]{\overset {\mbox{%
\lower1pt\hbox{${\scriptscriptstyle o}$}}}n^{\mbox{%
\lower2pt\hbox{$\scriptscriptstyle #1$}}}}

\newcommand{\vv}{\mathfrak v}

\newcommand{\LL}{{\mathbf L}}
\newcommand{\Ss}{{\mathbf S}}

\newcommand{\Spl}{\mathrm{Sp}}

\newcommand{\iMaslov}{\mathrm i_{\scriptscriptstyle{\textrm{Maslov}}}}
\newcommand{\ispec}{\mathrm i_{\scriptscriptstyle{\textrm{spectral}}}}
\newcommand{\iMorse}{\mathrm i_{\scriptscriptstyle{\textrm{Morse}}}}
\newcommand{\ind}{\mathrm{ind}}

\newcommand{\Scal}{\mathcal S}

\newcommand{\Fcal}{\mathcal F}
\newcommand{\Lcalsa}{\mathcal L^{\scriptscriptstyle{\textrm{sa}}}}
\newcommand{\Fcalsa}{\mathcal F^{\scriptscriptstyle{\textrm{sa}}}}

\newcommand{\spfl}{\mathrm{sf}}
\newcommand{\Hcal}{\mathcal H}
\newcommand{\Dcal}{\mathcal D}

\newcommand{\ddto}[1]{\left.\dfrac{\mathrm d^{#1}}{\mathrm
dt^{#1}}\right\vert_{t=0}}

\newcommand{\tdto}[1]{\tfrac{\mathrm d^{#1}}{\mathrm
dt^{#1}}\big\vert_{t=0}}

\newcommand{\tdtoo}[1]{\tfrac{\mathrm d^{#1}}{\mathrm
dt^{#1}}\big\vert_{t=t_0}}

\newcommand{\gr}{g^{\scriptscriptstyle+}}

\newcommand{\llangle}{\langle\!\!\!\langle}
\newcommand{\rrangle}{\rangle\!\!\!\rangle}

\documentclass[oneside,draft,a4paper,10pt]{amsart}

\usepackage{times} 
\usepackage{epsfig}

\numberwithin{equation}{section}

\title[On the Maslov index in the degenerate
case]{On the Maslov index
of Lagrangian paths that are not transversal to the
Maslov cycle.\\
Semi-Riemannian  index theorems in the degenerate~case.}
\author[R.\ Giamb\`o]{Roberto Giamb\`o}
\address{Dipartimento di Matematica e
Informatica\hfill\break\indent Universit\`a di Camerino
\hfill\break\indent Camerino, MC, Italy}
\email{roberto.giambo@unicam.it}
\author[P.\ Piccione]{Paolo Piccione}
\address{Departamento de Matem\'atica,\hfill\break\indent
Instituto de Matem\'atica e Estat\'\i stica\hfill\break\indent  Universidade de S\~ao
Paulo, \hfill\break\indent
Rua do Mat\~ao 1010, CEP 05508-900, S\~ao Paulo, SP\hfill\break\indent Brazil}
\email{piccione@ime.usp.br}
\urladdr{http://www.ime.usp.br/\~{}piccione}
\curraddr{Dipartimento di Matematica e Informatica,
Universit\`a di Camerino, Italy.}
\author[A.\ Portaluri]{Alessandro Portaluri}

\thanks{The research work exposed in this paper originated
from a suggestion given by the third author relating the spectral
flow with the partial signatures.
He pointed out several references and formulated the statement
that in the invertible endpoints case, the spectral flow of
a real analytic path of self-adjoint Fredholm operators is given by
the sum of the odd partial signatures at each degeneracy instant.
In the present version of the article, this statement is
part of Proposition~\ref{thm:central}; totally, the contribution
of the third author to the theory consists in Definition~\ref{thm:defroot},
formulas~\eqref{eq:Wk} and \eqref{eq:Bk}, the first statement of
Proposition~\ref{thm:proprbas}, parts of Remark~\ref{thm:remuseless},
Definition~\ref{thm:defpartsign}, and parts of the statement and
parts of the proof of Proposition~\ref{thm:central}. Apart from this original
contribution, the material contained in this paper was entirely
developed and written by the first two authors at the \emph{Universit\`a di Camerino}
(Italy), and at the \emph{Universidade de S\~ao Paulo} (Brazil).}

\subjclass[2000]{37G25, 37J05, 47A07, 47A53, 53C22, 53C50, 58E10}

\date{October 26th, 2004}

\begin{document}


\theoremstyle{plain}\newtheorem{teo}{Theorem}[section]
\theoremstyle{plain}\newtheorem{prop}[teo]{Proposition}
\theoremstyle{plain}\newtheorem{lem}[teo]{Lemma}
\theoremstyle{plain}\newtheorem{cor}[teo]{Corollary}
\theoremstyle{definition}\newtheorem{defin}[teo]{Definition}
\theoremstyle{remark}\newtheorem{rem}[teo]{Remark}
\theoremstyle{plain} \newtheorem{assum}[teo]{Assumption}
\theoremstyle{definition}\newtheorem{example}[teo]{Example}


\begin{abstract}
The Maslov index of a Lagrangian path, under a certain
transversality assumption, is given by an algebraic count of the
intersections of the path with a  subvariety of the Lagrangian
Grassmannian called the Maslov cycle. In these notes we use the
notion of {\em generalized signatures\/} at a singularity of a
smooth curve of symmetric bilinear forms to determine a formula
for the computation of the Maslov index in the case of a
real-analytic path having possibly non transversal intersections.
Using this formula we give a general definition of Maslov index
for continuous curves in the Lagrangian Grassmannian, both in the
finite and in the infinite dimensional (Fredholm) case, and having
arbitrary endpoints. Other notions of Maslov index are also
considered, like the index for pairs of Lagrangian paths, the
Kashiwara's triple Maslov index, and H\"ormander's four-fold
index.

We discuss some applications
of the theory, with special emphasis on the study of the
Jacobi equation along a semi-Riemannian geodesic. In this context, we prove an
extension of several versions of the Morse index theorems for geodesics
having possibly conjugate endpoints.

The main results of this paper were announced in
\cite{GiaPicPorCOMPTES}.
\end{abstract}

\maketitle
\tableofcontents


\typeout{}
\typeout{}
\typeout{______________________________________________________________________}
\typeout{}
\typeout{This is the LaTeX file of the notes: }
\typeout{}
\typeout{ 'On the Maslov index
of Lagrangian paths}
\typeout{  that are not transversal to the
Maslov cycle.}
\typeout{ Semi-Riemannian  index theorems in the degenerate case.'}
\typeout{}
\typeout{ by R. Giambo', P. Piccione and A. Portaluri}
\typeout{}
\typeout{Questions, comments or remarks should be addressed to Paolo Piccione}
\typeout{by either e-mail:}
\typeout{}
\typeout{ piccione@ime.usp.br}
\typeout{}
\typeout{or by regular mail:}
\typeout{}
\typeout{ Instituto de Matematica e Estatistica}
\typeout{ Universidade de Sao Paulo}
\typeout{ Caixa Postal 66281, Agencia Cidade de Sao Paulo}
\typeout{ Sao Paulo, SP, Brazil, CEP 05315-970.}
\typeout{}
\typeout{______________________________________________________________________}
\typeout{}


\begin{section}{Introduction}\label{sec:intro}

Recent extensions of classical variational theories, like for instance
Morse Theory, Ljusternik--Schnirelman theory  or bifurcation theory,
have shown the increasing relevance of topological techniques in finite or
infinite dimensional manifolds.
Solutions of different geometrical variational problems, like
geodesics, minimal surfaces, harmonic maps, solutions of Hamiltonian
systems, etc., are classified by an {\em index\/}, typically an integer number,
that carries both analytical and geometrical information on the solution.
Jumps of the appropriate
index function (Morse, Conley, Maslov, spectral, etc.), that
can only occur at degeneracy instants of the second variation
form,  are responsible for
topologically  non trivial changes in the solution space
of the given variational problem. In the  non positive definite  case,
a first order analysis of the degeneracy instant is no longer sufficient to
detect a jump of the index and to compute its value, so one
needs a more accurate analysis beyond the first derivative test,
capable of reckoning fine analytical details of the degeneracy.

In the symplectic world, a recurrent notion is that of Maslov
index, that appears naturally in many different contexts,
especially in relation with solutions of Hamiltonian systems. The
natural environment for the notion of Maslov index is the
Lagrangian Grassmannian $\Lambda$ of a finite or infinite
dimensional symplectic space; the Maslov index is a $\Z$-valued
homotopical invariant for continuous curves in $\Lambda$ that
gives an algebraic measure of the intersections with the Maslov
cycle, which is an algebraic subvariety of $\Lambda$. Under
generic circumstances, the intersections of a curve with the
Maslov cycle occur at its regular part, and they are transversal
(hence isolated). In this case  the computation of the Maslov
index is done via well established results obtained by
differential topological methods, and its relations with the
geometrical and analytical invariants  of the variational problem
are clear. Typically, the non transversal case is studied by
perturbative techniques, which allow to extend to this case the
results involving quantities that are stable by uniformly small
perturbations.

There are several reasons to develop a non perturbative analysis
of the non transversal intersections, and that motivated the
research exposed in this paper. In first place, perturbation
arguments do not work properly when non transversal intersections
occur at the endpoints; namely, in this case arbitrarily small
perturbations may destroy the intersection. Observe that the
arguments needed to prove the genericity of the transversality
property in specific contexts become drastically more involved (if
not impossible at all) if one  restricts to fixed endpoints
homotopies. For instance, in the case of periodic solutions of
Hamiltonian systems, in order to obtain multiplicity results one
needs to consider iteration formulas for the Maslov index, and it
is not evident whether transversality is generic in spaces of
inextensible paths.

Second, and more important, perturbative methods preserve global
quantities, but destroy the information concerning each single
intersection, which of course may be relevant in the problem under
consideration. This is the case, for instance, in bifurcation
theory, if one wants to know whether a given instant bifurcates or
not, and not just whether bifurcation occurs {\em somewhere}.
Also, in semi-Riemannian geometry the presence of conjugate
points, that correspond to intersections with the Maslov cycle
obtained from the flow of the Jacobi deviation equation, have
global geometrical implications on their own, and perturbing the
data (i.e., the metric) would not be a very meaningful procedure.
Along a semi-Riemannian geodesic, non transversal intersections
with the Maslov cycle (that may occur only if the metric is not
positive definite) correspond to degenerate conjugate points; the
presence of this kind of conjugate points is responsible for a
series of new and interesting phenomena in the semi-Riemannian
vs.\ the Riemannian world, hence it deserves a specific analysis.
The occurrence of degenerate intersections arising from
semi-Riemannian geodesics has been somewhat overlooked in the
classical literature, and only recently the relevance of such
occurrences has been recognized; for a more detailed discussion on
this topic, see for instance \cite[Subsection~5.4]{pacific} and
\cite{fechado}. The result of these notes contributes
understanding such degeneracy phenomena in the real-analytic case;
for instance, we give conditions under which the exponential map
is not one-to-one on any neighborhood of a possibly degenerate
semi-Riemannian conjugate point (Corollary~\ref{thm:expnot1to1}),
extending a classical Riemannian result of Morse and Littauer (see
\cite{War}).

In this paper we will be essentially interested in the notion of
Maslov index associated to geodesics  and, more generally, to
solutions of Hamiltonian systems.  The Hamiltonian flow on a
symplectic manifold preserves the symplectic  form, and by a
suitable trivialization of the tangent bundle of the manifold
along a given solution we get a smooth curve in the symplectic
group of a fixed symplectic space. Likewise, the evolution of
Lagrangian initial conditions by the flow of the Hamilton equation
produces a curve in the Grassmannian of all Lagrangian subspaces
of a fixed symplectic space. The local geometry of the Lagrangian
Grassmannian of a symplectic space, both in finite and the
infinite dimensional case, is the geometry of symmetric bilinear
forms, and intersections with the Maslov cycle correspond to
degeneracy instants of the forms; non transversality with the
Maslov cycle corresponds, in local charts, to the fact that the
corresponding curve $B(t)$ of symmetric bilinear forms has
derivative $B'(t_0)$ which is degenerate when restricted to
$\Ker\big(B(t_0)\big)$. The Maslov index of a curve measures at
each intersection the change of value of the index of the
symmetric form, and a computation of this jump is studied using
higher order derivatives and the introduction of the notions of
``generalized Jordan chains'' and of ``partial signatures'' at an
isolated degeneracy instant. Given a smooth curve in the
Lagrangian Grassmannian having only isolated intersections with
the Maslov cycle, using local charts we associate a sequence of
integers to each such intersection, and we give a formula to
compute the Maslov index of the curve in terms of these integers.
The theory works well under the assumption that the curve be
real-analytic, or, more generally, when it is possible to find a
smoothly varying Hilbert basis of the space that diagonalizes the
curve of symmetric forms. Important examples for the theory of
Maslov index of real-analytic curves are encountered  in  the
study of eigenvalue problems for some differential operators,
whose solutions depend analytically on the eigenvalue; in this
paper we will discuss the theory in the case of the
Morse--Sturm--Liouville equation.


The partial signatures at a degeneracy instant have appeared in
several different contexts in the literature, associated to
``jumps'' of integer valued invariants, like the spectral flow, or
the so-called eta-invariant  associated to an elliptic
self-adjoint operator. It is not an easy task to establish where
exactly the notion of partial signature at an isolated degeneracy
instant of a path of symmetric (or hermitian) forms has first
appeared. To the authors' knowledge, one of the first references
where such notion appears is \cite{GohSig}, where the partial
signatures are introduced as generalized multiplicities of
characteristic values of a meromorphic operator valued function.

The partial signatures appear in a paper of Rabier
(ref.~\cite{rabier}) as an evolution of Magnus' generalized
algebraic multiplicity for nonlinear eigenvalues of operators
between Banach spaces, aiming at results in bifurcation theory.
Our approach follows more closely the paper by Farber and Levine
\cite{FarLev}, where the partial signatures have been used to
determine a formula for the jumps of the so-called {\em
eta-invariant}. The link between the Maslov index and the theory
of partial signatures has been suggested to the authors by some
recent results of Fitzpatrick, Pejsachowicz and Recht on the
bifurcation theory for strongly indefinite functionals (see
\cite{FitzPejsaRecht}). The main result of Rabier in \cite{rabier}
relates bifurcation phenomena to generalized Jordan chains for a
sequence of symmetric linear operators on finite dimensional
vector spaces. The central result of \cite{FitzPejsaRecht} is,
roughly speaking, that bifurcation for a smooth family of
(strongly indefinite) functionals having a trivial branch of
critical point at which the Hessians are Fredholm, occurs at those
instants when the {\em spectral flow\/} of the path of second
variations jumps. Finally, a geometrical bifurcation problem in
the context of semi-Riemannian geodesics has been recently studied
in \cite{geobif}, where it is proven that the Maslov index of a
geodesic equals the spectral flow of the corresponding curve of
index forms.

As an application of our theory, we will discuss a very general
version of the semi-Riemannian Morse index theorem
(Theorem~\ref{thm:tuttouguale}), relating the spectral data of the
Jacobi differential operator with the spectral flow of the index
form along the geodesic, and with the conjugate points along the
geodesic. The theorem gives an equality of three different
indexes, called the {\em spectral index}, the {\em generalized
Morse index\/} and the {\em Maslov index\/} of the geodesic, which
is shown to hold without any assumption on the final endpoint of
the geodesic. The notion of spectral index, which plays a central
role in our theory, is defined as the spectral flow of the curve
of index form for the eigenvalue Jacobi differential problem as
the eigenvalue runs from $-\infty$ to $0$. The dependence of this
path of Fredholm forms on the eigenvalue $\lambda$ is analytic
(namely, affine), and the method of generalized Jordan chains can
be applied directly. An explicit formula for the spectral index of
a semi-Riemannian geodesic with possibly conjugate endpoint is
computed in Proposition~\ref{thm:formulaspindex}.

Motivated by infinite dimensional Morse theory, and
by the analysis of first order elliptic operators
on closed manifolds, the study of relations between  the Maslov
index and the spectral flow is a quite active research field,
and a rather extensive literature on this topic is available
nowadays.
The equality between the Maslov (or the Morse) index for paths
of self-adjoint Fredholm operators has been proven in several contexts.
Starting from the celebrated Morse index theorem in Riemannian
geometry, the Maslov index has been used by Duistermaat \cite{Duis}
to prove an index theorem for convex Hamiltonian systems.
Several versions of the Morse index theorem were proven in the last decade
in the case of geodesics in manifolds endowed with non positive
definite metrics and for non convex Hamiltonian systems
(see for instance \cite{BEE, asian, pacific, london, topology}).
Basic references on the notion of Maslov index, from which the authors
of the present paper have taken inspiration, are the articles
by  Robbin and  Salamon \cite{RobSal}, by de Gosson  \cite{deGos4},
and that of Cappell, Lee and Miller \cite{CapLeeMil}; analogies and
differences between the  results of \cite{CapLeeMil}, \cite{RobSal}
and  those of the present paper will be discussed at the beginning of
Section~\ref{sec:maslovindex}.
As to the notion of spectral flow for a path of self-adjoint Fredholm
operators, the literature available is enormous, starting from
the pioneering work of Atiyah,   Patodi and Singer \cite{APS}.
We found particularly simple and elegant the approach in a paper
by Phillips \cite{Phillips}, where functional calculus is used.
Phillips' definition of Maslov index has inspired the
proof of our Proposition~\ref{thm:Uiphii} that gives an abstract method
for constructing arbitrary group valued index functions for
continuous paths in topological spaces. However, for the
purposes of the present paper, the most appropriate definition
of spectral flow seems to be the one given in \cite{FitzPejsaRecht},
that uses the theory of relative dimension in Hilbert spaces
and the cogredient action of the general linear group;
this approach is particularly useful when dealing with
index and coindex of symmetric bilinear forms (see Proposition~\ref{thm:indrelindcoind}).
The notion of Maslov index has been related to the conjugate
points along a Riemannian geodesic by Morvan in \cite{morvan}, while
in the semi-Riemannian (i.e., non positive definite) case, Helfer
in \cite{Hel1} was the first to introduce the notion of Maslov index
of a geodesic.
The equality of the spectral index of a path $D(t)$ of Dirac operators
on a vector bundle over an odd dimensional Riemannian manifold $M$
split along an hypersurface $\Sigma$ and the
Maslov index of the pair of curves of Lagrangian
spaces $(L_1(t),L_2(t))$ obtained as the Cauchy data spaces of
$D(t)$ has been proven in \cite{Nico} by Nicolaescu
 in the case of nondegenerate endpoints, extending a previous result
of Yoshida \cite{Yoshi} in dimension  $3$; Yoshida--Nicolaescu's
theorem has been extended by Daniel in \cite{Dan}  to the possibly
degenerate case. A similar result is proven by Cappell, Lee and
Miller in
 \cite{CapLeeMil3}, where the authors obtain several formulas
expressing the spectral flow of a one-parameter family of self-adjoint
elliptic operators on a closed manifold as a sum of terms
computed from a decomposition of the manifold into two submanifolds.
In the infinite dimensional case, the relationship
between spectral flow and Maslov index for curves in the
Fredholm Lagrangian Grassmannian has been studied by
Booss--Bavnbek, Furutani and Otsuki in
\cite{BoosFuru, BoosFurub, FurOts1, FurOts2}, where index theorems
are proven under very general circumstances.

Our main interest in studying connections between
the Maslov index and the spectral flow comes from recent developments
of infinite dimensional Morse theory (see \cite{abbomajer}),
which follows the lines of the celebrated works of
Floer \cite{Flo1, Flo2, Flo3}.
The index theorem proven by Robbin and Salamon
in \cite{RobSal2} gives an equality between
the spectral flow of a one-parameter family $A(t)$ of unbounded
self-adjoint operators on a Hilbert space $H$ with the Fredholm
index of the densely defined operator $\frac{\mathrm d}{\mathrm dt}+A(t)
:L^2(\R,H)\to L^2(\R,H)$. Such result aims, as in our case,
to  Morse homology, where $A(t)$ is the Hessian
of a Morse function along an orbit of the gradient flow between two
critical points. In this case, the spectral flow of $A(t)$
gives the dimension of the intersection between the stable and
the unstable manifold at the endpoints; the finiteness of such
dimension is a central point in the construction of the
Morse--Witten complex. A crucial assumption
in Robbin and Salamon's proof is that the limits $\lim\limits_{t\to\pm\infty}
A(t)=A^\pm$ be hyperbolic, hence invertible. Assuming invertibility at the
endpoints is an obstruction for developing  Morse homology
in case of degenerate critical points, like for instance
in the case of closed  (semi-)Riemannian geodesics, in which case all
critical points are degenerate due to the equivariant action
of the group $\mathrm O(2)$. Our degenerate index theorem
aims at providing the tools for treating infinite dimensional
Morse theory in presence of equivariant, or otherwise
well-behaved, degenerate critical points.

It is very likely that one can obtain a somewhat simpler
proof of many of the index theorems mentioned above using
the theory of partial signature discussed in this paper
along the lines of the proof of Theorem~\ref{thm:tuttouguale}:
\begin{itemize}
\item  a direct (finite dimensional) homotopy argument proves
the equality between the Maslov index of the path of Lagrangians
associated to the path of self-adjoint operators
and the Maslov index of a real-analytic curve obtained
from the corresponding eigenvalue problem;
\item an infinite dimensional homotopy argument proves the
equality between the spectral flow and the ``spectral index''
of the data, which is the spectral flow of a real-analytic
curve depending on the spectral parameter;
\item the equality between the index of the two real-analytic curves
is obtained by proving the equality of the partial signatures
at each degeneracy instant.
\end{itemize}
Typically, the spectral index is the spectral flow of a path
depending affinely on the parameter, in which case the computation
of the partial signatures can be done explicitly
(see Subsection~\ref{sub:spflaffine}).

\smallskip

The paper is ideally split into two distinct parts.
Sections~\ref{sec:gensig} and \ref{sec:maslovindex} deal with the theory of
partial signatures and its applications to the
abstract index theory. Several notions of Maslov index
appearing in the literature are discussed in our framework
in the spirit of the article \cite{CapLeeMil},
including the Maslov index of pairs of Lagrangian paths
and for symplectic paths (Conley--Zehnder index),
the triple index (Kashiwara's index) and the four-fold index (H\"ormander index).
The notions of triple and four-fold index are given in terms
of the Maslov index for paths; it is an interesting observation
that, conversely, the Maslov index can be constructed
using only the axioms of Kashiwara's index (formulas~\eqref{eq:muq} and
\eqref{eq:mubarq}).
The second part of the paper,
contained in Section~\ref{sec:geodesics}, deals with the geometric
applications of the theory in the context of semi-Riemannian
geodesics. The central results are the Index Theorem~\ref{thm:tuttouguale}
and its geometrical version Theorem~\ref{thm:geomversion},
that generalize \cite[Theorem~5.2]{topology} and \cite[Theorem~3.3]{CPPT} to the case of possibly
conjugate endpoints.

Two appendices have been added at the end of the paper.
Appendix~\ref{sec:fredholmforms} contains the basic definitions
and some technical results concerning functional analytical
aspects of the index theory, with special emphasis on the
notion of commensurability of closed subspaces of Hilbert spaces, relative
dimension and relative index of Fredholm bilinear forms.
Appendix~\ref{sec:grasconn} contains the proof of the connectedness
of the Grassmannian of maximal negative subspaces of a finite dimensional
vector space endowed with a nondegenerate symmetric bilinear form;
such result is used in a homotopy argument employed in the proof
of a geometrical version of the index theorem (Subsection~\ref{sub:geomversion}).

In order to facilitate the reading, each section of the paper has been
divided into small subsections that should help the reader to keep
track of the several notions introduced and to localize cross references.

\end{section}


\begin{section}{Generalized signatures at an isolated
singularity}\label{sec:gensig}
Let $V$ be a real vector space; we will denote by
$\Bsym(V)$  be the vector space of all symmetric bilinear forms
$B:V\times V\to\R$  on $V$. When $V$ is endowed with a positive definite
inner product $\langle\cdot,\cdot\rangle$, we will
denote by $\Lsym(V)$   the vector
space of all linear maps $T:V\to V$ that are symmetric relatively to
$\langle\cdot,\cdot\rangle$. There is
an identification $\Bsym(V)\cong\Lsym(V)$ via the map $T\mapsto
B=\langle T\cdot,\cdot\rangle$, and such identification will
be made implicitly in many parts of the paper, although in some occasion
(especially when the choice of a fixed inner product is not done explicitly)
it will be convenient to maintain a conceptual distinction between
linear operators and bilinear forms.
For $B\in\Bsym(V)$, denote by $n^-(B)$,
$n_0(B)$ and $n^+(B)$ respectively the {\em index},  the {\em
degeneracy\/} (or {\em nullity}) and the {\em coindex\/}  of $B$, that are respectively the
number of $-1$'s,  the number of $1$'s and the number of $0$'s in the
canonical form of $B$ as given in Sylvester's Inertia Theorem.
For the purposes of the paper, it will be interesting to
introduce the notations \[{\noo+}(B)=n^+(B)+n_0(B),\quad\text{ and }\quad{\noo-}(B)=n^-(B)+n_0(B)\]
respectively for the {\em extended coindex\/} and the {\em extended
index\/} of $B$.
The {\em  signature\/} $\sgn(B)$ of a bilinear form\footnote{%
The use of the symbol $\sigma$ for the signature of a bilinear
form will be used quite frequently throughout the paper, and for
this reason the customary notation $\sigma(T)$ for the {\em
spectrum\/} of a linear operator $T$ will be  replaced with
$\mathfrak s(T)$.} $B$ is the difference
$n^+(B)-n^-(B)={\noo+}(B)-{\noo-}(B)$; $B$ is said to be {\em
nondegenerate\/} if $n_0(B)=0$. Clearly, $n^-(B)$, and  $n^+(B)$
are respectively the number of negative and of positive
eigenvalues of the symmetric linear map $T$ corresponding to $B$,
the degeneracy $n_0(B)$ is the multiplicity of $0$ as an
eigenvalue of $T$, and $B$ is nondegenerate exactly when $T$ is an
isomorphism. With a slight abuse of terminology, we will use the
notation $n^\pm$  (and ${\noo\pm}$) for symmetric linear maps,
meaning the (extended) index and the coindex of the corresponding
element in $\Bsym(V)$. Finally, by $\Ker(B)$ we will mean the
kernel of the corresponding linear map $T$, which can be described
also as \[\Ker(B)=\big\{w\in V:B(v,w)=0\ \text{for all}\ v\in
V\big\}.\]

In next two subsections we will assume that $V$ is a fixed finite dimensional
real vector space, and for simplicity we will assume that
$\langle\cdot,\cdot\rangle$ is a fixed positive definite
inner product on  $V$. Also in the infinite dimensional case,
we will deal with spaces endowed with a fixed Hilbert space inner product.
The choice of a specific inner product
is by no means essential, and we will prove at the end that
all the notions introduced are independent of such choice
(Remark~\ref{thm:indepinnpr}).

\subsection{Root functions and partial signatures at a
degeneracy instant}\label{sub:root}
\hfill\break
Let ${\LL}:\left[t_0-\varepsilon,t_0+\varepsilon\right]\to
\Lsym(V)$ be a curve such that $t=t_0$ is an isolated
singularity for ${\LL}(t)$. We are interested in determining the
jump of the functions $n^+\big({\LL}(t)\big)$ and
$n^-\big({\LL}(t)\big)$ as $t$ passes through $t_0$;
the following elementary result is well known:
\begin{prop}\label{thm:elementary}
Assume that the restriction $B_1$ of $\langle
{\LL}'(t_0)\cdot,\cdot\rangle$ to $\Ker\big({\LL}(t_0)\big)$ is
nondegenerate. Then:
\[\begin{split}
n^+\big({\LL}(t_0+\varepsilon)\big)-n^+\big({\LL}&(t_0)\big)=n^+(B_1),\quad
n^+\big(L(t_0)\big)-n^+\big({\LL}(t_0-\varepsilon)\big)=-n^-(B_1),\\ &
n^+\big({\LL}(t_0+\varepsilon)\big)-n^+\big({\LL}(t_0-\varepsilon)\big)=
\sgn(B_1).
\end{split}\]
\end{prop}
\begin{proof}
See for instance \cite[Proposition~2.5]{asian}.
\end{proof}
Attempts to generalize the result of
Proposition~\ref{thm:elementary} by replacing the nondegeneracy
assumption with an assumption concerning higher order
derivatives of ${\LL}(t)$ at $t=t_0$ in $\Ker\big({\LL}(t_0)\big)$
fail, as the following example shows:
\begin{example}\label{exa:doesnotwork}
Consider the curves $B,\widetilde B\colon\R\to\Bsym(\R^2)$ given by:
\[B(t)=\begin{pmatrix}1&t\\t&t^3\end{pmatrix},\quad
\widetilde B(t)=\begin{pmatrix}1&t^2\\t^2&t^3\end{pmatrix};\]
we have $B(0)=\widetilde B(0)$ and
$N=\Ker\big(B(0)\big)=\Ker\big(\widetilde B(0)\big)=\{0\}\oplus\R$. Observe that
$B(t)\vert_{N\times N}=\widetilde B(t)\vert_{N\times N}$ for all $t\in\R$, so that
the Taylor expansion of $B$ coincides with that of $\widetilde B$ in $N$;
on the other hand, for $\varepsilon>0$ sufficiently small, we have:
\begin{align*}
n^+\big(B(\varepsilon)\big)-n^+\big(B(-\varepsilon)\big)&=1-1=0,\\
n^+\big(\widetilde B(\varepsilon)\big)-n^+\big(\widetilde B(-\varepsilon)\big)&=2-1=1.
\end{align*}
\end{example}
We will now discuss a method for computing the jump of the coindex
involving higher order derivatives; let us recall that, given
$k\ge1$, a smooth  map
$v:\left]t_0-\varepsilon,t_0+\varepsilon\right[\to V$ is said to
have a {\em zero of order  $k$\/} at $t=t_0$ if
$v(t_0)=v'(t_0)=\ldots=v^{(k-1)}(t_0)=0$ and $v^{(k)}(t_0)\ne0$.
In order to set up properly our framework we give the following:
\begin{defin}
\label{thm:defroot}
A {\em root function\/} for ${\LL}(t)$ at $t=t_0$ is
a smooth map $u:[t_0-\varepsilon,t_0+\varepsilon]\to V$,
$\varepsilon>0$, such that $u(t_0)\in\Ker\big({\LL}(t_0)\big)$.
The {\em order\/} $\ord(u)$ of the root function $u$ is the
(possibly infinite)
order of zero at $t=t_0$ of the map $t\mapsto {\LL}(t)u(t)$.
\end{defin}
Associated to the degeneracy instant $t_0$ for the curve $\LL$, we
will now define a filtration of vector spaces $W_k\subset V$ and a
sequence of  bilinear forms $B_{k}:W_{k}\times W_{k}\to\R$,
$k\ge1$, as follows. Set:
\begin{equation}\label{eq:Wk}
W_k:=\big\{u_0\in V:\exists\ \ \text{a root
function $u$ with  $\ord(u)\ge k$ and $u(t_0)=u_0$}
\big\},
\end{equation}
and
\begin{equation}\label{eq:Bk}
B_k(u_0,v_0):=\frac1{k!}\big\langle
\tdtoo{k}\big[{\LL}(t)u(t)\big],v_0\big\rangle,\quad
u_0,v_0\in W_k,
\end{equation}
where $u$ is any root function with $\ord(u)\ge k$ and
$u(t_0)=u_0$. We will prove next that the right
hand side of the equality in \eqref{eq:Bk} does not indeed
depend on the choice of the root function $u$.
\begin{prop}\label{thm:proprbas}
$W_k$ is a subspace of $\Ker\big({\LL}(t_0)\big)$, and
$W_{k+1}\subseteq W_{k}$ for all $k\ge1$.
The bilinear forms $B_k$ are well defined and symmetric
for all $k\ge1$.
\end{prop}
\begin{proof}
The first statement follows trivially from the definition
of $W_k$.

For the second statement,
observe that for $u_0,v_0\in W_k$, if $u(t)$ and
$v(t)$ are root functions with $\ord(u),\ord(v)\ge
k$, $u(t_0)=u_0$, $v(t_0)=v_0$, keeping in mind that ${\LL}(t)u(t)$
and ${\LL}(t)v(t)$ have vanishing derivatives  of
order less than or equal to $k-1$ at $t=t_0$, then:
\begin{equation}\label{eq:WkBk}
\begin{split}
&\big\langle\tdtoo{k}\big[{\LL}(t)u(t)\big],v_0\big\rangle
=\tdtoo{k}\left\langle {\LL}(t)u(t),
v(t)\right\rangle\\&\qquad\qquad=
\tdtoo{k}\left\langle u(t), {\LL}(t) v(t)\right\rangle =
\big\langle u_0,\tdtoo{k}\big[{\LL}(t)v(t)\big]\big\rangle.
\end{split}
\end{equation}
The first term in \eqref{eq:WkBk} does not depend
on the choice of the root function $v(t)$ and the last
term does not depend on the choice of a root function $u(t)$.
This proves at the same time that $B_k$ is well defined and that
it is symmetric.
\end{proof}

\begin{rem}\label{thm:remuseless}
Here is an alternative definition of the
spaces $W_k$ and of the bilinear forms $B_k$.
Consider the Taylor expansion of ${\LL}(t)$ centered at $t=t_0$:
\[{\LL}_0+(t-t_0){\LL}_1+(t-t_0)^2{\LL}_2+\ldots+(t-t_0)^k{\LL}_k+\ldots,\]
where ${\LL}_k=\tfrac1{k!}{\LL}^{(k)}(t_0)\in \Lsym(V)$ for all $k$.
Then, the space $W_{k+1}$ can be described as the
set of those $u_0\in V$ such that there exist
$u_1,\ldots,u_k\in V$ for which the following system
of linear equations is satisfied:
\begin{equation}\label{eq:systgenJordch}
\left\{\begin{array}{l}{\LL}_0u_0=0,\\
{\LL}_1u_0+{\LL}_0u_1=0,\\{\LL}_2u_0+{\LL}_1u_1+{\LL}_0u_2=0,\\
\ \ \vdots\\
\sum_{j=0}^k{\LL}_{k-j}u_j=0.\end{array}\right.
\end{equation}
A root function $u(t)$ with $\ord(u)\ge k+1$ and $u(t_0)=u_0$
would be given in this case by $u(t)=\sum_{j=0}^k(t-t_0)^ju_j$; in
the language of \cite{rabier}, a sequence $(u_0,u_1,\ldots,u_k)$
satisfying \eqref{eq:systgenJordch} is called a {\em generalized
Jordan chain\/} for ${\LL}(t)$ at $t=t_0$ starting at $u_0$. The
{\em length\/} of a generalized Jordan chain is defined to be the
number of its elements; a generalized Jordan chain
$(u_0,\ldots,u_k)$ is said to be {\em extendible\/} if there
exists $u_{k+1}\in V$ such that $(u_0,\ldots,u_k,u_{k+1})$ is a
generalized Jordan chain. Thus, $W_{k+1}$ is the space of those
$u_0$ for which there exists a generalized Jordan chain of length
$k+1$ starting at $u_0$.

Non extendible Jordan chains are said to be {\em
maximal}; observe that maximal generalized Jordan chains starting
at the same element $u_0$ may have different lengths.%
\footnote{For instance, suppose that $(u_0,u_1,u_2)$ is
a generalized Jordan chain starting at $u_0$ and that
$\xi\in\Ker({\LL}_0)$ is such that ${\LL}_1\xi\not\in\Imm({\LL}_0)$.
Then $(u_0,u_1+\xi)$ is maximal.}

The first equation of \eqref{eq:systgenJordch}
tells us that $u_0\in\Ker({\LL}_0)$;
it is immediate to observe that, if $(u_0,\ldots,u_k)$
is a generalized Jordan chain for ${\LL}$ starting at $u_0$,
then it is extendible if and only if:
\[\sum_{j=0}^{k}{\LL}_{k+1-j}u_j\in\Imm\big({\LL}_0\big)=\Ker({\LL}_0)^\perp.\]
Moreover, the bilinear form $B_{k+1}$ can be defined
in terms of generalized Jordan chains by:
\begin{equation}\label{eq:descrBk+1}
B_{k+1}(u_0,v_0)=\sum_{j=0}^k
\big\langle {\LL}_{k+1-j}u_j,v_0\big\rangle,\end{equation}
where $(u_0,u_1,\ldots,u_k)$ is {\em any\/} generalized
Jordan chain for ${\LL}$ of length $k+1$ starting at $u_0$.
System \eqref{eq:systgenJordch} and formula \eqref{eq:descrBk+1}
appear in reference~\cite{FarLev} (see \cite[(23i)]{FarLev} and \cite[Corollary~3.13]{FarLev}).
\end{rem}

From the very definition of $B_k$, one sees immediately
that we have an inclusion
$W_{k+1}\subset
\Ker(B_k)$;
observe in particular that if $B_k$ is nondegenerate for some
$k$, then $W_j=(0)$ for all $j>k$.
We will show below that, in fact, $W_{k+1}=
\Ker(B_k)$ for all $k\ge1$ (Corollary~\ref{thm:Wk1kerBk}).

\begin{defin}\label{thm:defpartsign}
For all $k\ge1$, the integer numbers
\begin{equation}\label{eq:defsgnk}
n^-_k({\LL},t_0):=n^-(B_k),\quad n^+_k({\LL},t_0):=n^+(B_k),
\quad \sigma_k({\LL},t_0):=\sgn(B_k)
\end{equation}
are called respectively the {\em $k$-th partial index}, the {\em
$k$-th partial coindex}, and the {\em $k$-th partial signature\/}
of ${\LL}(t)$ at $t=t_0$. The {\em $k$-th partial extended
index\/} and {\em coindex\/} ${\noo\pm}_k({\LL},t_0)$ are defined
similarly.
\end{defin}

The integers $n^\pm_k(\LL,t_0)$, ${\noo\pm}_k(\LL,t_0)$ and
$\sigma_k(\LL,t_0)$ will be referred to collectively
as the ``partial signatures'' of the curve of bilinear forms $\LL$ at the
degeneracy instant $t_0$.

\begin{rem}\label{thm:W1B1}
It is immediate from the definition that
$W_1=\Ker({\LL}_0)$ and that $B_1$ coincides with the
restriction of $\langle {\LL}_1\,\cdot,\cdot\rangle$
(recall Proposition~\ref{thm:elementary}).
\end{rem}
\begin{example}
Let us compute explicitly the spaces $W_k$, the bilinear
forms $B_k$ and the partial signatures $\sigma_k$ for the
curves in Example~\ref{exa:doesnotwork} at the instant $t=0$.

For ${\LL}(t)=\begin{pmatrix}1&t\\t&t^3\end{pmatrix}$, one computes
easily:
\[W_1=\{0\}\oplus\R,\quad {\LL}_0=\begin{pmatrix}1&0\\
0&0\end{pmatrix},\quad {\LL}_1=
\begin{pmatrix}0&1\\ 1&0\end{pmatrix},
\quad {\LL}_2\begin{pmatrix}0&0\\ 0&0\end{pmatrix}.\]
For all $u_0=(0,a)$, the system ${\LL}_1u_0+{\LL}_0u_1$ has the unique
solution $u_1=(-a,0)$, i.e., $W_2=W_1=\{0\}\oplus\R$; moreover, it
is easily seen that the generalized Jordan chain
$\big((0,a),(-a,0)\big)$ is maximal, i.e., the system
${\LL}_2u_0+{\LL}_1u_1+{\LL}_0u_2$ has no solution $u_2\in\R^2$ and
therefore $W_3=(0)$.

Using \eqref{eq:descrBk+1} we compute:
\[\begin{split}
&B_1\big((0,a),(0,b)\big)=\big\langle
{\LL}_1u_0,v_0\big\rangle=0,
\quad\forall \,a,b\in\R,\quad n_1^\pm=0,\ \sigma_1=0\\
&B_2\big((0,a),(0,b)\big)=\big\langle
{\LL}_2u_0+{\LL}_1u_1,v_0\big\rangle=-ab,\quad n^+_2=0,\ n_2^-=1,\ \sigma_2=-1.
\end{split}\]
Clearly, $W_k=(0)$, $B_k=0$ and $\sigma_k=0$ for all $k\ge3$.

A similar computation for
${\LL}(t)=\begin{pmatrix}1&t^2\\t^2&t^3\end{pmatrix}$ gives:
\[W_1=\{0\}\oplus\R,\quad {\LL}_0=\begin{pmatrix} 1&0 \\ 0 &0
\end{pmatrix}, \quad {\LL}_1=\begin{pmatrix} 0&0 \\  0&0
\end{pmatrix},\quad {\LL}_2=\begin{pmatrix} 0&1 \\  1&0
\end{pmatrix}, \quad {\LL}_3\!=\!\begin{pmatrix} 0&0 \\  0&1
\end{pmatrix}\!.\]
Set $u_0=(0,a)\in W_1$; the system ${\LL}_1u_0+{\LL}_0u_1$ has
solutions $u_1=(0,\beta)$, with $\beta\in\R$; hence, $W_2=W_1$. The system
${\LL}_2u_0+{\LL}_1u_1+{\LL}_0u_2$ has the unique solution $u_2=(-a,0)$, and so $W_3=W_2$.
Finally, the generalized Jordan chain $(u_0,u_1,u_2)$
is maximal, because the system ${\LL}_3u_0+{\LL}_2u_1+{\LL}_1u_2+{\LL}_0u_3=0$
has no solution, i.e., $W_k=(0)$ for all $k\ge4$.

Using \eqref{eq:descrBk+1} we compute:
\[\begin{split}
&B_1\big((0,a),(0,b)\big)=\big\langle
{\LL}_1u_0,v_0\big\rangle=0, \quad n^\pm_1=0,\ \sigma_1=0,\\
& B_2\big((0,a),(0,b)\big)=\big\langle {\LL}_2u_0+{\LL}_1u_1,v_0
\big\rangle=\langle (a,0),(0,b)\rangle=0,\quad n_2^\pm=0,\ \sigma_2=0,\\
&
B_3\big((0,a),(0,b)\big)=\big\langle
{\LL}_3u_0+{\LL}_2u_1+{\LL}_1u_2,v_0\big\rangle=\big\langle(\beta,a),(0,b)\big\rangle=ab,\\
&\qquad\qquad\qquad\qquad\qquad\qquad\qquad\qquad\qquad\qquad\qquad\qquad n_3^+=1,\ n_3^-=0,\ \sigma_3=1.
\end{split}\]
\end{example}

\subsection{Computation of the partial signatures and the
spectral flow}\hfill\break When $V$ has a smoothly varying basis
of eigenvectors of ${\LL}(t)$, then the computation of the spaces
$W_k$ and of the bilinear forms $B_k$ can be simplified as
explained in the following:
\begin{prop}\label{thm:central}
Let ${\LL}:[t_0-\varepsilon,t_0+\varepsilon]\to\Lsym(V)$ be a
smooth curve having a unique degeneracy instant at $t=t_0$. Assume
that the following regularity condition for the eigensystem of
${\LL}(t)$ is satisfied: denoting by
$\lambda_1(t),\ldots,\lambda_n(t)$ the smooth functions of
eigenvalues of ${\LL}(t)$,
\begin{itemize}

\item[(a)] each non constant $\lambda_i$ has a zero of finite
order at $t=t_0$;

\item[(b)] there exist {\em smooth\/} functions
$\vv_i:[t_0-\varepsilon,t_0+\varepsilon]\to V$, $i=1,\ldots,n$, of
pairwise orthogonal unit eigenvectors of the $\lambda_i(t)$'s.
\end{itemize}
Then, the following statements hold:
\begin{enumerate}
\item\label{itm:span}
$W_k=\textrm{span}\Big\{\vv_i(t_0):\text{$i\in\{1,\ldots,n\}$ is
such that
$\lambda_i^{(j)}(t_0)=0$ for all $j<k$}\Big\}$;
\item\label{itm:eigenv} if $v\in W_k$ is an eigenvector of
$\lambda_i(t_0)$, where
$\lambda_i(t_0)\!=\!\lambda_i'(t_0)\!=\!\ldots\!=\!\lambda_i^{(k-1)}(t_0)\!=\!0$,
then $B_k(v,w)=\tfrac1{k!}\lambda_i^{(k)}(t_0)\langle
v,w\rangle$, for all $w\in W_k$;
\smallskip

\item\label{itm:formprinc}
$n^+\big({\LL}(t_0+\varepsilon)\big)-n^+\big({\LL}(t_0)\big)=
\sum\limits_{k\ge1}n^+_k(\LL,t_0)$,
\\ $n^+\big({\LL}(t_0)\big)-n^+\big({\LL}(t_0-\varepsilon)\big)=
-\sum\limits_{k\ge1}\big(n^-_{2k-1}(\LL,t_0)+n^+_{2k}(\LL,t_0)\big)$, and\\
$n^+\big({\LL}(t_0+\varepsilon)\big)-n^+\big({\LL}(t_0-\varepsilon)\big)=
\sum\limits_{k\ge 1}\sigma_{2k-1}(\LL,t_0)$,\\ where all the sums appearing in these
formulas have  at most a finite number of non zero terms.
 Similarly:\footnote{See also formulas~\eqref{eq:nuoveformule} below.}\\ \\
${\noo+}\big({\LL}(t_0+\varepsilon)\big)-
{\noo+}\big({\LL}(t_0)\big)=\sum\limits_{k\ge1}n^+_k(\LL,t_0)-
\Dim\big(\Ker(\LL(t_0))\big)$,\\
${\noo+}\big({\LL}(t_0)\big)-{\noo+}\big({\LL}(t_0-\varepsilon)\big)\!\!=\!
-\!\!\sum\limits_{k\ge1}\!\!\big(n^-_{2k-1}(\LL,\!t_0)+n^+_{2k}(\LL,\!t_0)\big)
+\Dim\big(\Ker(\LL(t_0))\big)$.
\end{enumerate}
\end{prop}
\begin{proof}
To prove \eqref{itm:span} argue as
follows. If $i\in\{1,\ldots,n\}$ is such that
$\lambda_i(t_0)=\lambda_i'(t_0)=\ldots=\lambda_i^{(k-1)}(t_0)=0$,
then setting $u_0=\vv_i(t_0)$ and $u(t)=\vv_i(t)$ it
follows immediately that $u$ is a root function
of order greater than or equal to $k$ at $t=t_0$ with $u(t_0)=u_0$.
This proves that $\vv_i(t_0)\in W_k$, i.e., that the span of such
$\vv_i(t_0)$'s is contained in $W_k$.
On the other hand, assume that $u_0\in W_k$; let $u(t)$ be
a root function with $\ord(u)\ge k$, $u(t_0)=u_0$, and
set $u(t)=\sum_i\mu_i(t)\vv_i(t)$. Then:
\[{\LL}(t)u(t)=\sum_i\mu_i(t)\lambda_i(t)\vv_i(t);\]
from the above equality it follows easily that
${\LL}(t)u(t)$ has a zero of order greater than or equal to
$k$ at $t=t_0$ if and only if the function $(\lambda_i\mu_i)$
has
a zero of order greater than or equal to $k$ at $t=t_0$
for all $i$.
In  particular, this implies that
$\mu_i(t_0)=0$  unless
$\lambda_i(t_0)=\lambda_i'(t_0)=\ldots=\lambda_i^{(k-1)}(t_0)=0$,
which proves \eqref{itm:span}.

The proof of part \eqref{itm:eigenv} is immediate using
the definition of $B_k$, taking a proper multiple
of $\vv_i$ as a root function for $v$.

In view of assumption~(a) in the hypotheses, part
\eqref{itm:span} and
\eqref{itm:eigenv} of the thesis, for each $k\ge1$, the index
$n^-(B_k)$  (resp., the coindex $n^+(B_k)$)
is given by the number of $i$'s in $\{1,\ldots,n\}$
such that $\lambda_i(t)$ has a zero of order $k$ at $t=t_0$
and whose $k$-th derivative is negative (resp., positive)
at $t=t_0$. The five formulas in part~\eqref{itm:formprinc} of
the thesis follow easily from this observation; note in
particular that the third formula is obtained by addition of the
first two.
\end{proof}
We observe that the assumptions of Proposition~\ref{thm:central}
are satisfied  when $t\mapsto {\LL}(t)$ is
 a real-analytic map with an
isolated degeneracy at
$t=t_0$. Namely, in this case both the eigenvalues and the
eigenvectors are real-analytic functions (see \cite[Chapter~2, \S~1]{kato}),
which implies immediately assumptions (a) and (b) in the
statement. Part \eqref{itm:span} and \eqref{itm:eigenv} of
Proposition~\ref{thm:central} are used as definition of the spaces
$W_k$ and the bilinear forms $B_k$ in \cite{KirKla}. Another
situation in which Proposition~\ref{thm:central} can be applied is
when $\Dim\big(\Ker(L_0)\big)=1$, in which case obviously the
entire statement and the proof can be simplified; the case of {\em
simple\/} eigenvalues can be treated in a much easier  way by
means of the Implicit Function theorem (see for instance
\cite{BGM, MasPic}), and the occurrence of such circumstance will
not be further commented in this paper.

Observe also that part \eqref{itm:formprinc} of the thesis
of Proposition~\ref{thm:central} is a generalization of
the result of Proposition~\ref{thm:elementary}.

\begin{cor}\label{thm:Wk1kerBk}
Let ${\LL}:[t_0-\varepsilon,t_0+\varepsilon]\to\Lsym(V)$ be any smooth curve
having a unique degeneracy instant at $t=t_0$. Then $W_{k+1}=\Ker(B_k)$
for all $k\ge1$.
\end{cor}
\begin{proof}
If ${\LL}$ is real analytic, the conclusion follows
immediately from
part~\eqref{itm:span} and part~\eqref{itm:eigenv}
of Proposition~\ref{thm:central}.

Assume that ${\LL}$ is smooth;
recall from Remark~\ref{thm:remuseless} that, for each $k\ge1$,
the definition of the space $W_k$ and of the bilinear form
$B_k$ depends only on ${\LL}_0$ and the first $k$ derivatives
${\LL}_1,\ldots,{\LL}_k$ of ${\LL}(t)$ at $t=t_0$.
In particular, the object $W_j$ and
$B_j$, $j=1,\ldots,k$
do not change if we replace ${\LL}$ by its $k$-th order Taylor
polynomial
$\widetilde {\LL}(t)=\sum_{j=0}^k{\LL}_j(t-t_0)^j$, which is a real-analytic
map and the first part of the proof applies.
\end{proof}

\begin{cor}\label{thm:sommank}
Under the assumptions of Proposition~\ref{thm:central},
\begin{equation}
\label{eq:sommank}
\sum_{k\ge1}\Big[n^+_k\big({\LL};t_0\big)+n^-_k\big({\LL};t_0\big)\Big]=\Dim\big(\Ker(\LL(t_0))\big).
\end{equation}
\end{cor}
\begin{proof}
By Corollary~\ref{thm:Wk1kerBk}, for all $k\ge1$:
\[n^+_k\big({\LL};t_0\big)+n^-_k\big({\LL};t_0\big)=\Dim(W_k)-\Dim(W_{k+1}).\]
The conclusion follows from an easy induction argument, keeping in
mind that, under the assumption of Proposition~\ref{thm:central},
$W_{k+1}=\{0\}$ for $k$ sufficiently large.
\end{proof}
Using \eqref{eq:sommank}, the last two formulas in
part~\eqref{itm:formprinc} of Proposition~\ref{thm:central}
can be rewritten as:
\begin{equation}\label{eq:nuoveformule}
\begin{split}
& {\noo+}\big({\LL}(t_0+\varepsilon)\big)-
{\noo+}\big({\LL}(t_0)\big)=-\sum_{k\ge1}n_k^-(\LL,t_0),\\
& {\noo+}\big({\LL}(t_0)\big)-{\noo+}\big({\LL}(t_0-\varepsilon)\big)=
\sum\limits_{k\ge1} \big(n^+_{2k-1}(\LL,t_0)+n^-_{2k}(\LL,t_0)\big).
\end{split}
\end{equation}

\begin{example}
For ${\LL}(t)=\begin{pmatrix}1&t\\t&t^3\end{pmatrix}$
one computes easily:
\[\lambda_1(t)=\frac12\big(1+t^3-\sqrt{1+4t^2-2t^3+t^6}\big),
\quad
\lambda_2(t)=\frac12\big(1+t^3+\sqrt{1+4t^2-2t^3+t^6}\big),\]
and:
\[\lambda_1(0)=\lambda_1'(0)=0,\ \lambda_1''(0)=-2,\
\vv_1(0)=(0,1),\
\lambda_2(0)=1.\]
Using Proposition~\ref{thm:central} one computes immediately:
\[\begin{split}
& W_1=\Ker\big({\LL}(0)\big)=\{0\}\oplus\R,\ B_1=0,\ n^\pm_1=0,\ \sigma_1=0\\
& W_2=\{0\}\oplus\R,\ B_2\big((0,a),(0,b)\big)=-ab,\ n^+_2=0,\ n^-_2=1,\
\sigma_2=-1,\\ & W_k=(0)\ \text{and}\ n^\pm_k=\sigma_k=0\ \text{for
all
$k\ge3$}.
\end{split}
\]

Taking now
${\LL}(t)=\begin{pmatrix}1&t^2\\t^2&t^3\end{pmatrix}$,
one has:
\[\lambda_1(0)=0,\ \lambda_1'(0)=0,\ \lambda_1''(0)=0,\
\lambda_1'''(0)=6,\ \lambda_2(0)=1, \vv_1(0)=(0,1).\]
Hence:
\[\begin{split}
&W_1=\{0\}\oplus\R, \ B_1=0,\ n^\pm_1=0,\ \sigma_1=0,
\\
&W_2=\{0\}\oplus\R, \ B_2=0,\ n^\pm_2=0,\ \sigma_2=0,
\\
&W_3=\{0\}\oplus\R, \ B_3=\big((0,a),(0,b)\big)=ab,\ n^+_3=1,\ n^-_3=0,\
\sigma_3=1,\\
&W_k=(0)\ \text{and}\ n^\pm_k=\sigma_k=0\ \text{for all}\ k\ge4.
\end{split}\]
\end{example}

Given a curve $[a,b]\ni t\mapsto B(t)$ of symmetric bilinear forms
on a finite dimensional vector space, the difference
\[\noo+\big(B(b)\big)-\noo+\big(B(a)\big)=n^-\big(B(a)\big)-n^-\big(B(b)\big)\]
is called the {\em spectral flow of the curve $B$ on $[a,b]$}, and
it will be denoted by $\spfl\big(B,[a,b]\big)$.
Such definition will be extended to the infinite dimensional
case in next subsection.

\subsection{Spectral flow and relative dimension in Hilbert spaces}
The result of Proposition~\ref{thm:central}
can be extended easily to analytic
paths of self-adjoint Fredholm operators, after proper
rephrasing of the statement.

Let $\Hcal$ denote a separable real Hilbert space, by $\Lcalsa$
the space of self-adjoint bounded linear operators on $\Hcal$, and
by $\Fcalsa(\Hcal)$ the open subset of $\Lcalsa(\Hcal)$ of
Fredholm operators. If ${\LL}:[a,b]\to\Fcalsa(\Hcal)$ is a
continuous curve, then an integer number is naturally associated
to ${\LL}$, called the {\em spectral flow\/} of ${\LL}$ over
$[a,b]$. Such number, denoted by:
\[\spfl\big({\LL},[a,b]\big)\]
is roughly speaking  the integer given by the number  of
negative eigenvalues of $\LL(a)$
that become nonnegative as the parameter $t$  goes from $a$ to $b$
minus the number of nonnegative eigenvalues of $\LL(a)$ that   become
negative.
Observe that for paths of strongly indefinite
self-adjoint operators\footnote{i.e., whose
essential spectrum is   contained neither in $\R^+$ nor
in $\R^-$.} both the
index and the coindex functions are infinite; we refer to \cite{Phillips} for a concise
introduction to the spectral flow for a continuous path
of self-adjoint Fredholm operators, although the reader
will find that in the literature
it is most frequently treated only the
case of paths having invertible endpoints.
We observe here that for paths with degenerate endpoints
there are several options for the choice of a definition of
spectral flow, depending on how one wants to consider the
contribution given by the kernel at the endpoints.
For the sake of consistency with our definition of Maslov
index in the finite dimensional case (Corollary~\ref{thm:esunmaslovindex}), we will
use the infinite dimensional analogue of the
variation of the {\em extended coindex\/} $\noo+$.
Let us recall a formula proven in \cite{FitzPejsaRecht}
that gives the spectral flow in terms of relative
dimension of closed subspaces\footnote{%
Two closed subspace $V,W\subset\Hcal$
 are {\em commensurable\/} if $P_V\vert_W:W\to V$ is a
Fredholm operator; if $V$ and $W$ are commensurable the {\em relative dimension\/}
$\mathrm {dim}_V(W)$ of $W$ with respect to $V$ is defined as:
\[\mathrm{dim}_V(W)=\mathrm{dim}(W\cap V^\perp)-\mathrm{dim}(W^\perp\cap V).\]
Some basic facts of the theory of relative dimension of closed
subspaces of Hilbert spaces will be recalled in
Appendix~\ref{sec:fredholmforms}.} of a Hilbert space. If
${\LL}(t)$ is of the form $\mathfrak J+K(t)$, where $\mathfrak J$
is a self-adjoint {\em symmetry\/} of $\Hcal$ (i.e., $\mathfrak
J^2=\mathrm{Id}$)
and $K(t)$ is a compact self-adjoint operator on $\Hcal$, then
\begin{equation}\label{eq:spflreldim}
\begin{split}
\spfl\big({\LL},&[a,b]\big)=\\  &\Dim_{V^{\scriptscriptstyle
-}(\mathfrak J+K(b))}\big(V^{\scriptscriptstyle -}(\mathfrak
J+K(a))\big)=\\ & \Dim_{V^{\scriptscriptstyle +}(\mathfrak
J+K(a))}\big(V^{\scriptscriptstyle +}(\mathfrak J+K(b))\big)
\!+\!\Dim\big(\Ker(\mathfrak J+K(b))\big)
\!-\!\Dim\big(\Ker(\mathfrak
J+K(a))\big),\end{split}\end{equation} where
$V^{\scriptscriptstyle -}(S)$ and $V^{\scriptscriptstyle +}(S)$
denote respectively  the {\em negative\/} and the {\em positive
eigenspace}\footnote{The negative (resp., the positive) eigenspace
of a self-adjoint operator $S$ can be defined, for instance, using
functional calculus as $\chi_{\left]-\infty,0\right[}\big(S\big)$
(resp., $\chi_{\left]0,+\infty\right[}\big(S\big)$), where
$\chi_I$ is the characteristic function.}
 of the   operator $S$. The second equality in \eqref{eq:spflreldim}
is proven easily using the result of Proposition~\ref{thMlsommadim}.
The computation of the spectral flow of an arbitrary continuous
path of self-adjoint Fredholm operators is then reduced to the
above case using the cogredient action of the general linear group
of $\Hcal$ (see Subsection~\ref{sub:invariance}).
\begin{rem}\label{thm:defspectralMorse}
Observe that if ${\LL}(a)$ (hence $\LL(t)$ for all $t$)
is an essentially positive operator, i.e.,
its essential spectrum is contained in $\left[0,+\infty\right[$, then
\[\spfl\big(\LL,[a,b]\big)=n^-\big(\LL(a)\big)-n^-\big(\LL(b)\big).\]
Likewise, if $\LL(a)$ is essentially negative, then
\[\spfl\big(\LL,[a,b]\big)=\noo+\big(\LL(b)\big)-\noo+\big(\LL(a)\big).\]
\end{rem}
The jumps of the spectral flow of a path occur precisely at the
degeneracy instants, and, using a Galerkin approximation, the
computation of the jump is reduced to a dimension counting of
finite rank projections. Observe also that, by the finite
dimensionality of the kernel, one can define the partial
signatures $n^\pm_k(\LL,t_0)$, ${\noo\pm}_k(\LL,t_0)$ and
$\sigma_k(\LL,t_0)$ at an isolated singularity $t=t_0$  of a path $\LL$ in
$\Fcalsa(\Hcal)$ exactly as in Definition~\ref{thm:defpartsign}.
Given a real-analytic path
${\LL}:[t_0-\varepsilon,t_0+\varepsilon]\to\Fcalsa(\Hcal)$ having a
unique degeneracy at $t=t_0$, then, for $t$ sufficiently close to
$t_0$ and for $a>0$ small enough, the intersection of $[-a,a]$ with
the spectrum of ${\LL}(t)$ consists of a finite number of
eigenvalues $\lambda_1(t),\ldots,\lambda_N(t)$ having bounded
multiplicity,  that are real-analytic functions of $t$ and that
vanish  (possibly) only at $t=t_0$. In this situation, the spectral
flow of ${\LL}$ over $[t_0-\varepsilon,t_0+\varepsilon]$ is computed, as
in the finite dimensional case, by looking at the change of sign
of the $\lambda_i$'s through $t_0$.   Using the same arguments
in proof of Proposition~\ref{thm:central} one obtains the
following:
\begin{cor}\label{thm:centralinfdim} Let
${\LL}:[t_0-\varepsilon,t_0+\varepsilon]\to\Fcalsa(\Hcal)$
be a real-analytic curve
of self-adjoint Fredholm operators on the real separable Hilbert space $\Hcal$
having a unique degeneracy instant at $t=t_0$. Then, the spectral flow of
${\LL}$ is computed as:
\begin{eqnarray*}&&
\begin{split}
\spfl\big({\LL},[t_0-\varepsilon,t_0]\big)=&\;
-\sum_{k\ge1}\big(n^-_{2k-1}(\LL,t_0)+n^+_{2k}(\LL,t_0)\big)+\Dim\big(\Ker(\LL(t_0))\big)
\\&=\sum_{k\ge1}\big(n^+_{2k-1}(\LL,t_0)+n^-_{2k}(\LL,t_0)\big),\end{split}\\
&&\spfl\big({\LL},[t_0,t_0+\varepsilon]\big)=\sum_{k\ge1}n^+_{k}(\LL,t_0)-\Dim\big(\Ker(\LL(t_0))\big)=-\sum_{k\ge1}n^-_{k}(\LL,t_0),
\\
&&\spfl\big({\LL},[t_0-\varepsilon,t_0+\varepsilon]\big)=\sum_{k\ge1}
\sigma_{2k-1}(\LL,t_0).\qed\end{eqnarray*}
\end{cor}

Explicit computations of the partial signatures of paths
of self-adjoint Fredholm operators will be done in Section~\ref{sec:maslovindex}
(see Proposition~\ref{thm:surpresa} and Corollary~\ref{thm:affine}).

We conclude this subsection with an elementary result concerning
the spectral flow of paths of self-adjoint Fredholm operators
having a common degeneracy instant at the right end point, and
at which the partial
signatures differ by the sign.
\begin{prop}\label{thm:cambiasegno}
Let $\Hcal_1$ and $\Hcal_2$ be real separable Hilbert spaces;
let be given two real analytic paths
$\left[t_0-\varepsilon,t_0\right]\ni t\mapsto\LL^{(1)}(t)\in\Fcalsa\big(\Hcal_1\big)$
and $\left[t_0-\varepsilon,t_0\right]\ni t\mapsto\LL^{(2)}(t)\in
\Fcalsa\big(\Hcal_2\big)$  having a unique
degeneracy instant at $t=t_0$. Assume that there exists an
isomorphism $Z:\Ker\big(\LL^{(1)}(t_0)\big)\to \Ker\big(\LL^{(2)}(t_0)\big)$
such that:
\begin{itemize}
\item[(a)] $Z\big(W_k(\LL^{(1)},t_0)\big)=W_k\big(\LL^{(2)},t_0\big)$;
\item[(b)] the pull-back $Z^*\big(B_k(\LL^{(2)},t_0)\big)$
coincides with $-B_k\big(\LL^{(1)},t_0\big)$,
\end{itemize}
for all $k\ge1$. Then,
\begin{equation}\label{eq:flussienucleo}
\spfl\big(\LL^{(1)},[t_0-\varepsilon,t_0]\big)=-
\spfl\big(\LL^{(2)},[t_0-\varepsilon,t_0]\big)+\Dim\big(\Ker(\LL^{(2)}(t_0))\big).
\end{equation}
\end{prop}
\begin{proof}
Using  Corollary~\ref{thm:centralinfdim},  the equalities $\Dim\big(\Ker(\LL^{(1)}(t_0))\big)=\Dim\big(\Ker(\LL^{(2)}(t_0))\big)$,
\[
n^-\big(B_{2k-1}(\LL^{(1)},t_0)\big)=n^+\big(B_{2k-1}(\LL^{(2)},t_0)\big),\
n^+\big(B_{2k}(\LL^{(1)},t_0)\big)=n^-\big(B_{2k}(\LL^{(2)},t_0)\big)
\]
and Corollary~\ref{thm:Wk1kerBk}, we  compute easily:
\begin{align*}
\spfl\big(\LL^{(1)},&[t_0-\varepsilon,t_0]\big)+
\spfl\big(\LL^{(2)},[t_0-\varepsilon,t_0]\big)=\\
=\;&2\,\Dim\big(\Ker(\LL^{(2)}(t_0)\big)-
\sum_{k\ge1}\big[n^+\big(B_{k}(\LL^{(2)},t_0)\big)+
n^-\big(B_{k}(\LL^{(2)},t_0)\big)\big]\\
=\;&2\,\Dim\big(\Ker(\LL^{(2)}(t_0)\big)-
\sum_{k\ge1}\Dim\big(W_k(\LL^{(2)},t_0)\big)+\sum_{k\ge1}
\Dim\big(\Ker(B_k(\LL^{(2)},t_0))\big)\\
=\;&2\,\Dim\big(\Ker(\LL^{(2)}(t_0)\big)-
\sum_{k\ge1}\Dim\big(W_k(\LL^{(2)},t_0)\big)+\sum_{k\ge1}
\Dim\big(W_{k+1}(\LL^{(2)},t_0)\big)\\
=\;&2\,\Dim\big(\Ker(\LL^{(2)}(t_0)\big)-\Dim\big(W_{1}(\LL^{(2)},t_0)\big)
=\Dim\big(\Ker(\LL^{(2)}(t_0)\big).\qedhere
\end{align*}
\end{proof}
\noindent\ \ Clearly, a similar statement holds for the spectral flow on an interval
of the form $[t_0,t_0+\varepsilon]$.
\subsection{Invariance properties of the partial  signatures}
\label{sub:invariance}
For the computation of the spectral flow of a curve
of self-adjoint Fredholm operators on a real separable
Hilbert space
$\Hcal$, it will be useful to consider the {\em cogredient action\/}
of the general linear group $\mathrm{GL}(\Hcal)$ of $\Hcal$ on
$\Fcalsa(\Hcal)$.  Recall that the cogredient action is the map:
\[\mathrm{GL}(\Hcal)\times\Fcalsa(\Hcal)\ni (S,T)\longmapsto
S^*TS\in\Fcalsa(\Hcal); \]
for instance, it is proven in \cite{FitzPejsaRecht}
that, given any curve ${\LL}:[a,b]\to\Fcalsa(\Hcal)$
of class $C^k$ and any symmetry
$\mathfrak J$ of $\mathcal H$, there exists a curve $M:[a,b]\to
\mathrm{GL}(\Hcal)$ of class  $C^k$ such that $M(t)^*{\LL}(t)M(t)-\mathfrak J$
is compact for all $t$. Recall also that the spectral flow
 is invariant by the cogredience
(see Corollary~\ref{thm:invspfl} below).

In view of these observations and of the result of
Corollary~\ref{thm:centralinfdim}, we are naturally
led to the   following:
\begin{prop}\label{thm:invcogredient}
Let ${\LL},\widetilde {\LL}:[a,b]\to\Fcalsa(\Hcal)$ be cogredient
smooth curves. If $t_0\in\left]a,b\right[$ is an isolated
degeneracy instant for the two curves, then ${\LL}$ and
$\widetilde{\LL}$ have the same partial signatures at $t=t_0$.
\end{prop}
\begin{proof}
Assume that ${\LL}(t)=M(t)^*\widetilde {\LL}(t)M(t)$, where $M:[a,b]\to
\mathrm{GL}(\Hcal)$ is a smooth curve; denote respectively by
$W_k,B_k$ and $\widetilde W_k,\widetilde B_k$ the objects
(vector space space and symmetric bilinear form)
defined in Subsection~\ref{sub:root} relatively
to the curves ${\LL}$ and $\widetilde {\LL}$ at the instant $t=t_0$.
Clearly,
$\widetilde W_1=\Ker\big(\widetilde
{\LL}(t_0)\big)=M(t_0)\big(W_1\big)$; we will prove that
$\widetilde W_k=M(t_0)\big(W_k\big)$ for all $k$, and that
$B_k$ is the pull-back of $\widetilde B_k$ by the isomorphism
$M(t_0)$.

To this aim, fix $k\ge1$, choose $u_0,v_0\in W_k$ and let $u(t)$
and $v(t)$ be root functions for ${\LL}(t)$ at $t=t_0$ with
$\ord(u),\ord(v)\ge k$, $u(t_0)=u_0$, $v(t_0)=v_0$. Set $\tilde
u_0=M(t_0)u_0$, $\tilde v_0=M(t_0)v_0$, $\tilde u(t)=M(t)u(t)$ and
$\tilde v(t)=M(t)v(t)$; in first place observe that:
\[\widetilde {\LL}(t)\tilde u(t)={M(t)^*}^{-1}{\LL}(t)u(t),\]
which implies that $\tilde u$ is a root function for $\widetilde
{\LL}(t)$ at $t=t_0$ with $\ord(\tilde u)\ge k$ and $\tilde
u(t_0)=\tilde u_0$. This says that $\widetilde
W_k=M(t_0)\big(W_k\big)$ for all $k\ge1$.

Moreover, using \eqref{eq:WkBk}, we compute easily:
\[\begin{split}
B_k(u_0,v_0)&\;=\frac1{k!}\ddto{k}\!\!\big\langle {\LL}(t)u(t),
v(t)\big\rangle=\frac1{k!} \ddto{k}\!\!\big\langle M(t)^*\widetilde
{\LL}(t) M(t)u(t),v(t)\big\rangle\\ &=
\frac1{k!}\ddto{k}\!\!\big\langle \widetilde {\LL}(t)\tilde u(t),
\tilde v(t)\big\rangle=\widetilde B_k(\tilde u_0,\tilde v_0),
\end{split}\]
which concludes the proof.
\end{proof}
\begin{rem}\label{thm:indepinnpr}
The notion of partial signatures and of spectral flow
can be given for curves of Fredholm symmetric bilinear
forms\footnote{%
A  bilinear form $B$ on a Hilbert space $\Hcal$ is said to be
Fredholm if it is realized by a Fredholm operator on $\Hcal$. Such
notion does not indeed depend on the choice of a Hilbert space
inner product on $\Hcal$.} on a Hilbert space, by considering the
associated curve of Fredholm self-adjoint operators that realize
the forms relatively to the inner product of $\Hcal$. If $\LL$ and
$\widetilde\LL$ are smooth curves in $\Fcalsa(\Hcal)$ associated
to the curve of Fredholm symmetric bilinear forms $\mathbf B$ on
$\Hcal$ using two different  Hilbert space inner products on
$\Hcal$, then $\LL$ and $\widetilde\LL$ are cogredient. It follows
that the notion of partial signatures at an isolated degeneracy
instant of a curve of Fredholm symmetric bilinear forms do not
depend on the choice of a Hilbert space inner product on $\Hcal$.
\end{rem}
Observe that, using Corollary~\ref{thm:centralinfdim},
Proposition~\ref{thm:invcogredient}, the density of
the set of real analytic curves in
$C^0\big([a,b],\Fcalsa(\Hcal)\big)$ and the homotopy invariance
of the spectral flow, one gets an alternative and simple proof
of the following (see \cite[Proposition~3.2]{FitzPejsaRecht}):
\begin{cor}\label{thm:invspfl}
The spectral flow is invariant by cogredience.\qed
\end{cor}
We will also be interested in a different invariance property of the
partial signatures (see Proposition~\ref{thm:3Lagr}):
\begin{prop}\label{thm:secinv}
Let ${\LL}:[a,b]\to\Fcalsa(\Hcal)$ be a smooth curve having an
isolated degeneracy instant at $t=t_0\in\left]a,b\right[$.
Suppose that $h:[a,b]\to\mathrm{GL}(\Hcal)$ is a smooth curve
of isomorphisms of $\Hcal$ such that:
\begin{enumerate}
\item\label{itm:a1} $h(t)^*{\LL}(t)={\LL}(t)h(t)$ for all $t\in[a,b]$;
\item\label{itm:a2} $h(t_0)$ is the identity on
$\Ker\big({\LL}(t_0)\big)$.
\end{enumerate}
Then, $t_0$ is an isolated degeneracy instant of $\widetilde{\LL}$,
and the curves ${\LL}$ and $\widetilde{\LL}:={\LL}\circ h$ have the same
partial signatures at $t=t_0$.
\end{prop}
\begin{proof}
Condition~\eqref{itm:a1} guarantees that $\widetilde {\LL}$ is a smooth curve
in $\Fcalsa(\Hcal)$. Condition~\eqref{itm:a2} tells us that $t_0$ is
an isolated degeneracy instant for $\widetilde {\LL}$, and that
\[W_1=\Ker\big({\LL}(t_0)\big)=\Ker\big(\widetilde {\LL}(t_0)\big)=
\widetilde W_1.\]
If $u_0\in W_1$ and $u(t)$ is a root function for ${\LL}$
starting at $u_0$, then $\tilde u(t)=h(t)^{-1}\circ u(t)$ is a root
function for $\widetilde {\LL}$ starting at $h(t_0)^{-1}u_0=u_0$, and $\ord(u)=\ord(\tilde u)$
because ${\LL}(t)u(t)=\widetilde {\LL}(t)\tilde u(t)$.
This proves that $W_k=\widetilde W_k$ for all $k$. Finally,
using \eqref{eq:Bk}, the equality $B_k=\widetilde B_k$ is easily checked by:
\[B_k(u_0,v_0)=\big\langle\tdto k{\LL}u,v_0\big\rangle=\big\langle\tdto
k\widetilde {\LL}\tilde u,
v_0\big\rangle=\widetilde B_k(u_0,v_0).\]
This concludes the proof.
\end{proof}

Arguing as above, we get the following invariance property
for the spectral flow:
\begin{cor}\label{thm:invSFLH}
Let ${\LL}$, $h$ and $\widetilde{\LL}$ be  as in Proposition~\ref{thm:secinv},
assuming that assumption \eqref{itm:a2} holds for all $t_0\in[a,b]$. Then,
$\spfl({\LL},[a,b])=\spfl(\widetilde {\LL},[a,b])$.
\qed\end{cor}

\subsection{Spectral flow of restrictions}
\label{sub:restrictions}
Let us conclude this chapter with a few simple observations
on the computation of the partial signatures of restrictions
of Fredholm bilinear forms.
Let us assume that ${\LL}:[-\varepsilon,\varepsilon]\to\Fcalsa(\Hcal)$ is
a real-analytic map and let us denote by $S(t)=\langle{\LL}(t)\cdot,\cdot\rangle$
the corresponding map of symmetric bilinear forms. Let $V_1,V_2\subset\Hcal$
be closed subspaces such that the restrictions $S_1(t)=S(t)\vert_{V_1}$
and $S_2(t)=S(t)\vert_{V_2}$ are Fredholm, and let us assume
that $t=0$ is a degeneracy instant for both $S_1$ and $S_2$.
Obviously, in general it will be $\spfl\big(S_1,[-\varepsilon,\varepsilon]\big)
\ne \spfl\big(S_2,[-\varepsilon,\varepsilon]\big)$, even in the case
that $\Ker\big(S_1(0)\big)=\Ker\big(S_2(0)\big)\subset V_1\cap V_2$.
On the other hand, if $W_k(S_1,0)=W_k(S_2,0)$ for all $k\ge1$,
then obviously   it must be $n^\pm_k(S_1,0)=n^\pm_k(S_2,0)$ for all $k\ge1$, hence
$\spfl\big(S_1,[-\varepsilon,\varepsilon]\big)
=\spfl\big(S_2,[-\varepsilon,\varepsilon]\big)$ for $\varepsilon>0$ small enough.
This follows easily from the observation that the derivatives of $S_i$
coincide with the restriction to $V_i$ of the derivatives of $S$, $i=1,2$.

Assuming that $S(t)$ has a degeneracy instant at $t=0$ on $\Hcal$,
it is an interesting question to ask how to find a ``minimal''
closed subspace $V$ of $\Hcal$ containing $\Ker\big(S(0)\big)$ and
with the property that the restriction of $S(t)$ to $V$ has the same
partial signatures at $t=0$. Denote by $\vv_i(t)$, $i=1,\ldots,N=\Dim\big(\Ker(S(0))\big)$,
a smooth orthonormal family of unit eigenvectors corresponding to the eigenvalues $\lambda_i(t)$ of $S_i(t)$
such that $\lambda_i(0)=0$. Examining the proof of Proposition~\ref{thm:central}
suggests that the desired subspace $V$ can be obtained by considering the span
of all the vectors $\vv_i(0)$,
together with their derivatives $\vv_i'(0)$, $\vv_i''(0)$, \ldots, $\vv_i^{(r_i)}(0)$,
where $i=1,\ldots,N$ and $r_i>0$ is such that $\lambda_i^{(r_i-1)}(0)$ and
$\lambda_i^{(r_i)}(0)\ne0$. Such a space $V$ is finite dimensional, and repeating the argument
at each degeneracy instant of a real-analytic path
in $\Fcalsa(\Hcal)$ we have obtained the following:
\begin{prop}\label{thm:reductionfindim}
Given any real-analytic path ${\LL}:[a,b]\to\Lcalsa(\Hcal)$ there exists a finite dimensional
subspace $V$ of $\Hcal$ such that, denoting by $S(t)$ the bilinear form
$\langle{\LL}(t)\cdot,\cdot\rangle$ on $\Hcal$ and by $\widetilde S(t)$ its
restriction to $V$, then $S$ and $\widetilde S$ have precisely the same
degeneracy instants in $[a,b]$, and $n_k^\pm(S,t_0)= n_k^\pm(\widetilde S,t_0)$
for all degeneracy instant $t_0$. In particular, $\spfl\big(S,[a,b]\big)=
\spfl\big(\widetilde S,[a,b]\big)$.\qed
\end{prop}

\begin{example}\label{thm:exarestriction}
Let us consider a real separable Hilbert space $\Hcal$, a
self-adjoint isomorphism $g:\Hcal\to\Hcal$ and a $g$-symmetric
Fredholm operator $T:\Hcal\to\Hcal$, i.e., $T$ is such that
$\langle gT\cdot,\cdot\rangle$ is a symmetric bilinear
form on $\Hcal$; assume that $0$ is (an isolated point)
in the spectrum of $T$. Consider the following real-analytic
path of Fredholm self-adjoint operators on $\Hcal$:
${\LL}(t)=gT-tg$, $t\in[-\varepsilon,\varepsilon]$;
then, $t=0$ is an isolated degeneracy instant for ${\LL}$.
It is not hard to see that $T$ has Fredholm index $0$, that
there exists $n_0>0$ such that $\Ker(T^n)=\Ker(T^{n_0})$ for
all $n\ge n_0$, with $\Ker(T^{n_0})$ a finite dimensional
subspace of $\Hcal$ (see Proposition~\ref{thm:surpresa}).
If $S(t)$ denotes the symmetric bilinear form $\langle{\LL}(t)\cdot,\cdot\rangle$
on $\Hcal$ and $\widetilde S(t)$ its restriction to $\Ker(T^{n_0})$, then
$n_k^\pm(S,0)=n_k^\pm(\widetilde S,0)$ for all $k\ge1$. To see this
observe that a sequence $(u_0,u_1,\ldots,u_k)$ of vectors in $\Hcal$,
with $u_0\in\Ker(T)$,
is a generalized Jordan chain for $S$ at $t=0$ iff   $Tu_r=u_{r-1}$
for all $r=1,\ldots,k$.
This in particular implies $T^ru_r\in\Ker(T)$, i.e., $u_r\in\Ker(T^{r+1})$ for
all $r$, which proves our assertion.
\end{example}
\end{section}

\begin{section}{On the Maslov index of  Lagrangian
paths}\label{sec:maslovindex} We will apply the result of
Section~\ref{sec:gensig} to the study of the Maslov index of a
path in the Lagrangian Grassmannian of a   symplectic space. Our
first goal is to give a general definition of Maslov index for
arbitrary continuous curves which is invariant by homotopies with
endpoints varying in a stratum of the Maslov cycle and additive by
concatenation. Such number will depend on the choice of a
Lagrangian $L_0$ that is used to define the Maslov cycle, except
for the case of closed paths (Corollary~\ref{thm:intfact2}). There
are several definitions of Maslov index available in the
literature, not always equivalent. Duistermaat's definition of
Maslov index in \cite{Duis} does not depend on the choice of
$L_0$, but it is not additive by concatenation. A semi-integer
valued Maslov index, which has the two properties above has been
introduced by Robbin and Salamon in \cite{RobSal} by considering
first regular curves having transversal intersections with the
Maslov cycle and then extending by homotopy invariance. The
definition of Maslov index given by de Gosson (see \cite{deGos1,
deGos2, deGos3, deGos4, deGosdeGos}) is based on the notion of
{\em Leray's index\/} for pairs in the universal covering of the
Lagrangian Grassmannian $\Lambda$. Booss--Bavnbek and Furutani
have given in \cite{BoosFuru} a functional analytical definition
of Maslov index, both in the finite and infinite dimensional case,
of a Fredholm Lagrangian paths. This is obtained by defining a
one-parameter operator family  of operators associated to the
path, whose spectrum oscillates on the unit circle around
$e^{i\pi}$ in the complex plane, and giving an appropriate
algebraic count of the passages through a fixed gauge. The
construction is done locally, and then patched together following
Phillips' definition of spectral flow in \cite{Phillips}. A simple
relation holds between the indices of Robbin/Salamon and
Booss--Bavnbek/Furutani for smooth transversal paths (see
\cite[Section~2]{BoosFuru}). We simplify greatly the approach of
\cite{BoosFuru} using a van Kampen type theorem for the
fundamental groupoid of a topological space
(Proposition~\ref{thm:Uiphii}), that reduces the proof of the
well-definiteness of the Maslov index to a simple local
compatibility condition in the space of symmetric bilinear forms
(Proposition~\ref{thm:3Lagr}). In turn, the latter compatibility
condition is proven as an immediate application of our partial
signature theory; it should also be emphasized that our
construction allows to define index functions taking values in
arbitrary groups. We start with an abstract result on how to
construct group-valued homomorphisms on the fundamental groupoid
of a topological space, and then we will study the case of
Lagrangian paths. Several other notions of Maslov index will be
discussed in the last part of the Section.

\subsection{An abstract index theory for continuous paths}
Let $X$ be a topological space. We denote by $\pi(X)$ the {\em fundamental homotopy groupoid\/}
of $X$, i.e., the set of all fixed-endpoint homotopy classes $[\gamma]$ of continuous paths $\gamma:[0,1]\to X$,
endowed with the partial binary operation of concatenation
$\diamond$, i.e., $[\gamma]\cdot[\mu]=
[\gamma\diamond\mu]$,
where $\gamma\diamond\mu$ is defined by
$\gamma\diamond\mu(t)=\gamma(2t)$, for
$t\in\big[0,\frac12\big]$ and $\gamma\diamond\mu(t)=\mu(2t-1)$,
for
$t\in\big[\frac12,1\big]$; obviously the concatenation $[\gamma]\cdot[\mu]$ is defined whenever
$\gamma(1)=\mu(0)$. If $\gamma:[a,b]\to X$ is a continuous path defined on an arbitrary compact
interval $[a,b]$, we will identify $\gamma$ with the curve $[0,1]\ni t\mapsto\gamma\big(t(b-a)+a\big)$,
and we will write $[\gamma]\in\pi(X)$. Given a group $G$ then by a {\em $G$-valued homomorphism\/}
$\phi$ on $\pi(X)$ we mean a map $\phi:\pi(X)\to G$ such that $\phi\big([\gamma]\cdot[\mu]\big)=\phi\big([\gamma]\big)\phi\big([\mu]\big)$,
for all $[\gamma],[\mu]\in\pi(X)$ with $\gamma(1)=\mu(0)$.

\begin{example}\label{exa:sodosextremos}
Given an arbitrary map $\tau:X\to G$ then one can define a
$G$-valued homomorphism on $\pi(X)$ by setting
$\phi\big([\gamma]\big)=\tau\big(\gamma(0)\big)^{-1}\tau\big(\gamma(1)\big)$,
for all $[\gamma]\in\pi(X)$.
\end{example}

We will now discuss the existence (and the uniqueness) of a group
valued homomorphism defined globally on the fundamental groupoid
of a topological space, once that the values of such homomorphism
are given on ``short'' curves. The following proposition  is a van
Kampen type result for the fundamental groupoid  of a topological
space; the proof presented takes inspiration on the construction
of the spectral flow given in \cite[Proposition~3]{Phillips}.
Recall that the classical statement of van Kampen's theorem in
elementary algebraic topology gives the fundamental group
$\pi_1(X)$ of a topological space $X$ as the quotient of a free
product of fundamentals groups $\pi_1(U_i)$, where $\{U_i\}_{i\in
I}$ is an open cover of $X$. Such quotient is characterized by a
universal property concerning the existence and uniqueness of an
extension  $\phi:\pi_1(X)\to G$ of a family of group homomorphisms
$\phi_i:\pi_1(U_i)\to G$, provided that some compatibility
assumption on the $\phi_i$ is satisfied. This is precisely the
spirit of our fundamental groupoid version of this result:

\begin{prop}\label{thm:Uiphii}
Let $X$ be a topological space, $G$ a group and $X=\bigcup_{i\in
I}U_i$  an open cover of $X$. Assume that it is given a $G$-valued
homomorphism $\phi_i$ on $\pi(U_i)$ for each $i\in I$ such that
$\phi_i\big([\gamma]\big)=\phi_j\big([\gamma]\big)$ for all
$[\gamma]\in\pi(U_i\cap U_j)$ and all $i,j\in I$. Then, there
exists a unique $G$-valued homomorphism $\phi$ on $\pi(X)$ such
that $\phi\big([\gamma]\big)=\phi_i\big([\gamma]\big)$ for all
$\gamma\in\pi(U_i)$ and all $i\in I$.
\end{prop}
\begin{proof}
Let $\gamma:[0,1]\to X$ be a continuous curve; by compactness
there exists a partition $0=t_0<t_1<\ldots<t_{N+1}=1$
of the interval $[0,1]$ and a map $r:\{0,1,\ldots,N\}\to
I$ such that $\gamma\big([t_k,t_{k+1}]\big)\subset
U_{r(k)}$ for all $k=0,\ldots,N$ (the choice of both the
partition and the map $r$ is highly non unique).
Define:
\begin{equation}\label{eq:deftildephi}
\widetilde\phi(\gamma):=\phi_{r(0)}
\big(\big[\gamma\vert_{[t_0,t_{1}]}\big]\big)\cdot\phi_{r(1)}
\big(\big[\gamma\vert_{[t_1,t_{2}]}\big]\big)\cdot\ldots\cdot
\phi_{r(N)}
\big(\big[\gamma\vert_{[t_N,t_{N+1}]}\big]\big);
\end{equation}
in order to have a well-defined map on $\pi(X)$
we need to show that the value $\widetilde\phi(\gamma)$ does
not depend on the choice of the partition $(t_k)_k$
and of the function $r$ as above.
In first place, one observes that passing to a finer
partition does not change the value of
$\widetilde\phi$, by the concatenation multiplicativity of
the maps $\phi_i$. Secondly, if for some $k$
it happens that $\gamma\big([t_k,t_{k+1}]\big)\subset
U_i\cap U_j$, then
$\phi_i\big(\big[\gamma\vert_{[t_k,t_{k+1}]}\big]\big)=
\phi_j\big(\big[\gamma\vert_{[t_k,t_{k+1}]}\big]\big)$, so that
the value of $\widetilde\phi$ does not depend on the
choice of the function $r$. Finally, given two
different partitions of
$[0,1]$,  one can find a third partition which is finer
than both of them, and by what has been observed, this
implies that $\widetilde\phi$ is well defined.

Uniqueness is also clear, since the concatenation multiplicativity
implies that \eqref{eq:deftildephi} must hold.

As to the homotopy invariance, observe first that the multiplicativity
by concatenation and the homotopy invariance of the $\phi_i$'s
imply that:
\begin{itemize}
\item if $\gamma$ is a constant curve, then
$\widetilde\phi(\gamma)=0$;
\item if $\gamma^-$ denotes the backwards reparameterization
of $\gamma$, then
$\widetilde\phi(\gamma^-)=\widetilde\phi(\gamma)^{-1}$.
\end{itemize}
Assume that $p$, $q$ are points of $X$ and $H:[0,1]\times[0,1]\to
X$ is a continuous map with $H(s,0)=p$ and $H(s,1)=q$ for all
$s\in[0,1]$. Set $\gamma_0=H(0,\cdot)$ and $\gamma_1=H(1,\cdot)$;
we want to prove that
$\widetilde\phi(\gamma_0)=\widetilde\phi(\gamma_1)$. Choose
partitions $0=s_0<s_1<\ldots<s_{N+1}=1$ and
$0=t_0<t_1<\ldots<t_{M+1}=1$ of the interval $[0,1]$
such that, for all $k=0,\ldots,N$ and all $j=0,\ldots,M$, the
image $H\big([s_k,s_{k+1}] \times [t_j,t_{j+1}]\big)$ is contained
in some $U_{i}$. For all $k=0,\ldots,N+1$, set
$\gamma_{s_k}=H(s_k,\cdot)$; we will prove that
$\widetilde\phi(\gamma_{s_k}\big)
=\widetilde\phi\big(\gamma_{s_{k+1}}\big)$ for all $k$, which will
conclude the argument. Define the following curves (see
Figure~\ref{fig:homotopy}):
\begin{figure}
\begin{center}
\psfull \epsfig{file=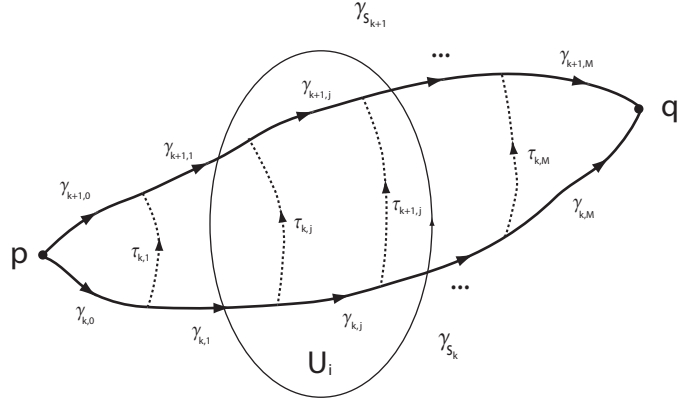}
\caption{The curves defined in
the proof of Proposition~\ref{thm:Uiphii}}\label{fig:homotopy}
\end{center}
\end{figure}
\[\begin{split}&\tau_{k,j}=H(\cdot,t_j)\Big\vert_{[s_k,s_{k+1}]},
\\
&\gamma_{k,j}=\gamma_{s_k}\big\vert_{[t_j,t_{j+1}]},\qquad\qquad
k=0,\ldots,N,\quad j=0,\ldots,M,\end{split}\] and
\[\widetilde\gamma_{s_k}=\gamma_{k,0}
\diamond\tau_{k,1}\diamond\tau_{k,1}^-\diamond
\gamma_{k,1}\diamond
\tau_{k,2}\diamond\tau_{k,2}^-\diamond
\ldots\diamond \tau_{k,M}\diamond\tau_{k,M}^-\diamond
\gamma_{k,M}.\]
Clearly, $\widetilde\phi(\gamma_{s_k})=\widetilde\phi(\widetilde
\gamma_{s_k})$, because
$\widetilde\phi(\tau_{k,j}\diamond\tau_{k,j}^-)=0$ for
all $j$. Moreover, $\gamma_{k,0}
\diamond\tau_{k,1}$ is homotopic with fixed endpoints
to $\gamma_{k+1,0}$ in some $U_i$,
 $\tau_{k,1}^-\diamond
\gamma_{k,1}\diamond
\tau_{k,2}$ is homotopic with fixed endpoints
to $\gamma_{k+1,1}$ in some
other $U_i$, etc., and $\tau_{k,M}^-\diamond
\gamma_{k,M}$ is homotopic with
fixed endpoints to $\gamma_{k+1,M}$
in still some other $U_i$. Since the $\phi_i$'s are
invariant by fixed endpoint homotopies, the multiplicativity by
concatenation once more implies that $\widetilde
\phi(\widetilde\gamma_{s_k})=\widetilde\phi(\gamma_{s_{k+1}})$.
The desired map $\phi$ is thus obtained by setting
$\phi\big([\gamma]\big)=\widetilde\phi(\gamma)$.
\end{proof}

\begin{example}\label{exa:taui}
Assume that $X=\bigcup_{i\in I}U_i$ is an open cover of
$X$ and that for each $i\in I$ one is given a map
$\tau_i:U_i\to G$ such that the map $U_i\cap U_j\ni
g\mapsto\tau_i(g)^{-1}\tau_j(g)\in G$ is constant on each
arc-connected component of $U_i\cap U_j$, for all $i,j\in I$.
By Example~\ref{exa:sodosextremos} and
Proposition~\ref{thm:Uiphii}, there exists a unique $G$-valued
homomorphism $\phi$ on $\pi(X)$ such that
$\phi\big([\gamma]\big)=\tau_i\big(\gamma(0)\big)^{-1}\tau_i\big(\gamma(1)\big)$
for all $\gamma\in\pi(U_i)$ and all $i\in I$.
\end{example}

\subsection{Charts on the Lagrangian Grassmannian of a
symplectic space}\label{sub:charts}
Throughout this subsection, $(V,\omega)$ will denote a
finite dimensional symplectic
space, i.e., $V$ is a real  vector space,
and $\omega$ is a nondegenerate skew-symmetric bilinear
form on $V$; set $\Dim(V)=2n$. The {\em symplectic group\/}
$\Spl(V,\omega)$ is the closed subgroup of $\mathrm{GL}(V)$
consisting of those linear maps on $V$ that preserve $\omega$.

Recall that a {\em Lagrangian subspace\/} of $V$ is an
$n$-dimensional subspace $L\subset V$ on which
$\omega$ vanishes. The set of all Lagrangian
subspaces of $V$, denoted by $\Lambda$, has the structure
of a compact, real-analytic submanifold of the Grassmannian
of all $n$-dimensional subspaces of $V$. The dimension
of $\Lambda$ equals $\frac12n(n+1)$, and a real-analytic
atlas on $\Lambda$ is given as follows.

For all $L\in\Lambda$ and $k\in\{0,\ldots,n\}$, set:
\[\Lambda_k(L)=\big\{L'\in\Lambda:\Dim(L\cap L')=k\big\};\]
in particular, $\Lambda_0(L)$ is the set of all Lagrangian
subspaces that are transversal to $L$, and it is a dense
open subset of $\Lambda$. Given a pair $L_0,L_1\in\Lambda$
of complementary Lagrangians, i.e., $L_0\cap L_1=\{0\}$, then
one defines a map:
\[\varphi_{L_0,L_1}:\Lambda_0(L_1)\longrightarrow\Bsym(L_0)\]
as follows. Any Lagrangian $L\in\Lambda_0(L_1)$ is the graph
of a unique linear map $T:L_0\to L_1$; then, $\varphi_{L_0,L_1}$
is defined to be the restriction of the bilinear map
$\omega(T\cdot,\cdot)$ to $L_0\times L_0$. It is easy to see
that, due to the fact that $L$ is Lagrangian, such bilinear map is symmetric.
Observe that:
\begin{equation}\label{eq:kerphil0l1}
\Ker\big(\varphi_{L_0,L_1}(L)\big)=L\cap L_0,\quad\forall\,L\in\Lambda_0(L_1).
\end{equation}

The collection of all $\varphi_{L_0,L_1}$, when $(L_0,L_1)$ runs
in the set of all pairs of complementary Lagrangians, is a
real-analytic atlas on $\Lambda$. If $L_1$ and $L_1'$ are
complementary to a given $L_0$, then the transition function:
\[\varphi_{L_0,L_1'}\circ\varphi_{L_0,L_1}^{-1}:\varphi_{L_0,L_1}
\big(\Lambda_0(L_1')\big)\subset\Bsym(L_0)\longrightarrow
\Bsym(L_0)\]
is given by:
\begin{equation}\label{eq:cambcarta}
\varphi_{L_0,L_1'}\circ\varphi_{L_0,L_1}^{-1}(B)=
B\circ\big(\mathrm{Id}+(\pi_0'\vert_{L_1})\circ\rho_{L_0,L_1}^{-1}\circ B\big)^{-1},
\end{equation}
where $\pi_0':L_0\oplus L_1'\to L_0$ is the projection onto
the first summand, and $\rho_{L_0,L_1}:L_1\to L_0^*$ is the map
$v\mapsto\omega(v,\cdot)\vert_{L_0}$, which is an isomorphism. Observe that in formula \eqref{eq:cambcarta}
the bilinear form $B$ is seen as a map $L_0\to L_0^*$.
The map $\mathrm{Id}+(\pi_0'\vert_{L_1})\circ\rho_{L_0,L_1}^{-1}\circ B$
is an automorphism of $L_0$ whose inverse, denoted by $\mathfrak h$,
satisfies $B\circ \mathfrak h=\mathfrak h^*\circ B$; observe that $\mathfrak h$ is the identity
on $\Ker(B)$.

Recall that every symplectic
space is  isomorphic $\R^n\oplus{\R^n}^*$ endowed with the {\em standard
symplectic form}:
\[\omega_0\big((v,\alpha),(w,\beta)\big)=\beta(v)-\alpha(w).\]
More generally, given Lagrangians $L_0,L_1\subset V$ with
$L_0\cap L_1=\{0\}$, then there exists a symplectic isomorphism
$\phi:V\to\R^n\oplus{\R^n}^*$ (i.e., the pull-back $\phi_*\omega_0$
coincides with $\omega$)
such that
$\phi(L_0)=\{0\}\oplus {\R^n}^*$ and
$\phi(L_1)=\R^n\oplus\{0\}$.

\subsection{Maslov index of a symplectic
path}\label{sub:maslovsymplpath}
Let now $L_0$ be a fixed Lagrangian in $\Lambda$.

The $L_0$-Maslov index of a continuous curve
$\gamma$ with endpoints in $\Lambda_0(L_0)$ can
be defined in terms of the first relative homology
group of the pair $\big(\Lambda,\Lambda_0(L_0)\big)$
(see \cite{pacific, topology}). We will give an
alternative and more general definition of Maslov
index in the case of continuous curves in
$\Lambda$ with arbitrary endpoints.

\begin{prop}\label{thm:3Lagr}
Given Lagrangians $L_0,L_1,L_1'$ in $V$ with
$L_1,L_1'\in\Lambda_0(L_0)$, then the map:
\[\Lambda_0(L_1)\cap\Lambda_0(L_1')\ni L\longmapsto
{\noo+}\big(\varphi_{L_0,L_1'}(L)\big)-
{\noo+}\big(\varphi_{L_0,L_1}(L)\big)\in\Z\]
is constant on each connected component of the open set
$\Lambda_0(L_1)\cap \Lambda_0(L_1')\subset\Lambda$.
\end{prop}
\begin{proof}
In first place we observe that $\Ker\big(\varphi_{L_0,L_1'}(L)\big)=L\cap L_0=\Ker\big(\varphi_{L_0,L_1}(L)\big)$,
hence
\[{\noo+}\big(\varphi_{L_0,L_1'}(L)\big)-
{\noo+}\big(\varphi_{L_0,L_1}(L)\big)=  n^+\big(\varphi_{L_0,L_1'}(L)\big)-
n^+\big(\varphi_{L_0,L_1}(L)\big).\]
As we have observed, it is not restrictive to assume
that $V=\R^n\oplus{\R^n}^*$ is the standard symplectic space,
 $L_0=\{0\}\oplus{\R^n}^*$ and $L_1=\R^n\oplus\{0\}$.
Then $L_1'$ is of the form:
\[L_1'=\big\{(v,Zv):v\in\R^n\big\},\]
where $Z:\R^n\to{\R^n}^*$ is some symmetric linear map.
Given $L\in\Lambda_0(L_1)$, we set $B=\varphi_{L_0,L_1}(L)
\in\Bsym({\R^n}^*)$;  it is easily checked that
$L\in\Lambda_0(L'_1)$ if and only if $\Id-ZB$ is invertible and
that $\widetilde B=\varphi_{L_0,L'_1}(L)$ is given by:
\begin{equation}\label{eq:BBtilde}
\widetilde B=B(\Id-ZB)^{-1}.
\end{equation}
The proof will be completed once we show that $n^+(\widetilde B)-n^+(B)$ is constant
on each connected component of the open set $\mathcal B\subset\Bsym({\R^n}^*)$ consisting
of those $B$ with $\Id-ZB$ invertible, where $\widetilde B$ is defined by \eqref{eq:BBtilde}.
Given $B_0$ and $B_1$ in some connected component of $\mathcal B$ then we can find a
real-analytic curve $B:[0,1]\to\mathcal B$, with $B(0)=B_0$,
$B(1)=B_1$. Setting $\widetilde B(t)=B(t)\circ h(t)$, with
$h(t):=\big(\Id-ZB(t)\big)^{-1}$, then
$\widetilde B:[0,1]\to\Bsym({\R^n}^*)$ is a
real-analytic curve of symmetric bilinear forms;
observe that $h(t)$ is the identity on $\Ker(B(t))$
for all $t\in[0,1]$. So, by Proposition~\ref{thm:secinv},
$n_k^\pm(B,t_0)=n_k^\pm(\widetilde B,t_0)$,
$\sigma_k(B,t_0)=\sigma_k(\widetilde B,t_0)$
for all $t_0\in[0,1]$ and for all $k\in\N$. From
Proposition~\ref{thm:central}, part~\eqref{itm:formprinc} it
follows that:
\[n^+\big(B(1)\big)-n^+\big(B(0)\big)=n^+\big(\widetilde B(1)\big)-n^+\big(\widetilde B(0)\big).\]
Hence $n^+\big(\widetilde B(1)\big)-n^+\big(B(1)\big)=n^+\big(\widetilde B(0)\big)-n^+\big(B(0)\big)$,
which concludes the proof.
\end{proof}
From Proposition~\ref{thm:Uiphii},
Example~\ref{exa:taui}
and Proposition~\ref{thm:3Lagr} we obtain immediately:

\begin{cor}\label{thm:esunmaslovindex}
For all $L_0\in\Lambda$, there exists a unique $\Z$-valued
groupoid homomorphism
$\mu_{L_0}$ on $\pi(\Lambda)$ such that:
\begin{equation}\label{eq:muL0tilde}
\mu_{L_0}\big([\gamma]\big)={\noo+}\big(\varphi_{L_0,L_1}(\gamma(1))\big)-
{\noo+}\big(\varphi_{L_0,L_1}(\gamma(0))\big)
\end{equation}
for all continuous curve $\gamma:[0,1]\to\Lambda_0(L_1)$ and for all
$L_1\in\Lambda_0(L_0)$.\qed
\end{cor}
\begin{rem}\label{thm:remcoindexorextendedcoindex?}
From what has been observed in the proof of Proposition~\ref{thm:3Lagr},
it is clear that the result of Corollary~\ref{thm:esunmaslovindex}  holds also
if one replaces the extended coindex ${\noo+}$ with the
coindex $n^+$ in \eqref{eq:muL0tilde}. When the endpoints of $\gamma$
are transversal to $L_0$, then such distinction does not affect
the value of $\mu_{L_0}$, while in the case of degenerate endpoints
the function $\mu_{L_0}$ obtained using the coindex would give a different
value. The choice of one or another definition is merely a matter of personal taste;
on the other hand, such choice should be made consistently with the choice
of a notion of spectral flow in the case of degenerate endpoints.
In our case, the Maslov index defined in Corollary~\ref{thm:esunmaslovindex}
is such that, when applied to the case of Riemannian geodesics
(see Subsection~\ref{sub:partsigconjinst}), it gives the total number of
conjugate points along the geodesic, including
that possibly  occurring at the final instant.
\end{rem}

\begin{defin}
\label{thm:defmuL0}
Given any continuous curve $\gamma$ in $\Lambda$, the integer
$\mu_{L_0}(\gamma):=\mu_{L_0}([\gamma])$ will be called
the {\em $L_0$-Maslov index} (or simply, the {\em Maslov index\/}
when the choice of $L_0$ is clear from the context) of the curve $\gamma$.
\end{defin}
Observe that the homotopy invariance implies, in particular,
that $\mu_{L_0}(\gamma)$ is independent on the parameterization
of $\gamma$.
We collect in the following statement the main properties of
the integer valued map $\mu_{L_0}$, whose proof follows almost immediately
from the definition:
\begin{lem}\label{thm:propmuL0}
The Maslov index $\mu_{L_0}$ satisfies the following properties:
\begin{enumerate}
\item\label{itm:propmuL01} replacing the symplectic form $\omega$ by $-\omega$ produces
a change in the sign  of $\mu_{L_0}$;\smallskip

\item\label{itm:propmuL02} given a continuos curve $\gamma:[a,b]\to\Lambda$, if $\Dim\big(\gamma(t)\cap L_0\big)$
is constant on $[a,b]$, then $\mu_{L_0}(\gamma)=0$;
\smallskip

\item\label{itm:propmuL03} if $H:[a,b]\times[c,d]\to\Lambda$ is a continuous map,
then:
\[\mu_{L_0}\big(t\mapsto H(t,c)\big)+\mu_{L_0}\big(s\mapsto H(b,s)\big)=
\mu_{L_0}\big(t\mapsto H(t,d)\big)+\mu_{L_0}\big(s\mapsto H(a,s)\big);\]

\item\label{itm:propmuL04} \emph{(symplectic invariance)} if $\phi:(V,\omega)\to(V',\omega')$ is a symplectomorphism and
$\gamma:[a,b]\to\Lambda(V,\omega)$ is continuous, then:
\[\mu_{L_0}(\gamma)=\mu_{\phi(L_0)}\big(\phi\circ\gamma\big);\]

\item \emph{(symplectic additivity)} if $\gamma:[a,b]\to(V,\omega)$ and $\widetilde\gamma:[a,b]\to(
\widetilde V,\widetilde \omega)$
are continuous, $L_0\in\Lambda(V,\omega)$ and $\widetilde L_0\in\Lambda(\widetilde
V,\widetilde\omega)$, then
\[\mu_{L_0\oplus\widetilde L_0}\big(\gamma\oplus\widetilde\gamma)=
\mu_{L_0}(\gamma)+\mu_{\widetilde L_0}(\widetilde \gamma).\qed\]
\end{enumerate}
\end{lem}
It will be useful to single out the following additional properties
of the Maslov index:

\begin{cor}\label{thm:admuL0}
The Maslov index $\mu_{L_0}$ satisfies also the following:
\begin{itemize}
\item[(a)] If $\gamma_1,\gamma_2:[a,b]\to\Lambda$ are continuous
curves   that are homotopic by a homotopy with free endpoints
in some  $\Lambda_k(L_0)$, i.e., there exists a continuous map
$H:[0,1]\times[a,b]\to\Lambda$ such that $H(0,\cdot)=\gamma_1$,
$H(1,\cdot)=\gamma_2$, and with $\Dim\big(H(s,a)\cap L_0\big)$ and
 $\Dim\big(H(s,b)\cap L_0\big)$ constant on $[0,1]$, then
$\mu_{L_0}(\gamma_1)=\mu_{L_0}(\gamma_2)$;
\smallskip
\item[(b)] if $\gamma_1$ and $\gamma_2$ are continuous loops
in $\Lambda$ that are {\em freely\/} homotopic, then
$\mu_{L_0}(\gamma_1)=\mu_{L_0}(\gamma_2)$;
\smallskip
\item[(c)] if $[a,b]\ni t\mapsto\phi(t)$ is a continuous curve
in the symplectic group $\Spl(V,\omega)$ such that
$\phi(t)\big(L_0)=L_0$ for all $t$, and $\gamma:[a,b]\to\Lambda$
is continuous, then:
\[\mu_{L_0}(\gamma)=\mu_{L_0}\big(t\mapsto\phi(t)(\gamma(t))\big).\]
\end{itemize}
\end{cor}
\begin{proof}
(a) and (b) follow easily from part~\eqref{itm:propmuL03}
of Lemma~\ref{thm:propmuL0}. Observe indeed that, if $H$ denotes
the given homotopy between $\gamma_1$ and $\gamma_2$, then
the curves $s\mapsto H(a,s)$ and $s\mapsto H(b,s)$ coincide
in the case (b) of freely homotopic loops; they are curves in
$\Lambda_k(L_0)$ for some fixed $k$ in the case (a), so their Maslov
index vanishes by Lemma~\ref{thm:propmuL0} part~\eqref{itm:propmuL02}.

To prove (c), observe that the curve $t\mapsto\phi(t)(\gamma(t))$ is homotopic
to the curve $t\mapsto\phi(a)(\gamma(t))$ by the homotopy:
\[H(s,t)=\phi\big((1-s)t\big)\big(\gamma(t)\big).\]
Since $\phi(t)(L_0)=L_0$ for all $t$, by part~\eqref{itm:propmuL04}
of Lemma~\ref{thm:propmuL0} it follows that
\[\mu_{L_0}\big(t\mapsto\phi(a)(\gamma(t))\big)=\mu_{L_0}
\big(t\mapsto\phi(t)(\gamma(t))\big).\]
Moreover,  the above homotopy is by curves with endpoints varying in some
fixed $\Lambda_k(L_0)$, and the conclusion follows from part (a).
\end{proof}
The above definition of Maslov index
generalizes that in \cite{pacific, topology},
where $\mu_{L_0}$ had been defined in terms of the first
relative homology group  $H_1(\Lambda,\Lambda_0(L_0))$.
Some properties of $\mu_{L_0}$ that were discussed in \cite{pacific, topology}
using homological techniques and functoriality properties will
have to be reproven in this more general context; to this aim
we will now discuss a method for computing $\mu_{L_0}(\gamma)$ when
$\gamma$ is a real-analytic curve.

\subsection{Partial signatures and Maslov index}
\label{sub:partsignMaslind}
We have seen that the computation of the
Maslov index of a curve in $\Lambda$ is reduced using local
charts to the study of the jumps of the extended coindex of curves
of symmetric bilinear forms, and this is where the result
of Proposition~\ref{thm:central} comes about.

Assume that $\gamma:[a,b]\to\Lambda$ is a smooth curve. From
\eqref{eq:kerphil0l1} it is clear that   the jumps of the coindex
of $\varphi_{L_0,L_1}(\gamma)$ occur precisely at those instants
$t\in\left]a,b\right[$ when $\gamma(t)$ intersects the set:
\[\LLL:=\bigcup_{k\ge1}\Lambda_k(L_0).\]
We will call  ${\LLL}$ the {\em Maslov cycle with vertex at
$L_0$}. We remark here that each $\Lambda_k(L_0)$ is an embedded
submanifold of $\Lambda$, but that ${\LLL}$ is {\em not\/} a
submanifold of $\Lambda$. Assume that $t_0\in\left]a,b\right[$ is
an isolated intersection of $\gamma$ with ${\LLL}$ and choose any
Lagrangian $L_1\in\Lambda_0(L_0)$ which is transversal also to
$\gamma(t_0)$. By continuity, $L_1$ is transversal to $\gamma(t)$
for $t$ near $t_0$, and we can define a smooth curve
$[t_0-\varepsilon,t_0+\varepsilon]\ni
t\mapsto\varphi_{L_0,L_1}\big(\gamma(t)\big)\in \Bsym(L_0)$, for
$\varepsilon>0$ small enough. As we have observed, such a curve
has an isolated degeneracy instant at $t=t_0$, and we can define
the sequences:
\begin{eqnarray*}n^-_k(\gamma,t_0;L_0,L_1),  \noo-_k(\gamma,t_0;L_0,L_1),\
n^+_k(\gamma,t_0;L_0,L_1),\\ \noo+_k(\gamma,t_0;L_0,L_1),\
\sigma_k(\gamma,t_0;L_0,L_1),\quad k\ge1,\end{eqnarray*}
respectively as the partial (extended) indexes, partial (extended)
coindexes and partial signatures of the curve of symmetric bilinear forms $ \varphi_{L_0,L_1}\circ \gamma  $
at $t=t_0$.
\begin{lem}\label{thm:nondip}
The integers $n_k^\pm(\gamma,t_0;L_0,L_1)$, $\noo\pm_k(\gamma,t_0;L_0,L_1)$ and
$\sigma_k(\gamma,t_0;L_0,L_1)$ do not depend
on the choice of the Lagrangian $L_1$.
\end{lem}
\begin{proof}
Choose two Lagrangian spaces $L_1,L_1'\in\Lambda_0\big(\gamma(t_0)\big)\cap \Lambda_0(L_0)$
and set $L(t)=\varphi_{L_0,L_1}\big(\gamma(t)\big)$,
$\widetilde L(t)=\varphi_{L_0,L_1'}\big(\gamma(t)\big)$
for $t$ sufficiently near $t_0$.

Using formula \eqref{eq:cambcarta} we get that $\widetilde L(t)=L(t)h(t)$,
where $h(t)$ is a curve of automorphisms of $L_0$ such that
$h(t_0)$ is the identity on $\Ker\big(L(t_0)\big)=\gamma(t_0)\cap L_0$.
The conclusion follows immediately from Proposition~\ref{thm:secinv}.
\end{proof}

\noindent\ \ We are now entitled to talk about the  $L_0$-partial
signatures $n^\pm_k(\gamma,t_0;L_0)$,
${\noo\pm}_k(\gamma,t_0;L_0)$ and $\sigma_k(\gamma,t_0;L_0)$ at an
isolated intersection of $\gamma$ with ${\LLL}$, without
specifying the choice of a Lagrangian $L_1$.

\begin{prop}\label{thm:maslovpartsign}
Let $\gamma:[a,b]\to\Lambda$ be a real-analytic curve which is not
entirely contained in the Maslov cycle ${\LLL}$. Then, the
$L_0$-Maslov index of $\gamma$ is given by:
\begin{multline}\label{eq:maslovpartsign}
\mu_{L_0}(\gamma)=\sum_{\substack{t_0\in\gamma^{-1}
\big(\LLL\big)\\
t_0\in\left]a,b\right[}}\left[\sum_{k\ge1}
\sigma_{2k-1}\big(\gamma,t_0;L_0)\right]-\sum_{k\ge1}n^-_k(\gamma,a;L_0)\\+ \sum_{k\ge1}
\big(n^+_{2k-1}(\gamma,b;L_0)+n^-_{2k}(\gamma,b;L_0)\big),
\end{multline}
where all the sums on the right hand side of
\eqref{eq:maslovpartsign} have a finite number of non zero
terms.
\end{prop}
\begin{proof}
In first place we observe that $\gamma$ has at most a finite
number of intersections with ${\LLL}$; namely, in local
coordinates $\varphi_{L_0,L_1}$, such intersections correspond to
zeroes of the real-analytic function
$t\mapsto\det\big(\varphi_{L_0,L_1}(\gamma(t)\big)$. Such function
is not identically zero because $\gamma$ is not entirely contained
in the Maslov cycle.

From the definition, given such a curve $\gamma$, its Maslov index
$\mu_{L_0}(\gamma)$ is given by the sum of the jumps of the
coindex function of $\varphi_{L_0,L_1}\circ\gamma$ at the instants
in $\gamma^{-1}\big({\LLL}\big)$, and the conclusion follows
readily from part~\eqref{itm:formprinc} of
Proposition~\ref{thm:central} and formulas~\eqref{eq:nuoveformule}, observing that
\[\Ker\big(\varphi_{L_0,L_1}\circ\gamma(t)\big)= \gamma(t)\cap
L_0.\qedhere\]
\end{proof}
We have proven in Corollary~\ref{thm:admuL0} that, given
a continuous curve of Lagrangians $\gamma:[a,b]\to\Lambda$ and a continuous
curve of symplectomorphisms $\phi:[a,b]\to\Spl(V,\omega)$ such that
$\phi(t)(L_0)=L_0$ for all $t$, the $L_0$-Maslov
index of the curves $t\mapsto\gamma(t)$ and $t\mapsto\phi(t)(\gamma(t))$ coincide.
We will conclude this subsection with the proof that, in the smooth case,
also all the partial signatures of the curves $t\mapsto\gamma(t)$
and $t\mapsto\phi(t)(\gamma(t))$ at each intersection with the Maslov cycle coincide:
\begin{lem}\label{thm:compsympl}
Let $\gamma:[a,b]\to\Lambda$ a continuous curve and $t_0\in[a,b]$
such that $\gamma(t_0)\in{\LLL}$, and let
$\phi:[a,b]\to\Spl(V,\omega)$ be a smooth curve of
symplectomorphisms such that $\phi(t)(L_0)=L_0$ for all $t$.
Denote by $\widetilde\gamma$ the  curve in $\Lambda$ given by
$\widetilde\gamma(t)=\phi(t)\big(\gamma(t)\big)$; then
$\widetilde\gamma(t_0)\in{\LLL}$  and all the partial signatures
of $\widetilde\gamma$ at $t_0$ coincide with those of $\gamma$.
\end{lem}
\begin{proof}
Since $\phi(t_0)$ is an isomorphism and $\phi(t_0)(L_0)=L_0$, then
$\phi(t_0)\big(\Lambda_k(L_0)\big)=\Lambda_k(L_0)$ for all $k$,
hence $\widetilde\gamma(t_0)\in{\LLL}$. To prove the invariance of
the partial signatures, let us choose a Lagrangian
$L_1\in\Lambda_0(L_0) \cap\Lambda_0\big(\gamma(t_0)\big)$; observe
that then $L_1'=\phi(t_0)(L_1)\in
\Lambda_0(L_0)\cap\Lambda_0\big(\phi(t_0)(\gamma(t_0))\big)$. The
partial signatures of $\gamma$ at $t_0$ are computed using the
curve of symmetric bilinear forms
$B(t)=\varphi_{L_0,L_1}\circ\gamma(t)$ on $L_0$, while the partial
signatures of $\widetilde\gamma$ at $t_0$ are computed using the
curve $\widetilde B(t)=\varphi_{L_0,L_1'}\circ\widetilde
\gamma(t)$, for $t$ near $t_0$.

Set $C(t)=\varphi_{L_0,\phi(t)(L_1)}\big(\phi(t)(\gamma(t))\big)\in\Bsym(L_0)$;
an immediate calculation using the very definition of the charts
$\varphi_{L_0,L_1}$ shows that $C(t)$ is the pull-back of the bilinear form
$\varphi_{L_0,L_1}\circ\gamma(t)$ by the isomorphism $\phi(t)^{-1}:L_0\to
L_0$. By Proposition~\ref{thm:invcogredient}, the curves
$B(t)$ and $C(t)$ have the same partial signatures at the
degeneracy instant $t_0$.

Recalling formula \eqref{eq:cambcarta}, we compute:
\[ \widetilde B(t)=\varphi_{L_0,L_1'}\circ
\widetilde\gamma(t)=\varphi_{L_0,L_1'}\circ
\varphi_{L_0,\phi(t)(L_1)}^{-1}\big(C(t)\big)=
C(t) \big(\mathrm{Id}+Z(t)  C(t)\big)^{-1},\]
where $Z(t):L_0^*\to L_0$ is a homomorphism depending smoothly on $t$.
Observing that $h(t)=\big(\mathrm{Id}+Z(t)  C(t)\big)^{-1}$ is
the identity on $\Ker\big(C(t)\big)$, the conclusion follows
from Proposition~\ref{thm:secinv}.
\end{proof}

\subsection{Infinite dimensional Lagrangian Grassmannians}
Let us now consider an infinite dimensional,
separable real Hilbert space $\Hcal$ with inner
product $\langle\cdot,\cdot\rangle$ endowed with a symplectic
form $\omega$; let $J:\Hcal\to \Hcal$ be the corresponding
almost complex structure, i.e., $\omega=\langle J\cdot,\cdot\rangle$.
Then, $J$ is an anti-symmetric bounded operator on $\Hcal$
such that $J^2=-1$. A (necessarily closed) subspace
$L$ of $\Hcal$ will be called {\em Lagrangian\/}
if $L=L^{\perp_\omega}=(JL)^\perp$.

Recall that a pair $(V_1,V_2)$ of  closed subspaces of $\Hcal$
is called a {\em Fredholm pair\/} if $V_1\cap V_2$ is finite dimensional
and if the sum $V_1+V_2$ is closed and it has finite codimension in
$\Hcal$.

Given a pair $(L_0,L_1)$ of complementary Lagrangian
subspaces of $\Hcal$ (i.e., $L_0\cap L_1=(0)$ and $L_0+L_1=\Hcal$),
then one can define a real-analytic chart $\varphi_{L_0,L_1}$ on
$\Lambda_0(L_1)$ and taking values in $\Lcalsa(L_0)$ by setting
$\varphi_{L_0,L_1}(L)=P_0JS$, where $S:L_0\to L_1$ is
the unique  bounded operator whose graph
in $L_0\oplus L_1$ is $L$, and $P_0:\Hcal\to L_0$ is the orthogonal
projection.

Observe that $\Ker\big(\varphi_{L_0,L_1}(L)\big)=L\cap L_0$;
more generally, it is not hard to prove that $(L_0,L)$ is a Fredholm
pair if and only if $\varphi_{L_0,L_1}(L)\in\Fcalsa(L_0)$ for
some, hence for all, $L_1\in\Lambda_0(L_0)\cap\Lambda_0(L)$.
In particular, the set $\Fcal\Lambda$ of all Lagrangian subspaces $L$ of
$\Hcal$ such that $(L_0,L)$ is  Fredholm, is an open submanifold
of the manifold $\Lambda$ of all Lagrangian subspaces of $\Hcal$.%
\footnote{In the infinite dimensional case, the full Lagrangian
Grassmannian $\Lambda$ is contractible (see for instance \cite{Nico}).}
The entire theory described in subsection~\ref{sub:charts}
concerning the Lagrangian Grassmannian of a symplectic space
extends from the finite dimensional case to the case of the Fredholm
Lagrangian Grassmannian of an infinite dimensional symplectic
space; in particular, formulas of the type \eqref{eq:cambcarta}
hold for the transition maps $\varphi_{L_0,L_1'}\circ\varphi_{L_0,L_1}^{-1}\
$.

In order to define the notion of Maslov index for continuous paths
in $\Fcal\Lambda$, we prove the analogue of Corollary~\ref{thm:esunmaslovindex}
for the infinite dimensional case:
\begin{prop}\label{thm:esmaslindinfdimcase}
There exists a unique $\Z$-valued groupoid homomorphism
$\mu_{L_0}$ on $\pi(\mathcal F\Lambda)$ such that:
\begin{equation}\label{eq:mloinfdim}
\mu_{L_0}\big([\gamma]\big)=\spfl\big(\varphi_{L_0,L_1}\circ\gamma,[a,b]\big),
\end{equation}
for all Lagrangian $L_1$ complementary to $L_0$
and for all continuous curve $\gamma:[a,b]\to\Fcal\Lambda$ having image contained in
$\Lambda_0(L_1)$.
\end{prop}
\begin{proof}
Using the result of Proposition~\ref{thm:Uiphii}, it suffices to show
that the right hand side in formula \eqref{eq:mloinfdim}
gives a well defined $\Z$-valued groupoid homomorphism
on $\pi\big(\Lambda_0(L_1)\big)$, and that
\[\spfl\big(\varphi_{L_0,L_1}\circ\gamma,[a,b]\big)=
\spfl\big(\varphi_{L_0,L_1'}\circ\gamma,[a,b]\big)\]
if $\gamma$ has image in $\Lambda_0(L_1)\cap \Lambda_0(L_1')$.

This last equality follows easily from the formula of
the transition maps \eqref{eq:cambcarta} and the corresponding
invariance property for the spectral flow stated in
Corollary~\ref{thm:invSFLH}.

As to the  first part of the proof, the claim
is equivalent to the fact that the spectral flow
is invariant by fixed-endpoints homotopies and additive by concatenation.
These facts are proven in \cite{FitzPejsaRecht}
in the case that $\gamma$ has endpoints in $\Lambda_0(L_0)$,
i.e., that $\varphi_{L_0,L_1}\circ\gamma$ has invertible endpoints.
The very same proof in \cite{FitzPejsaRecht}, which uses formula \eqref{eq:spflreldim}
and it is based
on an argument of homotopy lifting in fiber bundles, holds for the general
case of possibly degenerate endpoints.
\end{proof}
Following closely the theory in
subsection~\ref{sub:partsignMaslind} one defines partial indexes,
coindexes and signatures at an isolated intersection of a smooth
curve $\gamma:[a,b]\to\Fcal\Lambda$ with the Maslov cycle
${\LLL}$, and as in the finite dimensional case one obtains the
following:
\begin{prop}\label{thm:calcmiinfdimcase}
Let $\gamma:[a,b]\to\Fcal\Lambda$ be a real-analytic curve
whose image is not entirely contained in the Maslov cycle.
Then, the Maslov index $\mu_{L_0}(\gamma)$ is given by
\eqref{eq:maslovpartsign}.\qed
\end{prop}

\subsection{On the notion of Maslov index for pairs of Lagrangian paths}
\label{sub:pairs} There exists in the literature a slightly
different notion of Maslov index for {\em pairs\/}
$(\gamma_1,\gamma_2)$ of continuous curves
$\gamma_1,\gamma_2:[a,b]\to\Lambda$ (see for instance
\cite{CapLeeMil}), as an integer valued measure of the set of
instants $t\in[a,b]$ at which the Lagrangians $\gamma_1(t)$
and $\gamma_2(t)$ are not transversal.
Using our  partial signatures theory, we will discuss below the
definition of Maslov index for such pairs; we will consider, as in
\cite{CapLeeMil}, the case of arbitrary pairs
$(\gamma_1,\gamma_2)$ without any transversality assumption at the
endpoints.

For all $L_0\in\Lambda$, consider the real-analytic fibration
$\beta_{L_0}:\Spl(2n,\R)\longrightarrow\Lambda$:
\begin{equation}\label{eq:defbetaL0}
\beta_{L_0}(\phi)=\phi(L_0);\end{equation} given any curve
$\gamma:[a,b]\to \Lambda$ of class $C^k$,
$k=1,\ldots,\infty,\omega$, it can be lifted to a curve
$\psi:[a,b]\to\Spl(2n,\R)$ of class $C^k$, i.e.,
$\gamma(t)=\psi(t)(L_0)$ for all $t\in[a,b]$. We will call such a
curve $\psi$ a {\em $L_0$-lifting\/} of $\gamma$. Observe that,
given two curves $\gamma_1,\gamma_2:[a,b]\to\Lambda$ and any
$L_0$-lifting $\psi:[a,b]\to\Spl(V,\omega)$ of $\gamma_2$, the non
transversality instants for the curves $\gamma_1$ and $\gamma_2$
correspond to the intersections of the curve
$t\mapsto\psi(t)^{-1}(\gamma_1(t))$ with the Maslov cycle
${\LLL}$.

In order to define the Maslov index of a pair $(\gamma_1,\gamma_2)$
one needs the following:
\begin{lem}\label{thm:maslovpairs}
Let $L_0\in\Lambda$ be fixed, let
$(\gamma_1,\gamma_2):[a,b]\to\Lambda\times\Lambda$ be a pair of
continuous curves in $\Lambda$, and let $\psi:[a,b]\to\Spl(2n,\R)$
be any $L_0$-lifting of $\gamma_2$. Then, the Maslov index
$\mu_{L_0}\big(t\mapsto\psi(t)^{-1}(\gamma_1(t))\big)$ does not
depend on the choice of $\psi$. Moreover, if $\widetilde L_0$ is
another fixed Lagrangian and $\widetilde\psi:[a,b]\to\Spl(2n,\R)$
is a continuous curve such that
$\gamma_2(t)=\widetilde\psi(t)(\widetilde L_0)$, then
\[\mu_{L_0}\big(t\mapsto\psi(t)^{-1}\gamma_1(t)\big)=\mu_{\widetilde
L_0}\big(t\mapsto\widetilde\psi(t)^{-1}\gamma_1(t)\big).\]
Likewise, if $\gamma_1$ and $\gamma_2$ are smooth curves, and
$t_0$ is an isolated non transversality instant for $\gamma_1$ and
$\gamma_2$, then the $L_0$-partial signatures of the curve
$t\mapsto\psi(t)^{-1}\gamma_1(t)$ at $t=t_0$ do not depend on the
choice of the Lagrangian $L_0$ and of the $L_0$-lifting $\psi$ of
$\gamma_2$.
\end{lem}
\begin{proof}
If $\psi_1,\psi_2:[a,b]\to\Spl(2n,\R)$ are any two continuous $L_0$-liftings of
$\gamma_2$, then clearly $\phi(t)=\psi_2(t)^{-1}\psi_1(t)$ is a continuous
curve in $\Spl(2n,\R)$ such that $\phi(t)(L_0)=L_0$ for all $t$.
Using part (c) of Lemma~\ref{thm:admuL0} we obtain:
\[\mu_{L_0}\big(t\mapsto\psi_2(t)^{-1}(\gamma_1(t))\big)\!=\!
\mu_{L_0}\big(t\mapsto\phi(t)\psi_1(t)^{-1}(\gamma_1(t))\big)
\!=\!\mu_{L_0}\big(t\mapsto\psi_1(t)^{-1}(\gamma_1(t))\big),\]
which proves the independence on the lifting.

The proof of the second part of the statement follows from the
symplectic invariance of the Maslov index (part~\eqref{itm:propmuL04}
of Lemma~\ref{thm:propmuL0}). Namely, choose a symplectomorphism
$h\in\Spl(2n,\R)$ such that $L_0=h(\widetilde L_0)$,  let $\widetilde
\psi:[a,b]\to\Spl(2n,\R)$ be a continuous  $\widetilde L_0$-lifting of
$\gamma_2(t)$, and set $\psi(t)=\widetilde
\psi(t)h^{-1}$, so that $\psi$ is an $L_0$-lifting of $\gamma_2(t)$. Then:
\[\begin{split}\mu_{L_0}\big(t\mapsto\psi(t)^{-1}&(\gamma_1(t))\big)=
\mu_{L_0}\big(t\mapsto h\widetilde\psi(t)^{-1}(\gamma_1(t))\big)\\&=
\mu_{h^{-1}(L_0)}\big(t\mapsto\widetilde\psi(t)^{-1}(\gamma_1(t))\big)=
\mu_{\widetilde L_0}\big(t\mapsto\widetilde\psi(t)^{-1}(\gamma_1(t))\big),\end{split}\]
which concludes the proof of the first part of the statement.
The proof of the last statement is totally analogous, and it uses
the invariance property of the partial signatures discussed in
Lemma~\ref{thm:compsympl}.
\end{proof}
We can now give the following:
\begin{defin}\label{thm:defmaslovindpairs}
The Maslov index $\mu(\gamma_1,\gamma_2)$ of a pair
$(\gamma_1,\gamma_2):[a,b]\to\Lambda$ of continuous curves in
$\Lambda$ is the $L_0$-Maslov index of the curve
$t\mapsto\psi(t)^{-1}\gamma_2(t)$, where $L_0\in\Lambda$ is any
fixed Lagrangian and $\psi:[a,b]\to\Spl(V,\omega)$ is any
$L_0$-lifting of $\gamma_2$. If $\gamma_1$ and $\gamma_2$ are
smooth and $t_0\in[a,b]$ is an isolated non transversality instant
for $\gamma_1$ and $\gamma_2$, then the partial signatures of the
pair $(\gamma_1,\gamma_2)$ at $t=t_0$, denoted by
$n^\pm_k(\gamma_1,\gamma_2,t_0)$, $\noo\pm_k(\gamma_1,\gamma_2,
t_0)$ and $\sigma_k(\gamma_1,\gamma_2,t_0)$, are defined as the
corresponding $L_0$-partial signatures of the symplectic path
$t\mapsto\psi(t)^{-1}(\gamma_2(t))$ at $t=t_0$, where $L_0$ is any
Lagrangian and $\psi:[a,b]\to\Spl(V,\omega)$ any smooth
$L_0$-lifting of  $\gamma_2$.
\end{defin}

The Maslov index of a pair $(\gamma_1,\gamma_2)$ of continuous curves
in the Lagrangian Grassmannian of a symplectic space $(V,\omega)$
is better described
in terms of the Maslov index of a single curve in a suitably ``doubled''
symplectic space $V^2$. Given $(V,\omega)$, denote by $(V^2,\omega^2)$ the
symplectic  space $V\oplus V$ endowed with the symplectic form
$\omega^2=\omega\oplus(-\omega)$:
\[\omega^2((v_1\oplus v_2),(v_3\oplus v_4))=\omega(v_1,v_3)-\omega(v_2,v_4).\]
Clearly, if $L_1,L_2$ are two Lagrangians in $(V,\omega)$, then $L_1\oplus L_2$ is  Lagrangian
in $(V^2,\omega^2)$, i.e., there is an injection $\Lambda(V,\omega)\times
\Lambda(V,\omega)$ into $\Lambda(V^2,\omega^2)$; moreover, the diagonal
$\Delta\subset V\oplus V$ is Lagrangian in $(V^2,\omega^2)$.

We will look at the Maslov index of curves in $\Lambda(V^2,\omega^2)$
computed relatively to the Lagrangian $\Delta$.
An easy calculations shows that, if $L_0\in\Lambda(V,\omega)$ and
$\gamma:[a,b]\to\Lambda(V,\omega)$ is continuous, then,
\begin{equation}\label{eq:doppiosingolo}
\mu_\Delta\big(t\mapsto\gamma(t)\oplus L_0\big)=
\mu_{L_0}\big(t\mapsto\gamma(t)\big).
\end{equation}
\begin{prop}\label{thm:Deltadoppio}
Let $\gamma_1,\gamma_2:[a,b]\to\Lambda(V,\omega)$ be
a pair of continuous curves;
then,
\begin{equation}\label{eq:Deltadoppio}
\mu(\gamma_1,\gamma_2)=\mu_\Delta(\gamma_1\oplus\gamma_2).
\end{equation}
Moreover, if $\gamma_1$ and $\gamma_2$ are smooth and if
$t_0\in[a,b]$ is an isolated non transversality instant for
$\gamma_1$ and $\gamma_2$, then the partial signatures of the pair
$(\gamma_1,\gamma_2)$ at $t=t_0$ coincide with the corresponding
$\Delta$-partial signatures of the curve
$t\mapsto\gamma_1(t)\oplus\gamma_2(t)\in\Lambda(V^2,\omega^2)$.
\end{prop}
\begin{proof}
Choose a Lagrangian $L_0\in\Lambda(V,\omega)$ and
a continuous $L_0$-lifting  $\phi:[a,b]\to\Spl(V,\omega)$  of $\gamma_2(t)$.
Clearly, $\phi(t)^{-1}\oplus\phi(t)^{-1}\in\Spl(V^2,\omega^2)$ and
$\phi(t)^{-1}\oplus\phi(t)^{-1}(\Delta)=\Delta$ for all $t$. Hence, using
part (c) of Corollary~\ref{thm:admuL0} and formula \eqref{eq:doppiosingolo},
we get:
\begin{multline*}\mu_\Delta(\gamma_1\oplus\gamma_2)= \mu_\Delta\big(t\mapsto\phi(t)
^{-1}(\gamma_1(t))
\oplus\phi(t)^{-1}(\gamma_2(t))\big)\\ =\mu_\Delta\big(t\mapsto
\phi(t)^{-1}(\gamma_1(t))\oplus L_0\big) =\mu_{L_0}\big(t\mapsto
\phi(t)^{-1}(\gamma_1(t))\big)=\mu(\gamma_1,\gamma_2). \end{multline*}
This proves the first statement of the thesis; the last statement
is proven similarly, using Lemma~\ref{thm:compsympl}.
\end{proof}
We will  collect below a few properties of the map $(\gamma_1,\gamma_2)
\mapsto\mu(\gamma_1,\gamma_2)$:
\begin{prop}\label{thm:propinddoppio}
The Maslov index for pairs of curves in $\Lambda$ satisfies the following:
\begin{enumerate}
\item\label{itm:propmudoppio1} the induced map $\mu:\pi(\Lambda)\oplus\pi(\Lambda)\to\Z$ is a groupoid homomorphism;
\item\label{itm:propmudoppio2} $\mu$ is anti-symmetric: $\mu(\gamma_1,\gamma_2)=-\mu(\gamma_2,\gamma_1)$;
\item\label{itm:propmudoppio3} if $\gamma:[a,b]\to\Lambda$ is a continuous curve
and $L_0\in\Lambda$ is fixed, then $\mu(\gamma,L_0)=\mu_{L_0}(\gamma)$.
\end{enumerate}
\end{prop}
\begin{proof}
\eqref{itm:propmudoppio1} is proven observing that
$\mu=\mu_{\Delta}\circ \mathfrak i$, where $\mathfrak i:\pi(\Lambda)\oplus
\pi(\Lambda)\to\pi(\Lambda\times\Lambda)$ is the groupoid homomorphism
induced by the immersion $\Lambda(V,\omega)\times\Lambda(V,\omega)\hookrightarrow
\Lambda(V^2,\omega^2)$.

\eqref{itm:propmudoppio2} follows from \eqref{itm:propmuL01}
and \eqref{itm:propmuL04} of Lemma~\ref{thm:propmuL0}, observing
that the map $V^2\ni x\oplus y\mapsto y\oplus x\in V^2$ is
an anti-symplectomorphism that preserves $\Delta$.

\eqref{itm:propmudoppio3} follows immediately from
Proposition~\ref{thm:Deltadoppio} and formula \eqref{eq:doppiosingolo}.
\end{proof}
The homotopy invariance of the Maslov index for pairs of Lagrangian
paths can be used to obtain a series of interesting facts, otherwise
not so evident, about
the Maslov index for single curves in $\Lambda$.
Here is an example:
\begin{lem}\label{thm:lemint1}
Let $\gamma_1,\gamma_2:[a,b]\to\Lambda$ be a pair of continuous
curves. The following equality holds:
\begin{equation}\label{eq:lemint1}
\mu_{\gamma_1(a)}(\gamma_2)-\mu_{\gamma_1(b)}(\gamma_2)=
\mu_{\gamma_2(b)}(\gamma_1)-\mu_{\gamma_2(a)}(\gamma_1).
\end{equation}
\end{lem}
\begin{proof}
Consider the continuous map \[[a,b]\times[a,b]\ni(s,t)\longmapsto H(s,t)=\gamma_1(s)\oplus
\gamma_2(t)\in\Lambda(V^2,\omega^2)\] and apply part \eqref{itm:propmuL03}
of Lemma~\ref{thm:propmuL0}. The conclusion follows easily from
Proposition~\ref{thm:propinddoppio}.
\end{proof}
\begin{cor}\label{thm:intfact2}
If $\gamma:[a,b]\to\Lambda$ is a continuous loop, then the value of
$\mu_{L_0}(\gamma)$ does not depend on the choice of $L_0$.
\end{cor}
\begin{proof}
Choose any two Lagrangians $L_0,L_1\in\Lambda$, and any continuous
curve $\gamma_2:[a,b]\to\Lambda$ with $\gamma_2(a)=L_0$ and
$\gamma_2(b)=L_1$. Set $\gamma_1=\gamma$ and apply Lemma~\ref{thm:lemint1}
to $\gamma_1$ and $\gamma_2$. Since $\gamma_1$ is a loop, the
left hand side of \eqref{eq:lemint1} vanishes, yielding:
\[\mu_{L_0}(\gamma)=\mu_{L_1}(\gamma),\]
which was to be proven.
\end{proof}

\subsection{On the Maslov triple and four-fold indexes}
\label{sub:triple}
Let us now discuss a different notion of Maslov index, originally
due to Kashiwara (see \cite{LioVer}), and further investigated by Cappell, Lee
and Miller in \cite[Section~8]{CapLeeMil}.

Assume that $(V,\omega)$ is a fixed (finite dimensional) symplectic
space; given three Lagrangians $L_1,L_2,L_3\in\Lambda(V,\omega)$, the {\em Maslov
triple index\/} $\tau_V(L_1,L_2,L_3)$ is defined as the signature of
the (symmetric bilinear form associated to the) quadratic form $Q:L_1\oplus L_2\oplus L_3\to\R$ given by:
\[Q(x_1,x_2,x_3)=\omega(x_1,x_2)+\omega(x_2,x_3)+\omega(x_3,x_1).\]
It is proven in \cite[Section~8]{CapLeeMil} that $\tau_V$ is the unique
integer valued map on $\Lambda\times\Lambda\times\Lambda$ satisfying the following properties:
\begin{itemize}
\item[\textbf{[P1]}] (skew symmetry) If $\sigma$ is a permutation of
the set $\{1,2,3\}$,
\[\tau_V(L_{\sigma(1)},L_{\sigma(2)},L_{\sigma(3)})=\textrm{sign}(\sigma)
\,\tau_V(L_1,L_2,L_3);\]
\item[\textbf{[P2]}] (symplectic additivity) given symplectic spaces
$(V,\omega)$, $(\widetilde V,\widetilde\omega)$, and Lagrangians
$L_1,L_2,L_3\in\Lambda(V,\omega)$, $\widetilde L_1, \widetilde L_2,
\widetilde L_3
\in\Lambda(\widetilde V,\widetilde \omega)$, then:
\[\tau_{V\oplus \widetilde V}(L_1\oplus \widetilde L_1,L_2\oplus \widetilde L_2,
L_3\oplus \widetilde L_3)=\tau_V(L_1,L_2,L_3)+\tau_{\widetilde V}(\widetilde L_1,
\widetilde L_2,\widetilde L_3);\]
\item[\textbf{[P3]}] (symplectic invariance) if $\phi:(V,\omega)\to(\widetilde
V,\widetilde\omega)$ is a symplectomorphism, then:
\[\tau_V(L_1,L_2,L_3)=\tau_{\widetilde V}(\phi(L_1),\phi(L_2),\phi(L_3));\]
\item[\textbf{[P4]}] (normalization) if $V=\R^2$ is endowed with the
canonical symplectic form\footnote{i.e., $\omega((x_1,y_1),(x_2,y_2))=x_1y_2-x_2y_1$}, and $L_1=\R(1,0)$, $L_2=\R(1,1)$, $L_3=\R(0,1)$,
then\[\tau_V(L_1,L_2,L_3)=1.\]
\end{itemize}
We will now proceed to a geometrical description of the triple index
$\tau_V$ using the notion of Maslov index for paths; we will introduce to this
aim a four-fold index, i.e. a   map
\[\mathfrak q:\Lambda\times\Lambda\times\Lambda\times\Lambda\longrightarrow\Z.\]
\begin{lem}\label{thm:perdefq}
Given four Lagrangians $L_0,L_1,L_0',L_1'\in\Lambda$ and
any continuous curve $\gamma:[a,b]\to\Lambda$ such that
$\gamma(a)=L_0'$ and $\gamma(b)=L_1'$, then the value of
the quantity $\mu_{L_1}(\gamma)-\mu_{L_0}(\gamma)$ does
{\em not\/} depend on the choice of $\gamma$.
\end{lem}
\begin{proof}
An easy application of Corollary~\ref{thm:intfact2}.
\end{proof}
An analogous result has been proven by Robbin and Salamon
for their half-integer valued Maslov index (see \cite[Theorem~3.5]{RobSal}).
We are now entitled to give the following:
\begin{defin}\label{thm:deffourfold}
Given four Lagrangians $L_0,L_1,L_0',L_1'\in\Lambda$, the
{\em four-fold Maslov index\/} $\mathfrak q(L_0,L_1;L_0',L_1')$ is
the integer number $\mu_{L_1}(\gamma)-\mu_{L_0}(\gamma)$, where
$\gamma:[a,b]\to\Lambda$ is any continuous curve with
$\gamma(a)=L_0'$ and $\gamma(b)=L_1'$.
\end{defin}
The four-fold Maslov index $\mathfrak q$, also known in the literature as
the {\em H\"ormander's index},  satisfies some
symmetries that resemble those satisfied by the curvature
tensor of a symmetric connection:
\begin{prop}\label{thm:symfrakq}
Let $L_0,L_1,L_0',L_1', L\in\Lambda$ be five Lagrangians.
The following  identities hold:
\begin{itemize}
\item[(a)] $\mathfrak q(L_0,L_1;L_0',L_1')=-\mathfrak q(L_1,L_0;L_0',L_1')$;
\item[(b)] $\mathfrak q(L_0,L_1;L_0',L_1')=-\mathfrak q(L_0,L_1;L_1',L_0')$;
\item[(c)] $\mathfrak q(L_0,L_1;L_0',L_1')=-\mathfrak q(L_0',L_1';L_0,L_1)$;%
\footnote{the curvature tensor $R$ of a symmetric connection satisfies
the identity \[R(x_0,x_1;x_0',x_1')=R(x_0',x_1';x_0,x_1),\] and in turn, such
symmetry leads to the Bianchi identity for $R$. In the case of the Maslov four-fold
index, the anti-symmetry  (c) leads to the   {\em cocycle identity}
\eqref{eq:cocycle}.}
\item[(d)] $\mathfrak q(L_0,L_1;L_0',L)+\mathfrak q(L_0,L_1;L,L_1')=
\mathfrak q(L_0,L_1;L_0',L_1')$.
\end{itemize}
\end{prop}
\begin{proof}
(a) and (b) are obvious by the definition of $\mathfrak q$, while
(d) is simply the additivity by concatenation of the Maslov index.
Part (c) follows easily from Lemma~\ref{thm:lemint1}.
\end{proof}
A whole series of different identities satisfied by
the four-fold index $\mathfrak q$ are easily obtained
by combining the equalities above; for instance:
\begin{equation}\label{eq:skewmiddle}
\mathfrak q(L_0,L_1;L_0',L_1')=-\mathfrak q(L_1',L_0';L_1,L_0).
\end{equation}

We can now establish the relation between the four-fold Maslov
index $\mathfrak q$ and the Maslov triple index $\tau_V$. Let us define
$\overline{\mathfrak q}:\Lambda\times\Lambda\times\Lambda\to\Z$ by:
\begin{equation}\label{eq:defqbarra}
\overline{\mathfrak q}(L_0,L_1,L_2):=\mathfrak q(L_0,L_1;L_2,L_0).
\end{equation}
Observe that the function $\mathfrak q$ is completely determined
by $\overline{\mathfrak q}$, because of the following identity:
\begin{multline}\label{eq:relqqbar}\mathfrak q(L_0,L_1;L_0',L_1')=\mathfrak q(L_0,L_1;L_0',L_0)+
\mathfrak q(L_0,L_1;L_0,L_1')\\  =\overline{\mathfrak q}(L_0,L_1,L_0')-\overline{\mathfrak q}
(L_0,L_1,L_1'),\end{multline}
however, the symmetries of the four-fold Maslov index are easier to
detect thanks to the clear geometrical meaning of $\mathfrak q$.
\begin{prop}\label{thm:q=tauV}
The map $\overline{\mathfrak q}$ coincides with the Maslov
triple index $\tau_V$.
\end{prop}
\begin{proof}
By  uniqueness, it suffices to prove that $\overline{\mathfrak q}$
satisfies the properties [P1], [P2], [P3] and [P4] above.
[P2] and [P3] are easily checked using respectively the
symplectic additivity and the symplectic invariance of the Maslov
index of paths. [P4] is also easily checked by an explicit calculation,
whose details are omitted. Property [P1], the skew-symmetry, is
the non obvious part of the statement; it suffices to
prove the two equalities
$\overline{\mathfrak q}(L_0,L_1,L_2)=-\overline{\mathfrak q}(L_0,L_2,L_1)$
and $\overline{\mathfrak q}(L_0,L_1,L_2)=-\overline{\mathfrak q}(L_1,L_0,L_2)$.
The first of the two equalities is obtained using \eqref{eq:skewmiddle}, while
the second is obtained as follows:
\begin{multline*}
\overline{\mathfrak q}(L_1,L_0,L_2)= \mathfrak q(L_1,L_0;L_2,L_1)\\=
\mathfrak q(L_1,L_0;L_2,L_0)+\mathfrak q(L_1,L_0;L_0,L_1) =\mathfrak q(L_1,L_0;L_2,L_0)=
-\overline{\mathfrak q}(L_0,L_1,L_2).\end{multline*}
This concludes the proof.
\end{proof}
Using Proposition~\ref{thm:symfrakq}, it is easy to check that
$\overline{\mathfrak q}$ satisfies the following cocycle
identity (see \cite[p.\ 163]{CapLeeMil}):
\begin{equation}\label{eq:cocycle}
\overline{\mathfrak q}(L_1,L_2,L_3)=\overline{\mathfrak q}(L_1,L_2,L_4)+
\overline{\mathfrak q}(L_2,L_3,L_4)+
\overline{\mathfrak q}(L_3,L_1,L_4).
\end{equation}
Let us conclude our discussion on the triple and the four-fold
index with the observation that it is possible to give
an alternative construction of the Maslov index for Lagrangian
paths using only the function $\mathfrak q$ (or $\overline{\mathfrak q}$).
Namely, assume that one is given two fixed Lagrangians
$L_0,L_1\in\Lambda$ and a continuous curve $\gamma:[a,b]\to\Lambda$
whose image is contained in $\Lambda_0(L_1)$, i.e., $\gamma(t)$ is transversal
to $L_1$ for all $t\in[a,b]$. In this case $\mu_{L_1}(\gamma)=0$,
and thus the quadruple
index $\mathfrak q(L_0,L_1;\gamma(a),\gamma(b))$ coincides
with the negative Maslov index $-\mu_{L_0}(\gamma)$. In the general case,
the interval $[a,b]$ admits a finite partition $t_0=a<t_1<\ldots<t_M=b$
such that $\gamma\big([t_{i-1},t_{i}]\big)$ is contained in
$\Lambda_0(L_i)$, for some $L_i\in\Lambda$ and $i\in\{1,\ldots,M\}$, and
by the concatenation additivity:
\begin{equation}\label{eq:muq}
\mu_{L_0}(\gamma)=-\sum_{i=1}^M\mathfrak q\big(L_0,L_i;\gamma(t_{i-1}),
\gamma(t_i)\big).
\end{equation}
Clearly, the choice of the partition $(t_i)_{i=1}^{M-1}$
and of the Lagrangians $(L_i)_{i=1}^M$ is not unique.
Using \eqref{eq:relqqbar} we get:
\begin{equation}\label{eq:mubarq}
\mu_{L_0}(\gamma)=\sum_{i=1}^M\Big[\overline{\mathfrak q}
\big(L_0,L_i,\gamma(t_{i})\big)-\overline{\mathfrak q}
\big(L_0,L_i,\gamma(t_{i-1})\big)\Big].
\end{equation}
\subsection{Maslov index of   symplectic paths}
Let $\Phi:[a,b]\to\Spl(V,\omega)$ be a continuous curve;
it is easy to see that, for each $t\in[a,b]$, the graph
$\Gr\big(\Phi(t)\big)$ is a Lagrangian subspace of
the symplectic space $(V^2,\omega^2)$ defined in
Subsection~\ref{sub:pairs}. We can therefore give the following:
\begin{defin}\label{thm:defmindsimplgr}
Given a continuous curve $\Phi$ in the symplectic group
$\Spl(V,\omega)$, the Maslov index $\iMaslov(\Phi)$ of $\Phi$ is the $\Delta$-Maslov
index of the curve $t\mapsto\Gr\big(\Phi(t)\big)\in\Lambda(V^2,\omega^2)$:
\[\iMaslov(\Phi):=\mu_\Delta\big(t\mapsto\Gr(\Phi(t))\big).\]
\end{defin}
The $\Delta$-Maslov index of a symplectic path $\Phi$ is also known in the literature
as the {\em Conley--Zehnder index\/} of $\Phi$ (see \cite{ConZeh, RobSal, SalZeh}).

Recall that each $L_0\in\Lambda(V,\omega)$ gives a smooth
map $\beta_{L_0}:\Spl(V,\omega)\to\Lambda(V,\omega)$ (defined
in \eqref{eq:defbetaL0}), and, with the help of the four-fold
Maslov index $\mathfrak q$, we can compare the
Maslov index of the curve $\Phi$ with the $L_0$-Maslov
index of the curve $\beta_{L_0}\circ\Phi:[a,b]\to\Lambda(V,\omega)$.
To this aim, we first give the following:
\begin{lem}\label{thm:comind}
Let $\Phi:[a,b]\to\Spl(V,\omega)$ be a continuous curve and
let $L_0,L_1,L_1'\in\Lambda(V,\omega)$ be fixed. Then:
\[\mu_{L_0}\big(\beta_{L_1}\circ\Phi\big)-\mu_{L_0}\big(\beta_{L_1'}\circ\Phi\big)=
\mathfrak q\big(L_1,L_1';\Phi(a)^{-1}(L_0),\Phi(b)^{-1}(L_0)\big).\]
\end{lem}
\begin{proof}
Using the Maslov index for pairs and the symplectic
invariance, we compute as follows:
\[\mu_{L_0}\big(\beta_{L_1}\circ\Phi\big)=\mu\big(\beta_{L_1}\circ\Phi,L_0\big)
=\mu\big(L_1,t\mapsto\Phi(t)^{-1}(L_0)\big)=-\mu_{L_1}\big(t\mapsto\Phi(t)^{-1}(L_0)
\big).
\]
Similarly,
\[\mu_{L_0}\big(\beta_{L_1'}\circ\Phi\big)=
-\mu_{L_1'}\big(t\mapsto\Phi(t)^{-1}(L_0)\big).\]
The conclusion follows easily from the  definition
of $\mathfrak q$.
\end{proof}

\begin{prop}
\label{thm:compind2}
Let $\Phi:[a,b]\to\Spl(V,\omega)$ be a continuous curve and $L_0,\ell_0\in\Lambda(V,\omega)$ be
fixed. Then:
\[\iMaslov(\Phi)+\mu_{L_0}\big(\beta_{\ell_0}\circ\Phi\big)=
\mathfrak q\big(\Delta,L_0\oplus\ell_0;\Gr\big(\Phi(a)^{-1}\big),
\Gr\big(\Phi(b)^{-1}\big)\big).\]
In particular, if $\Phi$ is a loop, then $\iMaslov(\Phi)=-\mu_{L_0}(\beta_{\ell_0}\circ\Phi)$.
\end{prop}
\begin{proof}
We compute:
\[\iMaslov(\Phi)=\mu_\Delta\big(t\mapsto(\mathrm{Id}\oplus\Phi(t))(\Delta)\big)
\]
and, using the properties of the Maslov index for pairs of curves,
\[\mu_{L_0}\big(\beta_{\ell_0}\circ\Phi\big)=-\mu_{\Delta}\big(t\mapsto L_0\oplus
\beta_{\ell_0}\circ\Phi(t)\big)=-\mu_\Delta\big(t\mapsto(\mathrm{Id}\oplus
\Phi(t))(L_0\oplus\ell_0)\big).\]
The result follows now easily applying Lemma~\ref{thm:comind} to
the curve $t\mapsto\mathrm{Id}\oplus\Phi(t)\in\Spl(V^2,\omega^2)$ and
to the Lagrangians $\Delta,L_0\oplus\ell_0\in\Lambda(V^2,\omega^2)$.
\end{proof}

\subsection{Spectral flow of affine paths}
\label{sub:spflaffine}
Let us now discuss an application of Proposition~\ref{thm:maslovpartsign}
that will be used to compute the spectral flow of affine paths
(Proposition~\ref{thm:surpresa}, Corollary~\ref{thm:affine}), and
in the  spectral index theorem
(Subsection~\ref{sub:spectral}):
\begin{example}\label{exa:nilpotente}
Let $V$ be a real, finite dimensional vector space, $g:V\times V\to\R$ a nondegenerate
symmetric bilinear form, and $T:V\to V$ a nilpotent ($T^n=0$) linear $g$-symmetric endomorphism (i.e.,
such that the bilinear form $g^T:=g(T\cdot,\cdot)$ is symmetric) of $V$.
For a subspace $W\subset V$, denote by $W^{\perp_g}$ the $g$-orthogonal
space to $W$, defined by \[W^{\perp_g}=\big\{v\in V:g(v,w)=0\ \text{for all}\ w\in W\big\}.\]
For all $k=1,\ldots,n$, define
a vector space $W_k=T^{k-1}\big(\Ker(T^k)\big)$ and a  bilinear
form $B_k:W_k\times W_k\to\R$ by:
\[B_k(a,b)=g(c,b),\]
where $a=T^{k-1}(c)$. Using the identities:\footnote{The second identity in
\eqref{eq:kerimab} holds for any linear operator $T$ on
any vector space $V$. The first one needs the assumption
that $T$ is $g$-symmetric and, in addition,
that $\Dim(V)<+\infty$, or, if $V$ is an infinite dimensional Hilbert space,
that $g$ is {\em strongly nondegenerate\/} on $V$ (i.e.,
realized by a self-adjoint isomorphism of $V$) and that $T$ is a Fredholm operator of index $0$.
Identities~\eqref{eq:kerimab}  will be used also in the proof
of Propositions~\ref{thm:surpresa} and \ref{thm:surpresa2} under these more general assumptions.}
\begin{equation}\label{eq:kerimab}
\Ker(T^\alpha)^{\perp_g}=\Imm(T^\alpha),\quad
\Imm(T^\alpha)\cap\Ker(T^\beta)=T^\alpha\big(\Ker(T^{\alpha+\beta})\big),
\quad\forall\,\alpha,\beta=1,\ldots,n,\end{equation}
it is not hard to see directly that $W_k\supseteq W_{k+1}$,
that $B_k$ is a well defined  symmetric bilinear form on $W_k$, and that
$W_{k+1}=\Ker(B_k)$ for all $k$; in particular, $B_n$ is nondegenerate.
We get to the same conclusions
indirectly  observing that the spaces $W_k$ and the bilinear forms $B_k$ can
be obtained using the construction of Section~\ref{sec:gensig} (Remark~\ref{thm:remuseless})
applied to the real-analytic path of symmetric bilinear forms
\[[-\varepsilon,\varepsilon]\ni\lambda\longmapsto
B(\lambda)=g^T-\lambda g\in\Bsym(V)\] at the isolated singularity $\lambda=0$.

Using the homotopy invariance of the Maslov index, we will now show that
the following identities hold:
\begin{equation}\label{eq:formulespettrale}
\begin{split}&\sum_{k\ge1}\big(n^-(B_{2k-1})+n^+(B_{2k})\big)={\noo-}(g^T)-n^-(g)=n^+(g)-n^+\big(g^T\big),\\
&\sum_{k\ge1}n^+(B_k)={\noo-}(g^T)-n^+(g)=n^-(g)-n^+\big(g^T\big),\\ &
\sum_{k\ge1}\sigma(B_{2k-1})=-\sigma(g).\end{split}\end{equation}
By Corollary~\ref{thm:centralinfdim},
\[\begin{split}&\sum\limits_{k\ge1}\big(n^-(B_{2k-1})+n^+(B_{2k})\big)=
\Dim(\Ker T)-\spfl\big(B,[-\varepsilon,0]\big),\\
&\sum\limits_{k\ge1}n^+(B_k)=\Dim(\Ker T)+\spfl\big(B,[0,\varepsilon]\big),\\
&\sum\limits_{k\ge1}\sigma(B_{2k-1})=
\spfl\big(B,[-\varepsilon,\varepsilon]\big).\end{split}\]

To prove the  equalities \eqref{eq:formulespettrale}, consider the two-parameter map
\[[0,1]\times[-\varepsilon,\varepsilon]\ni(r,\lambda)\longmapsto \ell(r,\lambda)=rg^T-\lambda g\in\Bsym(V);\]
observe that $B(\lambda)=\ell(1,\lambda)$, and that the bilinear forms
$\ell(r,-\varepsilon),\ell(r,\varepsilon)$ are nondegenerate for all $r\in[0,1]$, because
$0$ is the unique eigenvalue of $T$.
It follows  (see Figure~\ref{fig:g1g2g3g4g5}))
\[\begin{aligned}&&\spfl\big(B,[-\varepsilon,0]\big)=
\spfl\big(\ell(0,\cdot),[-\varepsilon,0]\big)+
\spfl\big(\ell(\cdot,0),[0,1]\big),
\\ && \spfl\big(B,[0,\varepsilon]\big)=-\spfl\big(\ell(\cdot,0),[0,1]\big)+
\spfl\big(\ell(0,\cdot),[0,\varepsilon]\big).
\end{aligned}\]
\begin{figure}
\begin{center}
\psfull \epsfig{file=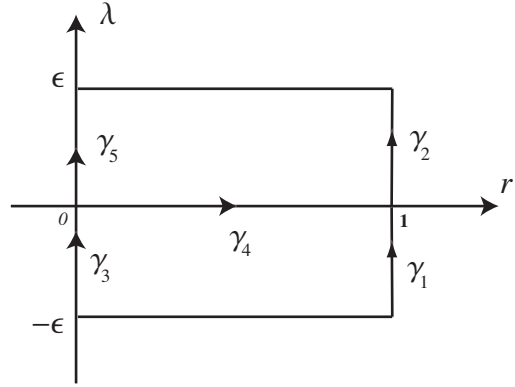} \caption{The curves of
symmetric bilinear forms used in Example~\ref{exa:nilpotente}. In
the picture, $\gamma_1$ is the curve
$[-\varepsilon,0]\ni\lambda\mapsto B(\lambda)$, $\gamma_2$ is the
curve $[0,\varepsilon]\ni\lambda\mapsto B(\lambda)$, $\gamma_3$ is
the curve $[-\varepsilon,0]\ni\lambda\mapsto \ell(0,\lambda)$,
$\gamma_4$ is the curve $[0,1]\ni r\mapsto \ell(r,0)$ and
$\gamma_5$ is the curve
$[0,\varepsilon]\ni\lambda\mapsto\ell(0,\lambda)$}\label{fig:g1g2g3g4g5}
\end{center}
\end{figure}
 An immediate computation gives:
\[\begin{split}&\spfl\big(\ell(0,\cdot),[-\varepsilon,0]\big)={\noo+}(0)-{\noo+}(\varepsilon g)=\Dim(V)-n^+(g)=n^-(g),\\ &
\spfl\big(\ell(\cdot,0),[0,1]\big)={\noo+}(g^T)-{\noo+}(0)={\noo+}(g^T)-\Dim(V)=-n^-(g^T),\\
&\spfl\big(\ell(0,\cdot),[0,\varepsilon]\big)={\noo+}(-\varepsilon g)-{\noo+}(0)=n^-(g)-\Dim(V)=-n^+(g),\end{split}\] from which
equalities
\eqref{eq:formulespettrale} follow easily.
\end{example}
\begin{prop}\label{thm:surpresa}
Let $\Hcal$ be a real, separable Hilbert space, $g:\Hcal\to \Hcal$ an invertible
self-adjoint linear operator and $T:\Hcal\to\Hcal$ a linear Fredholm operator
such that $gT=T^*g$. Then:
\begin{enumerate}
\item\label{itm:part1} $\Dim\big(\bigcup_{n\ge1}\Ker(T^n)\big)<+\infty$, i.e.,
there exists $n_0\ge0$ such that $\Ker(T^n)=\Ker(T^{n_0})$ for all
$n\ge n_0$ and $\Dim\big(\Ker(T^{n_0})\big)<+\infty$;
\item\label{itm:part2} the bilinear form $\langle g\cdot,\cdot\rangle$ is nondegenerate on $\Ker\big(T^{n_0}\big)$;
\item\label{itm:part3} there are no degeneracy instants $t\ne0$ near $0$ for the affine path
$t\mapsto {\LL}(t)=gT-t g$ of self-adjoint Fredholm operators, and for
$\varepsilon>0$ small enough:
\begin{equation}\label{eq:pp5}\begin{split}
&\spfl\big({\LL},[-\varepsilon,0]\big)=n^+\big(\mathcal B^T\big) -n^+(\mathcal
B)+\Dim(\Ker T)=\noo+\big(\mathcal B^T\big) -n^+(\mathcal
B),
\\ & \spfl\big({\LL},[0,\varepsilon]\big)=
n^-(\mathcal B)-n^+\big(\mathcal B^T\big)-\Dim(\Ker T)=n^-(\mathcal B)-\noo+\big(\mathcal B^T\big),\\ &
\spfl\big({\LL},[-\varepsilon,\varepsilon]\big)=-\sigma(\mathcal B),\end{split}\end{equation}
where $\mathcal B=\langle g\cdot,\cdot\rangle\vert_{\Ker(T^{n_0})}$ and
$\mathcal B^T=\langle gT\cdot,\cdot\rangle\vert_{\Ker(T^{n_0})}$.
\end{enumerate}
\end{prop}
\begin{proof}
We start observing that if $\Ker(T)=\{0\}$ then
the entire statement is trivial; we observe also that, since $g$ is an isomorphism,
the equality $gT=T^*g$ implies that $T$ is a Fredholm operator of index $0$.\footnote{%
Recall that the {\em index\/} $\ind(T)$ of a Fredholm operator $T$ is the integer
$\Dim\big(\Ker(T)\big)-\Codim\big(\Imm(T)\big)$; in particular, the index of an isomorphism
is zero. Composition of Fredholm operators is
Fredholm,  and the Fredholm index is additive by composition; moreover, $\ind(T^*)=-\ind(T)$.
The equality $gT=T^*g$ implies $\ind(T)=\ind(g)+\ind(T)=\ind(gT)=\ind(T^*g)=\ind(g)+\ind(T^*)=-\ind(T)$, hence
$\ind(T)=0$.}
If $\Ker(T)\ne\{0\}$, $t=0$ is  an isolated degeneracy instant of $\LL(t)$, because $0$ must be isolated
in the spectrum of the Fredholm operator $T$. As in Example~\ref{exa:nilpotente},
the spaces $W_k$ and the bilinear forms $B_k$ can be computed explicitly as:
\[W_k=T^{k-1}\big(\Ker(T^k)\big),\quad B_k(x,y)=g(z,y),\]
where $x,y\in W_k$ and $T^k(z)=x$. Since the path $\LL$ is real-analytic,
then  there exists $n_0\in\N$ such that $W_{n_0}\ne\{0\}$,
$B_{n_0}$ is nondegenerate on $W_{n_0}$, and $W_n=\{0\}$ for $n>n_0$.
The equality $W_n=\{0\}$ clearly implies $\Ker(T^n)=\Ker(T^{n-1})$, hence
$\Ker(T^n)=\Ker(T^{n_0})$ for all $n\ge n_0$, which proves part~\eqref{itm:part1}.
Consider now $u\in W_{n_0}$, $u\ne0$. If $x=T^{n_0-1}(u)\ne0$, then,
since $B_{n_0}$ is nondegenerate, there exists $y\in  W_{n_0}$ such that
$g(u,y)=B_{n_0}(x,y)\ne0$. On the other hand, if $u\in\Ker(T^{n_0-1})$
and $g(u,v)=0$ for all $v\in\Ker(T^{n_0})$, then, recalling the identities \eqref{eq:kerimab}:
\begin{multline*}u\in \Ker(T^{n_0-1})\cap \Ker(T^{n_0})^{\perp_g}=\\\Ker(T^{n_0-1})\cap\Imm(T^{n_0})=
T^{n_0}\big(\Ker(T^{2n_0-1})\big)=T^{n_0}\big(\Ker(T^{n_0})\big)=\{0\},\end{multline*}
which proves part~\eqref{itm:part2}. Finally, part~\eqref{itm:part3}
is obtained immediately from Example~\ref{exa:nilpotente} (formulas~\eqref{eq:formulespettrale})
and from Corollary~\ref{thm:centralinfdim}
by considering the restriction of $g$ and $T$ to the space $V=\Ker(T^{n_0})$
(which is clearly invariant by $T$), and observing that
the spectral flow of the path $t\mapsto gT-tg$ on $\Hcal$ coincides with the spectral
flow of its restriction\footnote{here, by restriction of a linear operator $L$
we mean the restriction of the corresponding bilinear form $\langle L\cdot,\cdot\rangle$.}
to $\Ker(T^{n_0})$ (see Example~\ref{thm:exarestriction}, Subsection~\ref{sub:restrictions}).
\end{proof}
We can now compute the spectral flow of arbitrary affine paths
of self-adjoint Fredholm operators:
\begin{cor}\label{thm:affine}
Let $A$ and $K$ be self-adjoint operators on $\Hcal$, with $A$ invertible
and $K$ compact. Assume that $\lambda_0\in\R\setminus\{0\}$ and that $\lambda_0^{-1}$
is in the spectrum of
$-A^{-1}K$, so that $\lambda_0$ is an isolated singularity of the
affine path $\lambda\mapsto {\Ss}(\lambda)=A+\lambda K$ in $\Fcalsa(\Hcal)$.
Then, for $\varepsilon>0$ small enough:
\begin{equation}\label{eq:affine}
\begin{split}&\spfl\big(\Ss,[\lambda_0-\varepsilon,\lambda_0]\big)=n^+(\mathcal B_2)-n^+(\mathcal
B_1)+\Dim\big(\Ker(A+\lambda_0K)\big)=\noo+(\mathcal B_2)-n^+(\mathcal
B_1),\\ &\spfl\big(\Ss,[\lambda_0,\lambda_0+\varepsilon]\big)=n^-(\mathcal B_1) -n^+(\mathcal
B_2)-\Dim\big(\Ker(A+\lambda_0K)\big)=n^-(\mathcal B_1) -\noo+(\mathcal
B_2),\\
&\spfl\big(\Ss,[\lambda_0-\varepsilon,\lambda_0+\varepsilon]\big)=-\sigma(\mathcal B_1),\end{split}
\end{equation}
where $\mathcal B_1=\langle A\cdot,\cdot\rangle\vert_{\Hcal_{\lambda_0}}$,
$\mathcal B_2=\langle (A+\lambda_0K)\cdot,\cdot\rangle\vert_{\Hcal_{\lambda_0}}$, and
$\mathcal H_{\lambda_0}$ is the finite dimensional subspace of $\Hcal$ given by
\[\mathcal H_{\lambda_0}=\bigcup_{n\ge1}\Ker\big(A^{-1}K+\tfrac1{\lambda_0}\mathrm{Id}\big)^n.\]
\end{cor}
\begin{proof}
 Assume first $\lambda_0>0$; by the cogredient invariance, for $\lambda$ near $\lambda_0$
the spectral
flow of ${\Ss}(\lambda)$ equals the spectral flow of
\[\tfrac1{\lambda}{\Ss}(\lambda)=(\tfrac1\lambda-\tfrac1{\lambda_0})A+A\big(A^{-1}K+\tfrac1{\lambda_0}\mathrm{Id}\big).\]
To obtain \eqref{eq:affine}, set $t=\tfrac1{\lambda_0}-\tfrac1{\lambda}$,
$g=A$, $T=A^{-1}K+\tfrac1{\lambda_0}\mathrm{Id}$
and apply Proposition~\ref{thm:surpresa} to this setup.
Observe that the identities \eqref{eq:kerimab} can now be used because
$T$ is Fredholm of index $0$ (it is a compact perturbation of the isomorphism
$\tfrac1{\lambda_0}\mathrm{Id}$).

The proof in the case $\lambda_0<0$ is obtained from the previous
case, replacing $\lambda_0$ with $-\lambda_0$ and $K$ with $-K$, observing
that the definition of the objects $\mathcal B_1$, $\mathcal B_2$ and
$\Hcal_{\lambda_0}$ are unchanged when both $\lambda_0$ and $K$ are
taken with the opposite sign.
\end{proof}

Using similar arguments, one proves the following
version of Proposition~\ref{thm:surpresa}, which is better suited
to study the case of Fredholm bilinear forms obtained from
unbounded operators paired with compact bilinear forms.
Recall that a densely defined linear operator
$T$ is said to be {\em discrete\/} if for some (hence for all)
$\lambda$ not in the spectrum $\mathfrak s(T)$ of $T$, the resolvent $(T-\lambda)^{-1}$
is a compact operator; the spectrum of a discrete operator $T$ is a discrete subset of $\C$,
and every element in the spectrum is an eigenvalue of $T$
(see \cite[Chapter~19]{DunSch}).
\begin{prop}\label{thm:surpresa2}
Let $\Hcal$ be a separable Hilbert space, $\Dcal_2\subset\Dcal_1\subset\Hcal$ be
dense algebraic subspaces, each of which is endowed with a Hilbert space
structure that makes each inclusion a bounded operator, and the inclusion of
$\Dcal_2$ into $\Hcal$ compact.
Let $G:\Hcal\to\Hcal$ be a self-adjoint isomorphism  and
$T:\Dcal_2 \to\Hcal$ a   $G$-symmetric Fredholm operator of index $0$.
Assume that the symmetric bilinear form $\langle GT\cdot,\cdot\rangle$
defined on $\Dcal_2$ admits a (bounded)
Fredholm extension $\widehat g^T $ to $\Dcal_1$. Then:
\begin{enumerate}
\item\label{itm:pp1}  $T$ is a discrete operator  (and therefore
its spectrum $\mathfrak s(T)$ consists of eigenvalues of $T$).
\end{enumerate}
Assume further that the following regularity property holds:
if $x\in\Dcal_1$ and $z\in\Hcal$ are
such that $\widehat g^T(x,y)=\langle Gz,y\rangle$ for all  $y\in\Dcal_2$, then
$x\in\Dcal_2$ (and necessarily $Tx=z$). Then:
\begin{enumerate}
\addtocounter{enumi}{1}
\item\label{itm:pp2} for all $\lambda\in\mathfrak s(T)\cap\R$, the generalized
eigenspace $\Hcal_\lambda=\bigcup_{n\ge1}\Ker\big((T-\lambda)^n\big)$ is finite dimensional,
and the restriction of $\widehat g=\langle G\cdot,\cdot\rangle$ to
$\Hcal_\lambda$ is nondegenerate;
\item\label{itm:pp3} for any compact interval $[a,b]$, the spectral flow of the path
of self-adjoint Fredholm operators
$[a,b]\ni\lambda\mapsto {\Ss}(\lambda)=\widehat g^T -\lambda\widehat g\in\Fcalsa(\Dcal_1)$
is given by:
\begin{equation}\label{eq:pp4}
\spfl\big([{\Ss},[a,b]\big)=n^-\big(\widehat g\vert_{\Hcal_a}\big)-{\noo+}\big(\widehat g_a^T\vert_{\Hcal_a}\big)
- \!\!\!\!\sum_{\lambda_0\in\mathfrak s(T)
\cap\left]a,b\right[}\!\!\!\sigma\big(\widehat g\vert_{\Hcal_{\lambda_0}}\big)
+{\noo+}\big(\widehat g_b^T\vert_{\Hcal_b}\big)-
n^+\big(\widehat g\vert_{\Hcal_b}\big)
\end{equation}
where $\widehat g_\lambda^T=\big\langle G(T-\lambda)\cdot,\cdot\big\rangle$.
\end{enumerate}
\end{prop}
\begin{proof}
Let $\lambda$ be an element\footnote{for the sake of precision,
the argument presented here works for only for {\em real\/} $\lambda$'s.
For the complex case, one needs to consider the complexification of
$\Hcal$ endowed with the inner product given by the sesquilinear extension
of $\langle\cdot,\cdot\rangle$, and the complex linear extensions of $T$ and of $G$.}  in
$\mathfrak s(T)$ and consider the Fredholm self-adjoint operator $S_\lambda=G(T-\lambda):\Dcal_2\to\Hcal$;
to prove that $T$ is discrete, observe that $(T-\lambda)^{-1}=S_\lambda^{-1}G$ and that
$G$ is a compact operator on $\Dcal_2$, due to the fact that the inclusion
of $\Dcal_2$ into $\Hcal$ is compact. This proves \eqref{itm:pp1}.

Parts~\eqref{itm:pp2} and \eqref{itm:pp3} are now proven repeating {\em verbatim\/}
the proof of Proposition~\ref{thm:surpresa}, keeping in mind the following:
\begin{itemize}
\item[(a)] the regularity property implies that $\Ker( {\Ss}(\lambda))=\Ker(T-\lambda)$
for all $\lambda\in[a,b]$, and that, for all $\lambda_0\in\mathfrak s(T)\cap[a,b]$, the
spaces $W_k$ obtained from the partial signatures construction for the curve
${\Ss}(\lambda)$ at the degeneracy
instant $\lambda_0$ are given  by
$(T-\lambda_0)^{k-1}\big(\Ker(T-\lambda_0)^k\big)$ for all $k\ge1$;
\item[(b)] the identities \eqref{eq:kerimab} can be used in this context
thanks to the assumption that $T$ is Fredholm of index $0$ and that $G$ is a self-adjoint isomorphism.
\end{itemize}
Formula \eqref{eq:pp4} is obtained easily from \eqref{eq:pp5} using the
additivity by concatenation of the spectral flow.
\end{proof}

\end{section}

\begin{section}{Semi-Riemannian geodesics}
\label{sec:geodesics}
As an application of the theory discussed, we will now
describe how the partial signature method can be
applied to the study of the Maslov index of a semi-Riemannian
geodesic.

The semi-Riemannian geodesic problem is a central example of strongly indefinite
variational problem to which some recent extension of the
classical Morse theory (see~\cite{abbomajer}), as well as of
the bifurcation theory (see~\cite{FitzPejsaRecht, geobif}),
can be applied to obtain global geometrical results.
In these theories, the Maslov index of a geodesic plays the role
of  a generalized Morse index, and it is an essential
point to understand how to compute it in terms of
the conjugate points along the geodesic.
In the Riemannian case, the well-known Morse index
theorem gives the equality between the number of conjugate
points along the geodesic and the Morse index of the geodesic
action functional, which is also equal to the
number of negative eigenvalues of the Jacobi differential
operator. We will establish a similar result for
the semi-Riemannian case, where suitable definitions
of ``Morse index'', ``number of conjugate points''
and ``number of negative eigenvalues'' have
to be introduced. More precisely, we will
define the following notions associated to a semi-Riemannian
geodesic:
\begin{itemize}

\item  {\em Maslov index}, as an appropriate count
of the conjugate points (Subsection~\ref{sub:conjpointsgeo});

\item {\em generalized Morse index\/}, as
the spectral flow of the path of index forms along
the geodesic (Subsection~\ref{sub:generalizedMorse});

\item {\em spectral index}, as an appropriate count of
the nonpositive eigenvalues of the Jacobi differential
operator (Subsection~\ref{sub:spectral}).

\end{itemize}
We will prove the equality of the three integer numbers
in the general case (Theorem~\ref{thm:tuttouguale}), extending the results of
\cite{CPPT, asian, Hel1, pacific, topology}
to the case of degenerate endpoints.

Two different notions of index are used
in the classical literature to count conjugate points.
The {\em geometric index}, i.e.,  the sum of the multiplicities
of the conjugate points, which is not a very meaningful notion
outside  the Riemannian or the causal Lorentzian context,
and the {\em conjugate index} (as in \cite{Hel1}, or
{\em focal index\/} in \cite{pacific}),
that coincides with the   Maslov index in the case
of nondegenerate (see Subsection~\ref{sub:conjpointsgeo})
conjugate points. In spite of the fact that such nondegeneracy
holds generically, the conjugate index is not
a good measure of the conjugate points (see \cite[\S5.4]{pacific} and
Subsection~\ref{sub:counterexample} below).

\subsection{Conjugate points and Maslov index}\label{sub:conjpointsgeo}
Let us recall briefly the definition of Maslov index
for a semi-Riemannian geodesic; the background material
for this section can be found in
references~\cite{Hel1, pacific, topology}.

Let $(M,g)$ be an $n$-dimensional semi-Riemannian manifold,
$\nabla$ the covariant derivative of the Levi--Civita
connection of $g$ and $R$ its curvature tensor, chosen with the
sign convention: $R(X,Y)=[\nabla_X,\nabla_Y]-\nabla_{[X,Y]}$.

Given a geodesic $\theta:[0,1]\to M$, the Jacobi equation
along $\theta$ is the second order linear equation
$V''=R(\dot\theta,V)\dot\theta$ for vector fields
$V$ along $\theta$; here prime means covariant differentiation
along $\theta$. Solutions of the Jacobi equation
are called Jacobi fields.
Let us recall that $t_0\in\left]0,1\right]$ is sad to be
a {\em conjugate instant along $\theta$} if there exists
a non zero Jacobi field $V$
such that $V(0)=V(t_0)=0$. The {\em multiplicity\/}
$\mul(t_0)$ of a conjugate instant $t_0$ is defined to be
the dimension of the vector space of all Jacobi fields
$V$ satisfying $V(0)=V(t_0)=0$; for all
conjugate instant $t_0$, $\mul(t_0)\le n-1$.

By a parallel trivialization of the tangent
bundle $TM$ along $\theta$ (or of the normal
bundle $\dot\theta^\perp$ in the non lightlike case), then the metric $g$
can be seen as a constant nondegenerate bilinear form
on $\R^n$, and the Jacobi equation becomes simply
$V''=RV$, where now the prime symbol denotes the standard derivative
of $\R^n$-valued maps, and $R(t)$ is a smooth curve of
$g$-symmetric endomorphisms of $\R^n$.
We will implicitly identify vector fields along $\theta$
with $\R^n$-valued maps via such trivialization.
Let us consider the flow of the Jacobi equation, which is the smooth
curve of isomorphisms \[\Phi_t:\R^{2n}\to\R^{2n}\]
defined by:
\[\Phi_t\big(V(0),V'(0)\big)=\big(V(t),V'(t)\big),\]
for all solutions $V$ of $V''(t)=R(t)V(t)$.
An immediate calculation shows that $\Phi_t$ preserves the symplectic form
$\omega_g\big((v_1,v_2),(w_1,w_2)\big)=g(v_1,w_2)-g(v_2,w_1)$,
hence we get a smooth curve in the Lie group
$\Spl(\R^{2n},\omega_g)$.

Setting $L_0=\{0\}\oplus\R^n$, which is Lagrangian
relatively to $\omega_g$, we get a smooth
curve $\gamma(t)=\Phi_t(L_0)$ in the Lagrangian
Grassmannian of $(\R^{2n},\omega_g)$.
Conjugate points along $\theta$ correspond to
instants $t_0$ at which $\gamma(t_0)$ is not transversal
to $L_0$, and the Maslov index of $\theta$, denoted by
$\iMaslov(\theta)$, is defined
to be the sum:
\begin{equation}\label{eq:defMaslovindexgeo}
\iMaslov(\theta)=\mu_{L_0}(\gamma)+n^-(g),
\end{equation}
where $\mu_{L_0}(\gamma)$ is the
$L_0$-Maslov index $\mu_{L_0}$ of the curve $\gamma$, as
defined in subsection~\ref{sub:maslovsymplpath}.
Such definition does not depend on the choice of a parallel
trivialization of $TM$ along $\theta$; a proof of this fact
in the case that the final  endpoint $\theta(1)$ is not
conjugate is proven in \cite{topology}, while for the general  case the
proof will be done in Lemma~\ref{thm:nondiptriv} below.

\begin{rem}\label{thm:rempuntoiniziale}
It should be observed
that the curve $\gamma$ obtained by the above construction from the Jacobi
equation along a semi-Riemannian geodesic  $\theta:[0,1]\to M$ is not entirely
contained in the Maslov cycle. This is due to the
well known fact that, although $\gamma$ always intersects
the Maslov cycle at the initial instant $t=0$ (i.e., $t=0$ is always a
conjugate instant in a trivial sense), there are no conjugate instants
in $\left]0, \varepsilon\right]$ for $\varepsilon>0$ small
enough (see for instance \cite[Proposition~2.7]{pacific}).
Based on this observation, the Maslov index of a semi-Riemannian
geodesic $\theta$ whose final endpoint $\theta(1)$ is not conjugate
had been defined in references \cite{Hel1, pacific, topology, fechado} as the Maslov index
of the restriction $\gamma\vert_{[\varepsilon,1]}$, in order to
exclude the contribution of the initial conjugate instant.
Such contribution can be computed easily using the theory
below, and it is equal to $-n^-(g)$ (see formula~\ref{eq:contributoiniziale}),
i.e., for $\varepsilon>0$ small enough:
\begin{equation}\label{eq:laltra}
\mu_{L_0}(\gamma)=\mu_{L_0}\big(\gamma\vert_{[\varepsilon,1]}\big)-n^-(g),\end{equation}
which shows that the definition of $\iMaslov(\theta)$ given
in \eqref{eq:defMaslovindexgeo} is consistent with that of
references~\cite{Hel1, pacific, topology, fechado}.
Observe that in the Riemannian case $\iMaslov(\theta)=\mu_{L_0}(\gamma)$
\end{rem}

It is well known that conjugate points along $\theta$ can
accumulate away from the instant $t=0$ (see \cite{fechado})
unless $g$ is positive (or negative!) definite.
Nevertheless, when $(M,g)$ is real-analytic, then so is
also every solution $V$ of the Jacobi equation along
$\theta$ and also the curve $\gamma$ above, and this
implies that there is only a finite number of conjugate points.

Let $\theta(t_0)$, $t_0\in\left]0,1\right]$, be a conjugate point
along $\theta$ and define:
\[\mathbb J[t_0]=\big\{J(t_0):J\in\mathbb J\big\},\]
where:
\begin{equation}\label{eq:defbbJ}
\mathbb J=\big\{J:\text{$J$ is a Jacobi field along $\theta$, with $J(0)=0$}\big\};
\end{equation}
observe that $\mathbb J$ is an $n$-dimensional vector space.
Then, since $\gamma(t_0)$ is conjugate, $\mathbb J[t_0]\ne\R^n$,
and the {\em signature\/} $\sgn(t_0)$ of the conjugate point
$\theta(t_0)$ is defined to be the signature of the restriction
of $g$ to $\mathbb J[t_0]^\perp$, where now $\perp$ denotes orthogonality
relative to $g$; an easy argument shows that the following equality holds:
\[\mathbb J[t_0]^\perp=\big\{J'(t_0):J\in\mathbb J,\ J(t_0)=0\big\}.\]
It is also easy to see that, for all conjugate instant $t_0$ along
$\theta$:
\[\mul(t_0)=\Dim\big(\mathbb J[t_0]^\perp\big)=\Codim\big(\mathbb J[t_0]\big).\]
When the restriction of $g$ to $\mathbb J[t_0]^\perp$ is
nondegenerate, then $\theta(t_0)$ is said to be a {\em
nondegenerate conjugate point}. Nondegenerate conjugate points
correspond to transversal intersections of the curve $\gamma$ with
the Maslov cycle ${\LLL}$; if all the conjugate points along
$\theta$ are nondegenerate, and if the final instant $t=1$ is not
conjugate, then the Maslov index $\iMaslov(\theta)$ equals the sum
of the signatures of all conjugate instants in $\left]0,1\right[$
along $\theta$. Using the theory developed in the present paper we
are now able to compute the Maslov index of any geodesic without
any nondegeneracy assumption and any assumption on the final
instant $t=1$.

\subsection{Partial signatures at a conjugate instant and Maslov index}
\label{sub:partsigconjinst}
It is natural to define   $n_k^-(\theta,t_0)$ ($\noo-_k(\theta,t_0)$),
$n_k^+(\theta,t_0)$ ($\noo+_k(\theta,t_0)$)
and $\sigma_k(\theta,t_0)$   respectively as the
$k$-th partial (extended) index, (extended) coindex and signature of the curve
$\gamma$ at $t_0$. For the sake of precision, we must show that these quantities
do not depend on the choice of a trivialization of $TM$:
\begin{lem}\label{thm:nondiptriv}
Let $\gamma,\widetilde\gamma:[0,1]\to\Lambda$ be curves   associated
to the semi-Riemannian geodesic $\theta$ by two different trivializations
of the tangent bundle $TM$ along $\theta$.
Then, for all conjugate instant $t_0$, the partial signatures
of $\gamma$ and $\widetilde\gamma$ at $t_0$ coincide.
\end{lem}
\begin{proof}
The conclusion follows easily from Lemma~\ref{thm:compsympl}, observing that
the curves $\gamma$ and $\widetilde\gamma$ are related by the formula
$\widetilde\gamma(t)=\varsigma\big(\gamma(t)\big)$, where
$\varsigma$ is a fixed symplectomorphism of $(\R^{2n},\omega_g)$.
More precisely, if $T:\R^n\to\R^n$ is the isomorphism
relating the two different trivializations of $TM$ along
$\theta$, then $\varsigma$ is given by:
\[\varsigma(x,y)=(Tx,g^{-1}{T^*}^{-1}y).\qedhere\]
\end{proof}
A method for computing the partial signatures of a geodesic at
a conjugate instant will be given in Subsection~\ref{sub:computation}.
Using Proposition~\ref{thm:maslovpartsign} and formulas \eqref{eq:defMaslovindexgeo},
\eqref{eq:laltra} we obtain immediately:

\begin{prop}\label{thm:maslovgeo}
If $(M,g)$ is a real-analytic semi-Riemannian manifold and
$\theta:[0,1]\to M$ is a geodesic in $M$, then  the Maslov
index of $\theta$ is given by:
\[\iMaslov(\theta)=  \sum_{t_0\ \text{conjugate instant in}\
\left]0,1\right[}\left[
\sum_{k\ge1}\sigma_{2k-1}(\theta,t_0)\right]+\sum_{k\ge1}\Big[n^+_{2k-1}(\theta,1)
+n^-_{2k}(\theta,1)\Big].\]
In particular,
the definition of Maslov index of $\theta$ does not depend on the choice
of the trivialization of $TM$ along $\theta$.\qed
\end{prop}
Observe that if $g$ is Riemannian, i.e., positive definite, then
for all conjugate instant $t_0$ along $\theta$, $B_1(\theta,t_0)$ is positive definite,
hence its signature coincides with the multiplicity of $t_0$ as a conjugate instant.
It follows that $B_k(\theta,t_0)=0$ for all $k\ge2$, and
the Maslov index of $\theta$ is equal to the sum of the multiplicities of all
the conjugate instants along $\theta$ in $\left]0,1\right]$.
Likewise, for any semi-Riemannian metric $g$, $B_1(\theta,0)$ coincides
with the metric $g$ on $T_{\theta(0)}M$, which is nondegenerate hence
$B_k(\theta,0)=0$ for all $k>1$, and
the contribution to the Maslov index $\mu_{L_0}(\gamma)$ given by the
initial instant $0$ can be computed from formula~\eqref{eq:maslovpartsign}:
\begin{equation}\label{eq:contributoiniziale}
n^+\big(B_1(\theta,0)\big)-\Dim\big(T_{\theta(0)}M\big)=
n^+(g)-n=-n^-(g).\end{equation}

\subsection{Computation of the partial signatures at a conjugate instant}
\label{sub:computation}

We will now give an operational method for computing
the partial signatures of a geodesic $\theta:[0,1]\to M$ at
a conjugate instant $t_0\in\left]0,1\right]$:
\begin{prop}\label{thm:operational}
Assume that the map $\mathbb J\ni J\mapsto J'(t_0)\in\R^n$ is injective (hence an isomorphism);
then the partial signatures of $\theta$ at $t_0$ coincide
with the partial signatures of the curve $\mathcal B$ of symmetric bilinear forms
on $\mathbb J$ given by:
\begin{equation}\label{eq:formula1}
\mathcal B_t(J_1,J_2)=g\big(J_1(t),J_2'(t)\big).
\end{equation}
Similarly, if $T:\R^n\to\R^n$ is a $g$-symmetric isomorphism of $\R^n$ such that the
map $\mathbb J\ni J\mapsto J(t_0)-TJ'(t_0)\in\R^n$ is injective,
then the partial signatures of $\theta$ at $t_0$ coincide
with the partial signatures of the curve $\mathcal B$ of symmetric bilinear forms
on $\mathbb J$ given by:
\begin{equation}\label{eq:formula2}
\mathcal B_t(J_1,J_2)=g\big(J_1(t),J_2'(t)-T^{-1}J_2(t)\big).
\end{equation}
\end{prop}
\begin{proof}
The condition that the map $S:\mathbb J\to\R^n$, $S(J)=J'(t_0)$,
be an isomorphism is equivalent to the transversality of the Lagrangian
\[\gamma(t_0)=\Phi_{t_0}(L_0)=\big\{\big(J(t_0),J'(t_0)\big):J\in \mathbb J\big\}\]
(recall that $\Phi$ has been defined in Subsection~\ref{sub:conjpointsgeo}
as the flow of the Morse--Sturm equation $V''=RV$)
with the Lagrangian $L_1=\R^n\oplus\{0\}$. Let us use the chart $\varphi_{L_0,L_1}$
around $\gamma(t_0)$ (Subsection~\ref{sub:charts}),
and let us identify the space $\mathbb J$ with $L_0$ via the map:
\begin{equation}\label{eq:identL0J}L_0\ni v\longmapsto J_v\in\mathbb J\end{equation}
where $J_v$ is the unique Jacobi field in $\mathbb J$
determined by the initial  condition $J'(0)=v$.
A straightforward calculation
shows that  $\varphi_{L_0,L_1}\circ\gamma$ is given by:
\begin{equation}\label{eq:formula1bis}
\varphi_{L_0,L_1}\circ\gamma(t)(v_1,v_2)=g\big(J_{S^{-1}(v_1)}(t_0),v_2\big),\quad
t\sim t_0.\end{equation}
The conclusion follows from Proposition~\ref{thm:invcogredient},
observing that \eqref{eq:formula1} and \eqref{eq:formula1bis}
are cogredient:
\[\mathcal B_t(v_1,v_2)=\varphi_{L_0,L_1}\circ\gamma(t)\big(Sv_1,Sv_2\big).\]

Similarly, the proof of the second part of the statement reduces
to a straightforward direct calculation of $\varphi_{L_0,\widetilde L_1}\circ\gamma$,
where  $\widetilde L_1=\big\{(Tw,w):w\in\R^n\big\}=\Gr(T^{-1})$. Note that
$\widetilde L_1$ is Lagrangian because $T$ is $g$-symmetric, transversal
to $L_0$ because $T$ is invertible, and transversal to $\gamma(t_0)$ due to
the assumption that $\mathbb J\ni J\mapsto J(t_0)-TJ'(t_0)\in\R^n$ is injective.
\end{proof}
Clearly, different choices of the isomorphism $T$ as
in the assumptions of Proposition~\ref{thm:operational} produce
cogredient curves of bilinear forms as in \eqref{eq:formula2}, so that
the partial signatures of \eqref{eq:formula2} do not depend on
the choice of the isomorphism $T$. Observe that the set of
isomorphisms $T$ as in the assumptions of Proposition~\ref{thm:operational}
is diffeomorphic to  $\Lambda_0\big(\gamma(t_0)\big)\cap
\Lambda_0(L_0)\cap\Lambda_0(L_1)$, which is a  dense  open subset of
$\Lambda$ by Baire's theorem.

Using Proposition~\ref{thm:operational}, the computation
of the Maslov index of a geodesic reduces to simple
computations involving the curvature tensor and
its derivatives. For instance, in the case of a simple
conjugate instant $t_0$ (i.e., a conjugate instant of multiplicity
one), if $J\in\mathbb J$ is a nontrivial Jacobi field
vanishing at $t_0$, then the derivatives of the
map $h(t)= g\big(J(t),J'(t)\big)$ at $t=t_0$
are given by:
\begin{multline*}h'(t_0)=g\big(J'(t_0),J'(t_0)\big),\  h''(t_0)=g\big(J'(t_0),R(t_0)J'(t_0)\big), \\
h^{(3)}(t_0)=g\big(J'(t_0),R'(t_0)J'(t_0)\big),
\ \dots,\ h^{(k)}(t_0)=g\big(J'(t_0),R^{(k-2)}(t_0)J'(t_0)\big),\dots\end{multline*}

From Proposition~\ref{thm:operational} we also get the following result:
\begin{cor}\label{thm:Wkconjinst}
 Let $t_0\in\left]0,1\right]$ be
a conjugate instant along $\theta$ and let $J_0\in\mathbb J$ be
such that $J_0(t_0)=0$; set $v_0=J_0'(0)$. Then, $v_0\in
B_k(\mathcal B,t_0)$ ($\mathcal B_t$ given in \eqref{eq:formula2})
if and only if there exists a smooth curve
$v:\left]t_0-\varepsilon, t_0+\varepsilon\right[\to
T_{\theta(0)}M$ with $v(t_0)=v_0$ such that the map
$\left]t_0-\varepsilon, t_0+\varepsilon\right[\ni t\mapsto
J_{v(t)}(t)$ has a zero of order greater than or equal to $k$ at
$t=t_0$.
\end{cor}

\begin{proof}
Choose a $g$-symmetric linear operator $T:\R^n\to\R^n $ as in the
second part of Proposition~\ref{thm:operational} and consider the
corresponding path of bilinear forms $\mathcal B_t$ as in
\eqref{eq:formula2}. Observe that when $J$ runs in $\mathbb J$,
$t\mapsto J'(t)-T^{-1}J(t)$ is a smooth map that takes arbitrary
values in $\R^n$ at $t=t_0$. Using this fact, an easy induction
argument on $k$ shows that a curve $t\mapsto J_t$ in $\mathbb J$
is a root function of order greater than or equal to $k$ for
$\mathcal B_t$ at $t=t_0$ if and only if the map $t\mapsto J_t(t)$
has a zero of order greater than or equal to $k$ at $t=t_0$. This
concludes the proof.
\end{proof}
\subsection{A generalized Morse index}
\label{sub:generalizedMorse}

For all $t\in\left]0,1\right]$, define  $H_t$ as the Sobolev
space $H^1_0\big([0,t],\R^n)$; let $S^\theta_t$ be the bounded
symmetric bilinear form on $H_t$ given by:
\begin{equation}\label{eq:defIt}
S^\theta_t(V,W)=\int_0^t\Big[g\big(V'(s),W'(s)\big)+g\big(R(s)V(s),W(s)\big)\Big]\,\mathrm ds,
\end{equation}
where $g$ is a nondegenerate symmetric bilinear form on $\R^n$
and $t\mapsto R(t)$ is a smooth curve of $g$-symmetric
endomorphisms of $\R^n$. The objects $g$ and $R$ are obtained
respectively from the metric $g$ and the curvature tensor $R$
of $g$ via a parallel trivialization of $TM$ along $\theta$.
By identifying the spaces $H_t$ with $H_1$ via the map
$h_t:H_t\to H_1$, $(h_tV)(s)=V(ts)$, we get a smooth
curve of bounded symmetric bilinear forms $\bar S^\theta_t$ on $H_1$ obtained
by the push-forward of $S^\theta_t$ by $h_t$. More explicitly:
\[\bar S^\theta_t(V,W)=\int_0^1\Big[\tfrac1t g\big(V'(s),W'(s)\big)+
tg\big(R(ts)V(s),W(s)\big)\Big]\,\mathrm ds,
\quad V,W\in H_1.\]
The kernel $\mathcal N_t$ of $\bar S^\theta_t$, which is clearly
given by the image of $\Ker(S^\theta_t)$ by $h_t$, consists of smooth
vector fields $V$ on $[0,1]$ such that $V(0)=V(1)=0$ and
satisfying the linear equation:
\[V''(s)=t^2R(ts)V(s),\quad s\in[0,1].\]
The map $\left]0,1\right]\ni t\mapsto \bar S^\theta_t$ is a smooth
map of Fredholm bounded symmetric bilinear forms on $H_1$, and the map:
$\mathcal S^\theta_t:=t\bar S^\theta_t$:
\begin{equation}\label{eq:defcalsthetat}
\mathcal S^\theta_t(V,W)=\int_0^1\Big[g\big(V'(s),W'(s)\big)+t^2g\big(R(ts)V(s),W(s)\big)\Big]
\,\mathrm ds,
\end{equation}
admits a real-analytic extension to $t=0$ obtained by setting:
\[\mathcal S^\theta_0(V,W)=\int_0^1g\big(V'(s),W'(s)\big)\,\mathrm ds.\]
It is easy to prove that, for $\varepsilon>0$ small enough,
$\mathcal S^\theta_t$ is nondegenerate for all
$t\in[0,\varepsilon]$, hence:
\[\Ker\big(\mathcal S^\theta_0\big)=\{0\}\quad\text{and}\quad
\spfl\big(\mathcal S^\theta,[0,1]\big)=\spfl\big(\mathcal
S^\theta,[\varepsilon,1]\big)=\spfl\big(\bar
S^\theta,[\varepsilon,1]\big).\]
\begin{defin}\label{thm:defgenMorseindex}
The {\em generalized Morse index\/} of $\theta$ is defined as:
\[\iMorse(\theta):=\Dim\big(\Ker(\mathcal S^\theta_1)\big)-\spfl\big(\mathcal S^\theta,[0,1]\big).\]
\end{defin}
Recalling Remark~\ref{thm:defspectralMorse}, if $(M,g)$ is Riemannian,
then the generalized Morse index of $\theta$ coincides with the
{\em extended Morse index\/} of the geodesic action functional
at the critical point $\theta$.

We observe here that, for all $t\in\left]0,1\right]$, $\mathcal
S^\theta_t$ is realized by a compact perturbation of $\mathcal
S^\theta_0$, which in turn is realized by a self-adjoint symmetry
of the Hilbert space $\Hcal^\theta$. Hence, formula
\eqref{eq:spflreldim} can be used to compute the spectral flow of
$\mathcal S^\theta$, and we get:
\begin{equation}\label{eq:spflreldimbis}
\iMorse(\theta)=\Dim\big(\Ker(\mathcal
S^\theta_1)\big)-\Dim_{V^{\!\scriptscriptstyle-}(\mathcal
S^\theta_1)}\big(V^-(\mathcal S^\theta_0)\big)=
\Dim_{V^{\!\scriptscriptstyle+}(\mathcal
S^\theta_0)}\big(V^+(\mathcal S^\theta_1)\big).
\end{equation}

\subsection{The spectral index}
\label{sub:spectral}
Important classes of examples where one can apply the Maslov index theory
for real-analytic curves arises naturally when one studies certain
eigenvalue problems for ODE's, whose solutions depend analytically on
the eigenvalue by standard regularity results.

We will  consider in what follows the case of Morse--Sturm--Liouville
equation in $\R^n$, whose
spectral index is given as the Maslov index of a certain curve parameterized
by the spectral parameter $\lambda$. In this case, each negative (real) eigenvalue
of the equation gives a contribution to the Maslov index, and it is possible to compute
explicitly the bilinear forms $B_k$ at each eigenvalue.

Let us consider a nondegenerate symmetric bilinear form $g:\R^n\times\R^n\to\R$
and a continuous map $[0,1]\ni t\mapsto R(t)$ of $g$-symmetric
endomorphisms of $\R^n$. The {\em Morse--Sturm--Liouville equation\/} with data
($g$ and) $R$ is given by:
\begin{equation}\label{eq:MSL}
-v''+(R-\lambda)v=0,\end{equation}
where $v:[0,1]\to\R^n$ and $\lambda\in\R$. The corresponding
differential operator, denoted by $\mathfrak J_\lambda$:
\[\mathfrak J_\lambda=-\frac{\mathrm d^2}{\mathrm dt^2}+(R-\lambda),\]
defined on the domain $\mathcal D=H^1_0\big([0,1],\R^n\big)\cap H^2\big([0,1],\R^n\big)$,
is an unbounded and, unless $g$ is positive or negative definite, non normal linear operator.
We will consider the following symmetric bounded and nondegenerate bilinear form
$\widehat g$ on $L^2\big([0,1],\R^n\big)$:
\[\widehat g(v,w)=\int_0^1g\big(v(t),w(t)\big)\,\mathrm dt.\]
It is easy to see that $\mathfrak J_\lambda$ is $\widehat g$-symmetric:
\[I^\theta_\lambda(v,w):=\widehat g\big(\mathfrak J_\lambda
v,w\big)\!=\!\int_0^1\Big[g\big(v'(t),w'(t)\big)+
g\big((R(t)-\lambda)v(t),w(t)\big)\Big]
\mathrm dt=\widehat g\big(v,\mathfrak J_\lambda
w\big),\] for all $v,w\in\mathcal D$; in particular,
$I^\theta_0=\widehat g(\mathfrak J_0\cdot,\cdot)$ coincides with the index form $S^\theta_1$ as defined in \ref{eq:defIt}.
Moreover, it is easy to see that $I^\theta_\lambda$ is a Fredholm bilinear form
for all $\lambda$ (it is a compact perturbation of $I^\theta_0$).

The spectral properties of $\mathfrak J_0$ have been studied in \cite{CPPT},
we will recall here some facts: $\mathfrak J_0$ is discrete (i.e., it has compact resolvent),
its spectrum $\mathfrak s(\mathfrak J_0)$ is a discrete subset of   the strip:
\[\big\{z\in\C:\Re(z)\ge -\Vert R\Vert_\infty,\ \vert \Im(z)\vert\le\Vert R\Vert_\infty\big\},\]
where $\Vert R\Vert_\infty$ denotes the supremum norm of $R$.

\begin{defin}\label{thm:defspectralindex}
The {\em spectral index\/} $\ispec(\theta)$ of the geodesic $\theta$ is
defined to be integer:
 \begin{equation}\label{eq:defspectralindex}
\ispec(\theta)=\Dim\big(\Ker(I^\theta_0)\big)-
\spfl\big(I^\theta_\lambda,[-M_0,0]\big),\end{equation}
where $M_0>\Vert R\Vert_\infty$.
\end{defin}
Both the generalized Morse index $\iMorse(\theta)$ and the
spectral index $\ispec(\theta)$ do not depend on the
choice of the parallel trivialization of $TM$ along $\theta$;
this fact can be proven directly, or obtained as a consequence
of Theorem~\ref{thm:tuttouguale}.

As an easy application of Proposition~\ref{thm:surpresa2}, we obtain the following:
\begin{prop}\label{thm:formulaspindex}
 The {\em generalized eigenspace\/}
$\Hcal_\lambda=\bigcup_{n\ge1}\Ker\big((\lambda-\mathfrak J_0)^n\big)$ is finite-dimensional
for all $\lambda\in\mathfrak s(\mathfrak J_0)$, and the restriction of $\widehat g$ to $\Hcal_\lambda$ is nondegenerate.
The spectral index of $\theta$ is given by:
\begin{equation}\label{eq:defspectraindex}
\ispec(\theta)=n^+\big(\widehat
g\vert_{\Hcal_0}\big) -n^+\big(I_0^\theta\vert_{\Hcal_0}\big) +\sum_{\substack{\lambda\in\mathfrak s(\mathfrak J_0)\\
\lambda\in\left]-\infty,0\right[}}\sgn\big(\widehat g\vert_{\Hcal_\lambda}\big).
\end{equation}
\end{prop}
\begin{proof}
Apply  Proposition~\ref{thm:surpresa2} to the Hilbert spaces \[\Hcal=L^2\big([0,1],\R^n\big),\quad
\Dcal_2=H^2\big([0,1],\R^n\big)\cap H^1_0\big([0,1],\R^n\big),\quad \Dcal_1=H^1_0\big([0,1],\R^n\big),\]
the operator $T=\mathfrak J_0$, and the self-adjoint isomorphism $G:\Hcal\to\Hcal$ given
by pointwise composition with the (constant) symmetric endomorphism $g$ of $\R^n$.

Observe that $-\frac{\mathrm d^2}{\mathrm dt^2}:\Dcal_2\to\Hcal$ is
and isomorphism, and $T$ is a compact perturbation of such isomorphism, hence
a Fredholm operator of index $0$.
As to the ``regularity'' condition assumed in the hypotheses of
Proposition~\ref{thm:surpresa2}, in our case it follows easily from
the fact that, using standard bootstrap arguments, if $v\in H^1_0\big([0,1],\R^n\big)$
is such that there exists $z\in H^1_0\big([0,1],\R^n\big)$ with
\[\int_0^1g(v',w')\,\mathrm dt=\int_0^1g(z,w)\,\mathrm dt\]
for all $w\in H^1_0\big([0,1],\R^n\big)$, then $v\in H^2\big([0,1],\R^n\big)$.

In order to obtain \eqref{eq:defspectraindex} from \eqref{eq:pp4} keep in mind that
the boundary term corresponding to $\lambda=-M_0$ is null, because
$-M_0\not\in\mathfrak s(\mathfrak J_0)$.
\end{proof}
Observe that when $g$ is Riemannian, then
$\mathfrak J_0$ is indeed self-adjoint,  $\Hcal_\lambda= \Ker(\lambda-\mathfrak J_0)$,
in particular $S^\theta_1\vert_{\Hcal_0}=0$ and thus
$n^+\big(\mathcal S_1^\theta\big\vert_{\Hcal_0}\big)=0$.
Moreover, since $\widehat g$ is positive definite, $n^+\big(\widehat g\vert_{\Hcal_\lambda}\big)=\sgn\big(\widehat
g\vert_{\Hcal_\lambda}\big)=\Dim\big(\Ker(\lambda-\mathfrak J_0)\big)$, and
$n^-(\widehat g)=0$; this shows that
in the Riemannian case the spectral index of $\theta$ coincides with
the  {\em extended Morse index\/} of the index form
$\mathcal S^\theta_1=I^\theta_0$.

\subsection{The index theorem}
\label{sub:indexthm}
Before we get into the aimed index theorem, we will need
to introduce the following notation. For all $\lambda\in\R$,
let $\Phi_\lambda:[0,1]\to\Spl(\R^{2n},\omega_g)$ denote
the flow of the Morse--Sturm--Liouville equation \eqref{eq:MSL},
i.e., $\Phi_\lambda(t):\R^{2n}\to\R^{2n}$ is the linear
isomorphism defined by:
\[\Phi_\lambda(t)\big(v(0),v'(0)\big)=\big(v(t),v'(t)\big),\]
for all solutions $v$ of \eqref{eq:MSL}.

We choose $M_0>\Vert R\Vert_\infty$
and we consider the curve $\ell:[-M_0,0]\to\Lambda$:
\begin{equation}\label{eq:elllambda}
\ell(\lambda)=\Phi_\lambda(1)\big(L_0\big).
\end{equation}
Observe that $\lambda$ is real-analytic; moreover, the
intersections of $\ell$ with the Maslov cycle occur
precisely at each real nonpositive eigenvalue of the Jacobi differential operators
$\mathfrak J_0$.

Everything is now ready to state and prove the following:

\begin{teo}[Index Theorem in the Degenerate Case]
\label{thm:tuttouguale}
Let $(M,g)$ be a semi-Riemannian manifold and $\theta:[0,1]\to M$ a
geodesic. Then:
\[\iMaslov(\theta)=\ispec(\theta)=\iMorse(\theta).\]
\end{teo}
\begin{proof}
The equality $\ispec(\theta)=\iMorse(\theta)$ is proven by an infinite dimensional
homotopy argument. Namely, consider the
two-parameter smooth map $C$ of symmetric bilinear forms  in $H^1_0\big([0,1],\R^n\big)$:
\begin{equation}\label{thm:Ctlambda}
C(t,\lambda)(V,W)=\int_0^1\Big[g\big(V'(s),W'(s)\big)+t^2g\big((R(ts)-\lambda)V(s),W(s)\big)\Big]\,\mathrm ds,
\end{equation}
$(t,\lambda)\in[0,1]\times[-M_0,0]$.

Observe that, by definition,
$-\ispec(\theta)+\Dim\big(\Ker(C(1,0))\big)$ equals the spectral flow of the curve $[-M_0,0]\ni\lambda\mapsto C(1,\lambda)$,
while $-\iMorse(\theta)+\Dim\big(\Ker(C(1,0))\big)$ equals the spectral flow of $[0,1]\ni t\mapsto C(t,0)$.
The equality $\ispec(\theta)=\iMorse(\theta)$ follows from
the fixed-endpoints homotopy invariance and the additivity by concatenation
of the spectral flow, observing
that the maps $[-M_0,0]\ni\lambda\mapsto C(0,\lambda)$ and
$[0,1]\ni t\mapsto C(t,-M_0)$ have null spectral flow, due to the fact that
$C(0,\lambda)$ and $C(t,-M_0)$ are always nondegenerate.

By a similar homotopy argument in $\Lambda$, one proves that $\iMaslov(\theta)$
equals the Maslov index of the curve $\ell$ defined in \eqref{eq:elllambda}. Namely,
$\iMaslov(\theta)$ is by definition the Maslov index of the
curve $\gamma$ given by $[\varepsilon,1]\ni t\mapsto\Phi_0(t)(L_0)$, and the two-parameter
map $[\varepsilon,1]\times[-M_0,0]\ni(t,\lambda)\mapsto\Phi_\lambda(t)(L_0)\in\Lambda$
gives a continuous homotopy between $\gamma$ and $\ell$.
In this case, observe that the curve $[\varepsilon,1]\ni t\mapsto\Phi_{-M_0}(t)(L_0)$
does not intersect the Maslov cycle for all $\varepsilon\ge0$, while the
curve $[-M_0,0]\ni\lambda\mapsto\Phi_\lambda(\varepsilon)(L_0)$ does not intersect
the Maslov cycle provided that $\varepsilon>0$ is chosen sufficiently small.

Finally, the crucial part of the proof consists in showing that
$\ispec(\theta)$ equals $\iMaslov(\ell)$; in this case
a direct homotopy argument cannot be used, because
$\ispec(\theta)$ is the negative spectral flow of the
path of Fredholm bilinear forms $\lambda\mapsto C(1,\lambda)$
on $H^1_0([0,1],\R^n)$, while $\iMaslov(\ell)$ is the Maslov
index of a curve in $\Lambda$. However, as we have observed
both curves are real-analytic, and they have precisely the same
degeneracy instants.

Observe that for this equality one can use the partial signatures
theory, since both integers are Maslov indexes of real-analytic paths:
for each eigenvalue $\lambda_0\in[-M_0,0]$, we prove that the spaces
$W_{k+1}$ and the bilinear forms $B_{k+1}$ obtained from the
two constructions coincide, up to the sign.

Let $\lambda_0\in\left[-M_0,0\right]$ be a degenerate value for
$\lambda\mapsto C(1,\lambda)$ and let
$T:\R^n\to\R^n$ be a $g$-symmetric linear endomorphism
such that the Lagrangian $L_1'=\Gr(T^{-1})$ is transversal
to $\Phi(1,\lambda_0)(L_0)$.
Let us consider the following real analytic path $\lambda\mapsto \mathfrak B_\lambda$
of Fredholm symmetric bilinear forms:
\[\mathfrak B_\lambda(V,W)=\int_0^1\big[g\big(V',W')+g\big((R-\lambda)V,W\big)\big]\mathrm ds-
g\big(V(1),T^{-1}W(1)\big)\]
defined on the Hilbert space:
\[\Hcal=\big\{V\in H^1([0,1],\R^n):V(0)=0\big\}.\]
Observe that the restriction of $\mathfrak B_\lambda$ to $H^1_0\big([0,1],\R^n\big)$
coincides with $C(1,\lambda)$, while the restriction
of $\mathfrak B_\lambda$ to the finite dimensional space:
\[\mathbb J_\lambda:=\big\{J\in C^2\big([0,1],\R^n\big):
J''=(R-\lambda)J,\ J(0)=0\big\}\subset\Hcal\]
gives the symmetric bilinear form $D_\lambda:\mathbb J_\lambda\times\mathbb J_\lambda\to\R$:
\[D_\lambda(J_1,J_2)=g\big(J_1(1),J_2'(1)-T^{-1}J_2(1)\big).\]
Thus, recalling Proposition~\ref{thm:operational}, if  $\kappa_\lambda:L_0\to\mathbb J_\lambda$ is the isomorphism
$v\mapsto J_v$, where $J_v$ is the unique element of $ \mathbb J_\lambda$
satisfying $J_v'(0)=v$, we consider the pull-back $(\kappa_\lambda)^*(D_\lambda)$,
which is the symmetric bilinear form
$\widetilde D_\lambda:L_0\times L_0\to\R$ given by:
\[\widetilde D_\lambda(v,w):=B_\lambda\big(\kappa_\lambda(v),
\kappa_\lambda(w)\big)=
\varphi_{L_0,L_1'}\big(\ell(\lambda)\big)(v,w),\] for all $v,w\in
L_0$. We observe that:
\[\Ker\big(C(1,\lambda_0)\big)=\mathbb J_{\lambda_0}\cap
H^1_0\big([0,1],\R^n\big)=\Ker(D_{\lambda_0})=\kappa_{\lambda_0}\big(\Ker(\widetilde
D_{\lambda_0})\big).\]

Our aim is to show that, for all $k\ge0$,
\begin{equation}\label{eq:Wkuguali}
W_{k}\big(C(1,\lambda);\lambda_0\big)=\kappa_{\lambda_0}\big(W_{k}
\big(\widetilde D_\lambda;\lambda_0\big)\big),\end{equation}
and that:
\begin{equation}\label{eq:Bkuguali}B_{k}\big(C(1,\lambda);\lambda_0\big)=-
(\kappa_{\lambda_0})^*\big(B_{k}
\big(\widetilde D_\lambda;\lambda_0\big)\big),\end{equation}
for all eigenvalue $\lambda_0\in[-M_0,0]$.
Recalling Proposition~\ref{thm:cambiasegno}, the equality
$\ispec(\theta)=\iMaslov(\ell)$ will follow at once from
\eqref{eq:Wkuguali} and \eqref{eq:Bkuguali}.

Let $k\ge1$ be fixed, and choose $v_0\in L_0$; if $\lambda\mapsto v_\lambda$
is a root function of order greater than or equal to $k$ for
$\widetilde D_\lambda$ at $\lambda=\lambda_0$, and such that $v_{\lambda_0}=v_0$,
then by Corollary~\ref{thm:Wkconjinst} the map
$\lambda\mapsto \kappa_\lambda(v_\lambda)(1)$ has a zero of order greater
than or equal to $k$ at $\lambda=\lambda_0$. Thus, the
map $\lambda\mapsto\widetilde J_\lambda\in H^1_0\big([0,1],\R^n\big)$
defined by:
\[\widetilde J_\lambda(t)=J_\lambda(t)-tJ_\lambda(1),\]
where
\[J_\lambda=\kappa_\lambda(v_\lambda),\] is a root function of order
greater than or equal to $k$ for $C(1,\lambda)$ at
$\lambda=\lambda_0$, with $\widetilde
J_{\lambda_0}=\kappa_{\lambda_0}(v_0)$. For, an easy computation
shows that for all $W\in H^1_0\big([0,1],\R^n\big)$:
\[C(1,\lambda)\big(\widetilde J_\lambda,W\big)=-\int_0^1tg\big((R(t)-\lambda)
J_\lambda(1),W(t)\big)\,\mathrm dt.\] This shows that we have an
inclusion that $W_{k}\big(C(1,\lambda);
\lambda_0\big)\supset\kappa_{\lambda_0}\big(W_{k} \big(\widetilde
D_\lambda;\lambda_0\big)\big)$.

To prove the opposite inclusion, we use the result of
part~\eqref{itm:span} in Proposition~\ref{thm:central}
arguing as follows. Choose a positive definite inner product
$\gr$ in $\R^n$ and let $A:\R^n\to\R^n$ the $g$-symmetric automorphism
such that $\gr=g(A\cdot,\cdot)$; define the following  Hilbert space
inner product in $H_0^1\big([0,1],\R^n\big)$:
\[\llangle V,W\rrangle:=\int_0^1\gr\big(V'(t),W'(t)\big)\,\mathrm dt=
\int_0^1g\big(AV'(t),W'(t)\big)\,\mathrm dt.\]
Let $\lambda\mapsto{\LL}_\lambda\in\Fcalsa\big(H^1_0\big([0,1],\R^n\big)\big)$
be the real-analytic path of self-adjoint operators that realize
$C(1,\lambda)$ with respect to the above inner product, and let
$\lambda\mapsto\sigma(\lambda)\in\R$,
$\lambda\mapsto V_\lambda\in H^1_0\big([0,1],\R^n\big)$ be  real-analytic
maps such that:
\begin{itemize}
\item  ${\LL}_\lambda(V_\lambda)=\sigma(\lambda)V_\lambda$ for all
$\lambda$;
\item $\sigma$ has a zero of order greater than or equal to
$k$ at $\lambda=\lambda_0$.
\end{itemize}
Then, for all $W\in H^1_0\big([0,1],\R^n\big)$ the following equality
holds:
\[\int_0^1\Big[g\big((\mathrm{Id}-\sigma(\lambda)A)V_\lambda'(t),W'(t)\big)+
g\big(R(t)-\lambda)V_\lambda(t),W(t)\big)\Big]\,\mathrm dt=0,\]
from which it follows that $V_\lambda$ is a map of class
$C^2$ that satisfies the ``perturbed'' Jacobi equation:
\[V_\lambda''=\big(\mathrm{Id}-\sigma(\lambda)A\big)^{-1}(R-\lambda)V_\lambda.\]
Observe that for $\lambda$ near $\lambda_0$, the operator $\mathrm{Id}-\sigma(\lambda)A$
is invertible.
Observe that also $V_{\lambda_0}$ is in $\Ker\big(C(1,\lambda_0)\big)$, and
by Proposition~\ref{thm:central}, the space $W_{k+1}\big(C(1,\lambda);\lambda_0\big)$
is generated by such functions   $V_{\lambda_0}$.
The map $\lambda\mapsto v_\lambda=V_\lambda'(0)\in L_0$ is real-analytic,
and we claim that it is a root function for $\widetilde D_\lambda$ at
$\lambda=\lambda_0$ of order greater than or equal to $k$, with
$\kappa_{\lambda_0}(v_{\lambda_0})=V_{\lambda_0}$.

\noindent\  By Corollary~\ref{thm:Wkconjinst}, to prove our
assertion it suffices to show that the map
$\lambda\mapsto\kappa_{\lambda}(v_\lambda)(1)\in\R^n$ has a zero
of order greater than or equal to $k$ at $\lambda=\lambda_0$. To
prove this fact, observe   that
$J_\lambda=\kappa_\lambda(v_\lambda)$ is the solution of the
equation $J_\lambda''=(R-\lambda)J_\lambda$ satisfying
$J_\lambda(0)=0$ and $J_\lambda'(0)=v_\lambda$. For $\lambda$ near
$\lambda_0$, say $\lambda\in[\lambda_0-\delta,\lambda_0+\delta]$,
$\delta>0$,  we set:
\[B(\lambda)=\sum_{n=1}^\infty\sigma(\lambda)^{n-1}A^n,\]
so that
$(\mathrm{Id}-\sigma(\lambda)A)^{-1}=\mathrm{Id}+\sigma(\lambda)B(\lambda)$;
then, for $s\in[0,1]$:
\[\begin{split}
\big\vert V_\lambda(s)-&J_\lambda(s)\big\vert\le\int_0^s\big\vert
V_\lambda'(\tau)-J_\lambda'(\tau)\big\vert\,\mathrm
d\tau\le\int_0^s\mathrm d\tau\int_0^\tau\big\vert
V_\lambda''(r)-J_\lambda''(r)\big\vert\,\mathrm dr\\&
=\int_0^s\mathrm
d\tau\int_0^\tau\big\vert(R(r)-\lambda)(V_\lambda(r)-J_\lambda(r))+
\sigma(\lambda)B(\lambda)(R(r)-\lambda)V_\lambda\big\vert\,\mathrm
dr\\
&=\int_0^s\mathrm
dr\int_r^s\big\vert(R(r)-\lambda)(V_\lambda(r)-J_\lambda(r))+
\sigma(\lambda)B(\lambda)(R(r)-\lambda)V_\lambda\big\vert\,\mathrm
d\tau\\
&=\int_0^s\
(s-r)\big\vert(R(r)-\lambda)(V_\lambda(r)-J_\lambda(r))+
\sigma(\lambda)B(\lambda)(R(r)-\lambda)V_\lambda\big\vert\,\mathrm
dr\\&\le\int_0^s\big\vert
(R(r)-\lambda)(V_\lambda(r)-J_\lambda(r))\big\vert+\vert\sigma(\lambda)\vert\big\vert
B(\lambda)(R(r)-\lambda))V_\lambda(r)\big\vert\,\mathrm dr\\& \le
d_0\int_0^s\big\vert (V_\lambda(r)-J_\lambda(r))\big\vert\;\mathrm
dr+d_1\vert\sigma(\lambda)\vert,
\end{split}\]
where: \[d_0=\max_{\substack{r\in[0,1]\\
\lambda\in[\lambda_0-\delta,\lambda_0+\delta]}}\Vert
R(r)-\lambda\,\mathrm{Id}\Vert,\quad d_1=\max_{\substack{r\in[0,1]\\
\lambda\in[\lambda_0-\delta,\lambda_0+\delta]}}\Vert
B(\lambda)(R(r)-\lambda\,\mathrm{Id})\Vert\cdot \vert
V_\lambda(r)\vert.\] From Gronwall's Lemma, we then obtain:
\[\Vert V_\lambda-J_\lambda\Vert_\infty\le d_1e^{d_0}\vert\sigma(\lambda)\vert,\]
and since $V_\lambda(1)=0$, it follows that $J_\lambda(1)$ has a
zero of order greater than or equal to $k$.
This argument shows that we have an inclusion $W_{k}\big(C(1,\lambda);
\lambda_0\big)\subset\kappa_{\lambda_0}\big(W_{k}
\big(\widetilde D_\lambda;\lambda_0\big)\big)$ and \eqref{eq:Wkuguali}
is proven.

Finally, we will now prove equality \eqref{eq:Bkuguali};
it will suffice to show that, given root functions
$\lambda\mapsto u_\lambda,v_\lambda\in L_0$ of order greater than or equal to
$k$ for $\widetilde D_\lambda$
at $\lambda=\lambda_0$, then:
\begin{equation}\label{eq:derivataCD}\frac{\mathrm d^k}{\mathrm d\lambda^k}\Big\vert_{\lambda=\lambda_0}
\widetilde D_\lambda\big(u_\lambda,v_\lambda\big)=-
\frac{\mathrm d^k}{\mathrm d\lambda^k}\Big\vert_{\lambda=\lambda_0}
C(1,\lambda)\big(\widetilde J_\lambda^{(1)},\widetilde J_\lambda^{(2)}\big),
\end{equation}
where:
\[\widetilde
J_\lambda^{(1)}(t)=J_\lambda^{(1)}(t)-tJ_\lambda^{(1)}(1), \quad
\widetilde
J_\lambda^{(2)}(t)=J_\lambda^{(2)}(t)-tJ_\lambda^{(2)}(1),\] and
\[J_\lambda^{(1)}=\kappa_\lambda(u_\lambda),\quad
J_\lambda^{(2)}=\kappa_\lambda(v_\lambda).\] A direct computation
gives:
\[\widetilde D_\lambda\big(u_\lambda,v_\lambda\big)=g\big(
J_\lambda^{(1)}(1),(J_\lambda^{(2)})'(1)-T^{-1}J_\lambda^{(2)}(1)\big),\]
and since $g\big(J_\lambda^{(1)}(1),
T^{-1}J_\lambda^{(2)}(1)\big)$ has a zero of order greater than or
equal to $2k$ at $\lambda=\lambda_0$, it follows:
\begin{equation}\label{eq:derivataD}
\frac{\mathrm d^k}{\mathrm
d\lambda^k}\Big\vert_{\lambda=\lambda_0} \widetilde
D_\lambda\big(u_\lambda,v_\lambda\big)= \frac{\mathrm d^k}{\mathrm
d\lambda^k}\Big\vert_{\lambda=\lambda_0}
g\big(J_\lambda^{(1)}(1),(J_\lambda^{(2)})'(1)\big).\end{equation}
On the other hand, integration by parts yields the following:
\[
C(1,\lambda)\big(\widetilde J_\lambda^{(1)},\widetilde
J_\lambda^{(2)}\big)=
-\int_0^1t\,g\big((R(t)-\lambda)J_\lambda^{(1)}(1),\widetilde
J_\lambda^{(2)}(t) \big)\,\mathrm dt;\] and, again because of the
fact that $g\big(J_\lambda^{(1)}(1), J_\lambda^{(2)}(1)\big)$ has
a zero of order greater than or equal to $2k$ at
$\lambda=\lambda_0$, we get:
\begin{equation}\label{eq:derivataC}
\begin{split}
\frac{\mathrm d^k}{\mathrm
d\lambda^k}&\Big\vert_{\lambda=\lambda_0}
C(1,\lambda)\big(\widetilde J_\lambda^{(1)},\widetilde
J_\lambda^{(2)}\big)\\&= -\frac{\mathrm d^k}{\mathrm
d\lambda^k}\Big\vert_{\lambda=\lambda_0}
\int_0^1tg\big((R(t)-\lambda)J_\lambda^{(1)}(1),
J_\lambda^{(2)}(t)\big)\,\mathrm dt\\&= -\frac{\mathrm
d^k}{\mathrm d\lambda^k}\Big\vert_{\lambda=\lambda_0}
\int_0^1tg\big(J_\lambda^{(1)}(1),(R(t)-\lambda)J_\lambda^{(2)}(t)\big)\,\mathrm
dt
\\&=-\frac{\mathrm d^k}{\mathrm d\lambda^k}\Big\vert_{\lambda=\lambda_0}
\int_0^1tg\big(J_\lambda^{(1)}(1),(J_\lambda^{(2)})''(t)\big)\,\mathrm
dt
\\&=
\frac{\mathrm d^k}{\mathrm
d\lambda^k}\Big\vert_{\lambda=\lambda_0}
\Big[\int_0^1g\big(J_\lambda^{(1)}(1),
(J_\lambda^{(2)})'(t)\big)\,\mathrm dt-g\big(J_\lambda^{(1)}(1),
(J_\lambda^{(2)})'(1)\big)\Big]
\\&= \frac{\mathrm d^k}{\mathrm d\lambda^k}\Big\vert_{\lambda=\lambda_0}
\Big[g\big(J_\lambda^{(1)}(1),
J_\lambda^{(2)}(1)\big)-g\big(J_\lambda^{(1)}(1),
(J_\lambda^{(2)})'(1)\big)\Big]\\&= -\frac{\mathrm d^k}{\mathrm
d\lambda^k}\Big\vert_{\lambda=\lambda_0} g\big(J_\lambda^{(1)}(1),
(J_\lambda^{(2)})'(1)\big).
\end{split}\end{equation}
Comparison of \eqref{eq:derivataD} and \eqref{eq:derivataC}
give \eqref{eq:derivataCD}, and the proof is concluded.
\end{proof}

\subsection{On a counterexample for the equality of the conjugate and the Maslov
index}\label{sub:counterexample}
The sum of the signatures of
the conjugate points along a  geodesic in a real-analytic
semi-Riemannian manifold  is called the {\em conjugate index\/} in
reference \cite{Hel1} (or {\em focal index}
in reference \cite{pacific}). It was erroneously stated in
\cite{Hel1} that the conjugate index is equal to the
Maslov index of a semi-Riemannian geodesic, and in \cite[Subsection~5.4]{pacific}
the authors have given a counterexample to such equality,
occurring in the case of a degenerate conjugate point along a Lorentzian
spacelike geodesic.
Using the results of the present paper we are now able
to have a better view of the phenomenon.

\noindent
Recall that the counterexample mentioned consists in a spacelike geodesic
$\theta\!:\![-\varepsilon,\varepsilon]\!\to~M$, where $M$ is a real-analytic three-dimensional Lorentzian
manifold,  $\varepsilon>0$, having a unique
conjugate point at $t=0$. By a parallel transport of the normal
bundle $\dot\theta^\perp$ along $\theta$, and using suitable
coordinate systems in the Lagrangian Grassmannian as
explained above, the Maslov index of $\theta$ is
computed as the spectral flow through $t=0$ of the curve of symmetric
bilinear forms on $\R^2$ given by:
\[[-\varepsilon,\varepsilon]\ni t\longmapsto{\LL}(t)=\begin{pmatrix}x(t)&z(t)\\ z(t)&y(t)\end{pmatrix},\]
where:
\[x(t)=-2t^3-\frac{54}5t^5,\quad y(t)=-1-6t+18t^2-54t^3,\quad z(t)=-3t^2.\]
One computes easily:
\[{\LL}_0=\begin{pmatrix}0&0\\ 0&-1\end{pmatrix},\quad {\LL}_1=\begin{pmatrix}0&0\\ 0&-6\end{pmatrix},
\quad {\LL}_2=\begin{pmatrix}0&3\\ 3&18\end{pmatrix},\quad {\LL}_0=\begin{pmatrix}-2&0\\ 0&-54\end{pmatrix},\]
hence $W_1=\Ker({\LL}_0)=\R\oplus\{0\}$ and $B_1=\langle L_1\cdot,\cdot\rangle\vert_{W_1}=0$,
and $t=0$ is a degenerate conjugate instant having signature equal to $0$.
Nevertheless, the Maslov index of $\theta$, which is equal to the spectral flow of ${\LL}$,   is
easily computed as:
\[n^+\big({\LL}(\varepsilon)\big)-n^+\big({\LL}(-\varepsilon)\big)=-1,\]
providing a counterexample for the equality between the
Maslov and the conjugate index in the case of degenerate
conjugate points. An elementary direct computation gives:
\[W_2=\R\oplus\{0\}, \ B_2=0,\quad W_3=\R\oplus\{0\},\ B_3\big((\alpha,0),(\beta,0)\big)=-2\alpha\beta,\]
hence $\sigma_1=\sigma_2=0$, $\sigma_3=-1$ and $\sigma_k=0$ for all $k>3$.
Obviously,   $\iMaslov(\theta)$ equals the sum
$\sum_{k\ge1}\sigma_k(\theta,0)$, as it must be in compliance
with Proposition~\ref{thm:maslovgeo}.
\subsection{A geometrical version of the semi-Riemannian index theorem}
\label{sub:geomversion}
Using an abstract result on the computation of the relative index of Fredholm
bilinear forms, we will now give a geometrical
version of the index theorem, in the spirit of the semi-Riemannian index
theorem in \cite{topology}. The Maslov index of a semi-Riemannian
geodesic can be computed as the difference between the index and
the coindex of suitable restrictions of the index form $S^\theta_1$  \eqref{eq:defIt}; this version
of the index theorem gives a link with variational problems
and aims at developments of Morse homology in the infinite dimensional
Hilbert manifold of all paths in $M$ joining two fixed points, as in \cite{asian, calcvar}.

Let $(M,g)$ be a semi-Riemannian manifold, set $k=n^-(g)$.
Let $\theta:[0,1]\to M$ be a geodesic; a {\em maximal negative distribution\/}
along $\theta$ is a smooth family $\mathcal D_t\subset T_{\theta(t)}M$
of $k$-dimensional subspaces, $t\in[a,b]$, such that $g\vert_{\mathcal D_t}$ is
negative definite for all $t$. By ``smooth'', we mean that $\mathcal D_t$
is the span of $Y_1(t),\ldots,Y_k(t)$ for all $t\in[0,1]$, where $Y_1,\ldots,Y_k$ is a family
of smooth vector fields along $\theta$; such a family $Y_1,\ldots,Y_k$ will be called
a {\em frame\/}  for $\mathcal D$.
Associated to each choice of a maximal
negative distribution $\mathcal D$ along $\theta$ one can define
two closed spaces of variational vector fields along $\theta$:
the space of vector fields along $\theta$ taking values in $\mathcal D$, denoted
by $\mathcal Q$, and the space of vector fields  along $\theta$ that
are ``Jacobi in the directions of $\mathcal D$'', denoted by $\mathcal K$.
More precisely, denote by $\Hcal^\theta$ the space of all vector
fields of Sobolev class $H^1$ along $\theta$ vanishing at the endpoints;
fix a frame $Y_1,\ldots,Y_k$ for $\mathcal D$ and define:
\[\begin{split}&\mathcal Q=\big\{v\in \Hcal^\theta:v(t)\in\mathcal D_t,\ \text{for all}\  t\in[0,1]\big\},\\
&\mathcal  K=\big\{v\in\Hcal^\theta:g(v',Y_i)\ \text{is of class $H^1$,}\  g(v',Y_i)'=g(v',Y_i')+g\big(R
(\dot\theta,v)\,\dot\theta,Y_i\big)\ \forall\,i\big\}.\end{split}\]
Observe that a vector field $v$ of class $C^2$ along $\theta$ belongs to
$\mathcal K$ if and only if $v''-R(\dot\theta,v)\dot\theta$ is pointwise orthogonal
to $\mathcal D$, i.e., fields in $\mathcal K$ are interpreted as ``Jacobi fields
in the directions of $\mathcal D$''; geometrical and analytical descriptions of the
spaces $\mathcal Q$ and $\mathcal K$ can be found in \cite{topology}.
\begin{prop}\label{thm:geomversion}
The restriction $S^\theta_1$ to $\mathcal K$ has finite index, and
the restriction of $S^\theta_1$ to $\mathcal Q$ has finite
coindex; moreover, the following equality holds
\begin{equation}\label{eq:p34-b}
\iMorse(\theta)=n^-\big(S^\theta_1\vert_{\mathcal
K}\big)-n^+\big(S^\theta_1\vert_{\mathcal Q}\big)+
\Dim\big(\Ker(S^\theta_1)\big).\end{equation}
\end{prop}
\begin{proof}
The restriction of $S^\theta_1$ to $\mathcal Q$ is realized by a
compact perturbation of a negative isomorphism of $\mathcal Q$
(\cite[Corollary~5.25]{topology}), hence it has finite coindex.
The restriction of $S^\theta_1$ to $\mathcal K$ is realized by a
compact perturbation of a positive isomorphism of $\mathcal K$
(\cite[Lemma~2.6.6]{london}), hence it has finite coindex.
Moreover, an immediate calculation shows that $\mathcal K=\mathcal
Q^{\perp_{S^\theta_1}}$.

By Lemma~\ref{thm:RCP-comm}, $\mathcal Q$ is commensurable with
$V^{\scriptscriptstyle-}
(\mathcal S_1^\theta)$, and using the abstract result of Proposition~\ref{thm:indrelindcoind}, the relative
dimension $\Dim_{\mathcal Q}\big(V^{\scriptscriptstyle-}(\mathcal S^\theta_1)\big)$
can be computed as:
\begin{equation}\label{eq:eeqq1}
\Dim_{\mathcal Q}\big(V^{\scriptscriptstyle-}(\mathcal
S^\theta_1)\big)= n^-\big(S^\theta_1\vert_{\mathcal
K}\big)-n^+\big(S^\theta_1\vert_{\mathcal Q}\big),
\end{equation}
where $V^{\scriptscriptstyle-}(\mathcal S^\theta_1)$ is the
negative eigenspace of the realization of $\mathcal S^\theta_1$
relatively to any Hilbert structure on $\Hcal^\theta$. In order to
compute the left hand side in equality \eqref{eq:eeqq1}, we will
first show that its value does not depend on the choice of a
maximal negative distribution along $\theta$.

The idea to prove the independence of the relative index
from the choice of a maximal negative distribution consists
in showing that any two maximal negative distributions
can be joined by a ``continuous'' selection of maximal negative
distributions, and that the relative index depends ``continuously''
on such selection. Let us make the argument formal, as follows.
In first place, using a parallel trivialization of the tangent bundle
$TM$ along $\theta$, the problem is reduced to studying the
equality of the relative dimensions:
\[\Dim_{\mathcal Q^0}\big(\Hcal_0)\quad\text{and}\quad
\Dim_{\mathcal Q^1}\big(\Hcal_0),\]
where
\begin{itemize}
\item $\Hcal_0$ is a  closed subspace ($\cong V^{\scriptscriptstyle-}
(\mathcal S^\theta_1)$) of
the Sobolev space $\Hcal=H^1_0\big([a,b],\R^n\big)$ ($\cong \Hcal^\theta$)
\item $g$ is a fixed nondegenerate symmetric bilinear form
on $\R^n$ of index $k$;
\item denoting by $G_k^-(n)$ the Grassmannian of all $k$-dimensional
subspaces of $\R^n$ on which $g$ is negative definite,
$\Dcal^0,\Dcal^1:[a,b]\to G_k^-(n)$ are continuous curves;
\item $\mathcal Q^i=\big\{v\in\Hcal:v(t)\in\Dcal^i(t)\ \text{for all}\
t\in[a,b]\big\}$,
$i=0,1$, are closed subspaces of $\Hcal$
that are
commensurable to $\mathcal Q_0$ and $\mathcal Q_1$.
\end{itemize}
The set $G_k^-(n)$ is open in the Grassmannian $G_k(n)$, and
it is an arc-connected set (see Appendix~\ref{sec:grasconn}); hence,
the curves $\Dcal^0$ and $\Dcal^1$ are homotopic. Choose
a continuous map $[a,b]\times[0,1]\ni(t,s)\mapsto\Dcal_{t,s}\in G_k^-(n)$
such that $\Dcal_{t,i}=\Dcal^i_t$ for $i=0,1$ and for
all $t\in[a,b]$, and for all $s\in[0,1]$ set:
\[\mathcal Q^s=\big\{v\in\Hcal:v(t)\in\Dcal_{s,t}\
\text{for all}\ t\in[a,b]\big\};\]
using the argument above, $\mathcal Q^s$ is commensurable
with $V^{\scriptscriptstyle-}(\mathcal S_1^\theta)$ for all $s$.

Let us prove that the map $s\mapsto\mathcal Q^s$ is a {\em continuous family
of closed subspaces\/} of $\Hcal$ (see Appendix~\ref{sec:fredholmforms}),
i.e., that  $\mathcal Q^s$ is the image of a fixed closed
subspace $\mathcal Q^*$ of $\Hcal$ via a continuous
family of isomorphisms $\phi_s:\Hcal\to\Hcal$.
For all fixed $S^*\in G_k^-(n)$, the map $\mathrm{GL}(n,\R)\ni U\mapsto U(S^*)\in G_k(n)$
is a smooth fibration; choose a continuous lifting
$[a,b]\times[0,1]\ni(t,s)\mapsto U_{t,s}\in\mathrm{GL}(n,\R)$ of
the map $(t,s)\mapsto\Dcal_{t,s}$, i.e.,
\[\phantom{\qquad\forall\;(t,s)\in[a,b]\times[0,1].}U_{t,s}(S^*)=\Dcal_{t,s},\qquad\forall\;(t,s)\in[a,b]\times[0,1].\]
Finally, define
\[\mathcal Q^*=\big\{v\in\Hcal:v(t)\in S^*\ \text{for all}\ t\in[a,b]\big\},\]
and for all $s\in[0,1]$ let $\phi_s:\Hcal\to\Hcal$ be the isomorphism given by:
\[\phi_{s}(x)(t)=U_{t,s}\big(x(t)\big),\quad x\in\Hcal;\]
clearly, $\phi_s(\mathcal Q^*)=\mathcal Q^s$ for all $s$, and
$s\mapsto\phi_s\in\mathrm{GL}(\Hcal)$ is continuous. This proves
that $s\mapsto\mathcal Q^s$ is a continuous family of closed subspaces
of $\Hcal$, and thus, by Corollary~\ref{thm:dimrelconst},
$\Dim_{\mathcal Q^s}(\Hcal_0)$ is constant.

Once the independence on the choice of the maximal
negative distribution has been established,
to prove equality \eqref{eq:eeqq1} we will now choose a maximal negative distribution $\mathcal
D^-=\{D^-_t\}_{t\in[0,1]}$ which is obtained by the parallel
transport along $\theta$ of a maximal negative subspace of
$T_{\theta(0)}M$; let us also denote by
$\Dcal^+=\{\Dcal^+_t\}_{t\in[0,1]}$ a maximal positive
distribution along $\theta$ which is obtained by the parallel
transport along $\theta$ of a maximal positive subspace of
$T_{\theta(0)}M$. In order to compute the left hand side of
\eqref{eq:eeqq1} we will also have to choose a Hilbert space inner
product in $\Hcal^\theta$; to this aim, a convenient choice is to
set:
\[\langle V,W\rangle=\int_0^1\gr_t\big(V'(t),W'(t)\big)\,\mathrm dt,\qquad V,W\in\Hcal^\theta,\]
where $\gr_t$ is the unique positive definite inner product on
$T_{\theta(t)}M$ for which $\mathcal D_t^+$ and $\mathcal D^-_t$
are orthogonal spaces, that coincides with $g$ on $\mathcal D_t^+$
and with $-g$ on $\mathcal D^-_t$. It is easy to see that, with
such choice, the space $\mathcal Q$ is precisely the negative
eigenspace of $\mathcal S_0^\theta$, and recalling
\eqref{eq:spflreldimbis}, we compute easily:
\begin{equation}\label{eq:p34-a}
\begin{split} \Dim_{\mathcal
Q}\big(V^-(\Scal^\theta_1)\big)&\;=\Dim_{V^-(\Scal^\theta_0)}\big(V^-(\Scal^\theta_1)\big)=-
\Dim_{V^-(\Scal^\theta_1)}\big(V^-(\Scal^\theta_0)\big)\\&\;=\iMorse(\theta)-\Dim\big(\Ker(\Scal^\theta_1)\big).
\end{split}\end{equation} Equality \eqref{eq:p34-b} follows immediately
from \eqref{eq:eeqq1} and \eqref{eq:p34-a}.
\end{proof}

\subsection{Bifurcation of geodesics at a   conjugate instant}
As a further application of our theory, we will discuss briefly
a simple consequence of Proposition~\ref{thm:maslovgeo} obtained using
recent results in bifurcation theory for strongly indefinite
variational problems (see \cite{FitzPejsaRecht}).

Let us recall the definition of bifurcation for a smooth
family of functionals at a common critical point.
Given a   family  $\{f_r\}_{r\in [c,d]}$ of
smooth functionals on some Hilbert space $\Hcal$ depending smoothly
on the parameter $r$, assume that $x=0$ is a critical point
for $f_r$ for all $r\in[c,d]$. Consider the set $C=\big\{(x,r):\mathrm df_r(x)=0\big\}$
endowed with the relative topology of $\Hcal\times[c,d]$ and assume that the
segment $D=\{0\}\times[c,d]$ is entirely contained in $C$, i.e., that $x=0$ is
a critical point for $f_r$ for all $r\in[c,d]$. An instant $r_0\in[c,d]$ is said
to be a {\em bifurcation point\/} for the family $\{f_r\}$ if $(0,r_0)$ is an accumulation
point for $C\setminus D$, i.e., if there exists a sequence $(r_n)_n$ in $[c,d]$
tending to $r_0$ and a sequence $(x_n)_{n\in\N}\subset\Hcal\setminus\{0\}$ tending to $0$
such that $\mathrm df_{r_n}(x_n)=0$ for all $n$. Assume that for some (hence for all)
$r\in[c,d]$, the Hessian $\mathrm d^2f_r(0)$ is a (self-adjoint) Fredholm operator
on $\Hcal$. A sufficient condition for the existence of a bifurcation
instant  in the strongly indefinite case has been proven
recently in \cite{FitzPejsaRecht}: if $\spfl\big(\mathrm d^2f_r(0),[c,d]\big)\ne0$,
then there exists a bifurcation instant for $\{f_r\}_r$ in $[c,d]$.
This result can  be applied to the case of bifurcation of semi-Riemannian
geodesics; let us recall from \cite{geobif} the definition of bifurcation point
along a semi-Riemannian geodesic.
\begin{defin}\label{thm:defbifpoint}
Let $(M,g)$ be a semi-Riemannian manifold, $\theta:[0,1]\to M$ be a geodesic
in $M$ and $t_0\in\left]0,1\right[$. The point $\theta(t_0)$ is said to be a
{\em bifurcation point along $\theta$\/}   if there exists a sequence
$\theta_n:[0,1]\to M$ of geodesics in $M$ and a sequence
$(t_n)_{n\in\N}\subset\left]0,1\right[$
satisfying the following properties:
\begin{enumerate}
\item\label{itm:1defbifpt} $\theta_n(0)=\theta(0)$ for all $n$;
\item\label{itm:2defbifpt} $\theta_n(t_n)=\theta(t_n)$ for all $n$;
\item\label{itm:3defbifpt} $\theta_n\to\theta$ as $n\to\infty$;
\item\label{itm:4defbifpt} $t_n\to t_0$ (and thus $\theta_n(t_n)\to\theta(t_0)$) as $n\to\infty$.
\end{enumerate}
\end{defin}
An immediate application of the inverse function theorem tells us
that bifurcation points occur necessarily at conjugate instants,
however, in the non Riemannian (or causal Lorentz\-ian) case it is not
clear which conjugate points determine bifurcation.
Note that if $\theta(t_0)$ is a bifurcation point along $\theta$, then
the exponential map $\exp_{\theta(0)}$ in not one-to-one
on {\em any\/} neighborhood of $t_0\theta'(0)$ in $T_{\theta(0)}M$.

Conditions for the bifurcation at a nondegenerate conjugate instant
have been discussed in \cite{geobif}; using the theory of partial
signatures it is now an easy game to extend the result to the  possibly degenerate
real-analytic case:
\begin{prop}
\label{thm:geobif}
Let $(M,g)$ be a real-analytic semi-Riemannian manifold,  let
$\theta:[0,1]\to M$ be a geodesic and let $t_0\in\left]0,1\right[$ be
a conjugate instant along $\theta$. If the sum of the odd partial signatures
$\sum_{k\ge1}\sigma_{2k-1}(\theta,t_0)$ is
different from $0$, then $\theta(t_0)$ is a bifurcation point along $\theta$.
\end{prop}
\begin{proof}
By a standard local construction (see \cite[Section~5.1]{geobif} for details), the
geodesic bifurcation problem is cast into a bifurcation problem
for a smooth family of functionals defined in a neighborhood of
$0$ in a {\em fixed\/} Hilbert space $\Hcal$.
The spectral flow of the corresponding path of second variations
is precisely generalized Morse index, whose jumps, by Theorem~\ref{thm:tuttouguale},
occur at those conjugate instants giving a non zero contribution
to the Maslov index, i.e., those conjugate instants $t_0\in\left]0,1\right[$
along $\theta$ such that, for
$\varepsilon>0$ small enough,
$\iMaslov\big(\theta\vert_{[0,t_0-\varepsilon]}\big)\ne
\iMaslov\big(\theta\vert_{[0,t_0+\varepsilon]}\big)$.
The conclusion follows easily from Proposition~\ref{thm:maslovgeo}.
\end{proof}
The result of Proposition~\ref{thm:geobif} gives an important link
between the theory of bifurcation, for which the method of partial
signatures was originally conceived, and the theory of Maslov index
in the context of semi-Riemannian and symplectic geometry.
Moreover, as an easy application of Proposition~\ref{thm:geobif},
we get an extension of a classical result of Morse and Littauer (see Warner's proof
in \cite{War}) that the exponential map of a Riemannian manifold
is never one-to-one on any neighborhood of a conjugate point.
\begin{cor}\label{thm:expnot1to1}
Let $(M,g)$ be a real-analytic semi-Riemannian manifold,
let $\theta:[0,1]\to M$ be a geodesic having a conjugate instant
$t_0\in\left]0,1\right[$ such that $\sum_{k\ge1}\sigma_{2k-1}(\theta,t_0)\ne0$.
Then, the exponential map $\exp_{\theta(0)}$ is not injective on any neighborhood
of $t_0\dot\theta(0)$.

In particular, the result holds true if $(M,g)$ is Riemannian, or
if $(M,g)$ is Lorentzian and $\theta$ is nonspacelike.\qed
\end{cor}
\end{section}

\appendix
\begin{section}[Relative index of Fredholm bilinear forms on
Hilbert spaces]{Relative index of Fredholm bilinear forms\\ on
Hilbert spaces}
\label{sec:fredholmforms}
The  goal of this appendix  is
to provide the reader with a formal proof of a result  (Proposition~\ref{thm:indrelindcoind})
that gives the relative index of a form  as the difference between the
index and the coindex of suitable restrictions  of the form.
A large portion of the material presented is borrowed from \cite[Section~2]{geobif}.

\smallskip

Let $\Hcal$ be a Hilbert space with inner product
$\langle\cdot,\cdot\rangle$; we will denote by
$\mathrm{Lin}(\Hcal)$ the space of all bounded linear operators on
$\Hcal$ and by $\mathrm{GL}(\Hcal)$ the group of invertible
operators in $\mathrm{Lin}(\Hcal)$. Let $B$ a bounded symmetric
bilinear form on $\Hcal$; there exists a unique self-adjoint
bounded operator $S:\Hcal\to \Hcal$ such that $B=\langle
S\cdot,\cdot\rangle$, that will be called the {\em realization of
$B$\/} (with respect to $\langle\cdot,\cdot\rangle$). $B$ is
nondegenerate if its realization is injective, $B$ is strongly
nondegenerate if $S$ is an isomorphism. If $B$ is strongly
nondegenerate, or if more generally $0$ is not an accumulation
point of the spectrum of $S$ (for instance, if $S$ is Fredholm),
we will call the {\em negative eigenspace\/} (resp., {\em the
positive eigenspace}) of $B$ the closed subspace
$V^{\scriptscriptstyle-}(S)$ (resp., $V^{\scriptscriptstyle+}(S)$)
of $\Hcal$ given by $\chi_{\left]-\infty,0\right[}(S)$ (resp.,
$\chi_{\left]0,+\infty\right[}(S)$), where $\chi_I$ denotes the
characteristic function of the interval $I$. The spaces
$V^{\scriptscriptstyle-}(S)$ and $V^{\scriptscriptstyle+}(S)$ are
$S$-invariant, and they are both orthogonal and $B$-orthogonal. We
will say that $B$ is {\em Fredholm\/} if $S$ is Fredholm, or that
$B$ is {\em RCPPI}, {\em realized by a compact perturbation of a
positive isomorphism}, (resp., {\em RCPNI}) if $S$ is of the form
$S=P+K$ (resp., $S=N+K$) where $P$ is a positive isomorphism of
$\Hcal$ ($N$ is a negative isomorphism of $\Hcal$) and $K$ is
compact. The properties of being Fredholm, RCPPI or RCPNI do not
depend on the inner product, although the realization $S$ and the
spaces $V^\pm(S)$ do.

If $B$ is RCPPI (resp., RCPNI), then both its nullity $n_0(B)$ and its index $n^-(B)$ (resp., and its
coindex $n^+(B)$) are finite numbers.
Given a closed subspace $W\subset \Hcal$, the {\em $B$-orthogonal complement of $W$}, denoted by
$W^{\perp_B}$, is the closed subspace of $\Hcal$:
\[W^{\perp_B}=\big\{x\in \Hcal:B(x,y)=0\ \text{for all}\ y\in W\big\}=S^{-1}(W^\perp).\]
If $B$ is  Fredholm, and let $S$ be its realization and $W\subset \Hcal$
is any subspace, then the following properties hold:
\begin{itemize}
\item $B$ is nondegenerate iff it is strongly nondegenerate;
\item $n_0(B)<+\infty$;
\item $(W^{\perp_B})^{\perp_B}=\overline W+\mathrm{Ker}(S)$;
\item if $W$ is closed and $B\vert_W$ (i.e., the restriction of $B$ to $W\times W$)
in nondegenerate, then also $B\vert_{W^{\perp_B}}$ is nondegenerate and $\Hcal=W\oplus W^{\perp_B}$.
\end{itemize}

Let us now recall a few basic things on the notion of commensurability of
closed subspaces (see reference \cite{Abbo1} for more details).
Let $V,W\subset \Hcal$ be closed subspaces and let $P_V$ and
$P_W$ denote the orthogonal projections respectively onto $V$ and $W$. We say that
$V$ and $W$ are {\em commensurable\/} if the restriction
to $V$ of the projection $P_W$ is a Fredholm operator from $V$ to $W$.
It is an easy exercise to show that commensurability is an equivalence
relation in the Grassmannian of all closed subspaces of $\Hcal$;
observe in particular that, identifying each Hilbert space with its own
dual, the adjoint of the operator
$P_W\vert_V:V\to W$ is precisely $P_V\vert_W:W\to V$.
If $V$ and $W$ are commensurable the {\em relative dimension\/}
$\Dim_W(V)$ of $V$ with respect to $W$ is defined as the Fredholm
index $\mathrm{ind}\big(P_V\vert_W:W\to V\big)$, which is equal to:
\[\Dim_W(V)=\mathrm{ind}\big(P_W\vert_V:V\to W\big)=\mathrm{dim}(W^\perp\cap V)-\mathrm{dim}(W\cap V^\perp).\]
Clearly, if $V$ and $W$ are commensurable, then $V^\perp$ and $W^\perp$ are commensurable, and:
\[\Dim_{W^\perp}(V^\perp)=-\Dim_W(V)=\Dim_V(W).\]
The commensurability of closed subspaces and the relative dimension
do {\em not\/} depend on the choice of a Hilbert space inner
product on $\Hcal$.

Using the basic properties of Fredholm index, it
is easy to prove the following:
\begin{prop}\label{thMlsommadim}
Let $V,W\subset\Hcal$ be closed commensurable subspaces, and
let $V',W'$ be finite dimensional subspaces such that $V'\subset V^\perp$
and $W'\subset W^\perp$. Then, $V+V'$ and $W+W'$ are commensurable, and
\[\Dim_{W+W'}(V+V')=\Dim_W(V)+\Dim(V')-\Dim(W').\qed\]
\end{prop}
Let us recall that a family $\{S_t\}_{t\in[a,b]}$ of closed
subspaces of a Hilbert space is said to be {\em continuous\/}
if for all $t_0\in[a,b]$ there exists $\varepsilon>0$, a closed subspace
$S\subseteq H$ and a continuous map $\left]t_0-\varepsilon,
t_0+\varepsilon\right[\ni t\mapsto\phi_t\in\mathrm{GL}(\Hcal)$
such that $\phi_t(S)=S_t$ for all $t$.
Such a map $t\mapsto\phi_t$ will be called a {\em local trivialization\/}
of the family $S_t$ around $t_0$.
The relative dimension
is continuous with respect to this notion of continuity for
paths of closed subspaces. In order to prove this, it is useful
to characterize continuous families  of closed subspaces in terms of graphs
of continuous families of linear operators:
\begin{lem}\label{thm:grafici}
Let $[a,b]\ni t\mapsto S_t$ be a continuous family of closed
subspaces of $\Hcal$. Then, for all $t_0\in[a,b]$ there
exists an orthogonal decomposition $\Hcal=W_0\oplus W_1$
of $\Hcal$ into closed subspaces, and a continuous map
defined around $t_0$, $t\mapsto L_t\in\mathrm{Lin}(W_0,W_1)$
of bounded linear operators from $W_0$ to $W_1$ such that
$S_t=\mathrm{Gr}(L_t)$ for all $t$.
\end{lem}
\begin{proof}
Choose a closed subspace $S\subset\Hcal$, a
local trivialization $t\mapsto\phi_t$ of $S_t$ around $t_0$,
with $\phi_t(S)=S_t$;
set $W_0=S_{t_0}$, $W_1=W_0^\perp$. Denote by $\pi_i:\Hcal\to W_i$ the
orthogonal projection, $i=0,1$, and define:
\[L_t=\big(\pi_1\circ\phi_t\vert_S\big)\circ
\big(\pi_0\circ\phi_t\vert_S\big)^{-1}.\]
Observe that $\pi_0\circ\phi_{t_0}\vert_S:S\to W_0$ is an
isomorphism, and by continuity, $\pi_0\circ\phi_{t}\vert_S:S\to W_0$
is an isomorphism for $t$ near $t_0$.
Clearly, $t\mapsto L_t$ is continuous, and an easy
computation shows that $\mathrm{Gr}(L_t)=S_t$ for all $t$.
\end{proof}
\begin{lem}\label{thm:contproiezioni}
Let $\Hcal=W_0\oplus W_1$ be an orthogonal direct sum.
The map $\mathrm{Lin}(W_0,W_1)\to \mathrm{Lin}(\Hcal)$
given by $L\mapsto P_{\mathrm{Gr}(L)}$ is continuous.
\end{lem}
\begin{proof}
A straightforward computation gives:
\[P_{\mathrm{Gr}(L)}(x_0+x_1)=\big((\mathrm{I}+L^*L)^{-1}(x_0+L^*x_1),
L(\mathrm{I}+L^*L)^{-1}(x_0+L^*x_1)\big),\]
for $x_0\in W_0$, $x_1\in W_1$. The conclusion follows easily.
\end{proof}
\begin{cor}
\label{thm:proiezcontinue}
If $t\mapsto S_t\subset\Hcal$ is a continuous family of closed subspaces,
then the map $t\mapsto P_{S_t}\in\mathrm{Lin}(\Hcal)$ is continuous.
\end{cor}
\begin{proof}
It follows immediately from Lemmas~\ref{thm:grafici}
and \ref{thm:contproiezioni}.
\end{proof}
\begin{cor}\label{thm:dimrelconst}
Let $[a,b]\ni t\mapsto V_t\subset\Hcal$ and
$[c,d]\ni s\mapsto W_s\subset\Hcal$ be continuous families
of closed subspaces, and let $(t_0,s_0)\in[a,b]\times [c,d]$
be such that $V_{t_0}$ is commensurable with $W_{s_0}$.
Then, $V_t$ is commensurable
with $W_s$  and the relative dimension $\Dim_{W_s}(V_t)$ is
constant  for $(s,t)$ near $(s_0,t_0)$; if $V_t$ is commensurable
with $W_s$ for all $s$ and $t$, then $\Dim_{W_s}(V_t)$
is constant on $[a,b]\times[c,d]$.
\end{cor}
\begin{proof}
Choose closed subspaces $V,W\subset\Hcal$ and local trivializations
$\phi_t$ and $\psi_s$ such that $\phi_t(V)=V_t$ and $\psi_s(W)=W_s$
for $(t,s)$ near $(t_0,s_0)$. Since the Fredholm index
is additive by composition, and the Fredholm index of an isomorphism
is $0$, then:
\[\Dim_{W_s}(V_t)=\mathrm{ind}\big(P_{W_s}\vert_{V_t}:V_t\to W_s\big)=\mathrm{ind}
\big(\psi_s^{-1}\circ P_{W_s}
\circ\phi_t\vert_V:V\to W\big).\]
The conclusion follows observing that the set of Fredholm
operators is open, the map $(s,t)\mapsto \psi_s^{-1}\circ P_{W_s}
\circ\phi_t$ is continuous (Corollary~\ref{thm:proiezcontinue}), and the Fredholm index
is locally constant.
\end{proof}
Compact perturbations preserve the commensurability class of
positive and negative eigenspaces of Fredholm operators:
\begin{prop}
\label{thm:era1.6}
Let $S,T$ be linear bounded self-adjoint operators on $\Hcal$ whose difference
$K=S-T$ is compact. Then $V^{\scriptscriptstyle-}(S)$ (resp., $V^{\scriptscriptstyle+}(S)$)
is commensurable with $V^{\scriptscriptstyle-}(T)$ (resp., with $V^{\scriptscriptstyle+}(T)$).
Conversely, assume that $S$ is a bounded self-adjoint Fredholm operator on $\Hcal$, and let
$\Hcal=W^{\scriptscriptstyle-}\oplus W^{\scriptscriptstyle+}$ be an orthogonal decomposition
of $\Hcal$ such that $W^{\scriptscriptstyle-}$ is commensurable with $V^{\scriptscriptstyle-}(S)$
and $W^{\scriptscriptstyle+}$ is commensurable with $V^{\scriptscriptstyle+}(S)$.
Then there exists an invertible self-adjoint operator $T$ on $\Hcal$ such that
$V^{\scriptscriptstyle-}(T)=W^{\scriptscriptstyle-}$,
$V^{\scriptscriptstyle+}(T)=W^{\scriptscriptstyle+}$ and such that $S-T$ is compact.
\end{prop}
\begin{proof}
See \cite[Proposition~2.3.2 and Proposition~2.3.5]{Abbo1}.
\end{proof}
Let us now study under which conditions a closed subspace $W\subset\Hcal$
is commensurable with the negative eigenspace of a Fredholm symmetric bilinear form
$B$ on $\Hcal$.
\begin{lem}\label{thm:WWperpB}
Let $\Hcal$ be a Hilbert space with inner product $\langle\cdot,\cdot\rangle$,
let $S:\Hcal\to\Hcal$ be a self-adjoint Fredholm operator and
set $B=\langle S\cdot,\cdot\rangle$.
Let $W\subset\Hcal$ be a closed subspace with the following
properties:
\begin{itemize}
\item[(a)] $B\vert_W$ is {\em strongly negative definite}, i.e.,
there exists $k>0$ such that $-B(x,x)\ge k\Vert x\Vert^2$ for
all $x\in W$;
\item[(b)] $B\vert_{W^{\perp_B}}$ is positive semi-definite, i.e.,
$B(x,x)\ge0$ for all $x\in W^{\perp_B}$.
\end{itemize}
Then $W$ is commensurable with $V^{\scriptscriptstyle-}(S)$,
and $\Dim_W\big(V^{\scriptscriptstyle-}(S)\big)=0$.
\end{lem}
\begin{proof}
Set $H_-=V^{\scriptscriptstyle-}(S)$ and $H_+=\Ker(S)\oplus
V^{\scriptscriptstyle+}(S)$, so that $\Hcal=H_-\oplus H_+$, the
direct sum being orthogonal and also $B$-orthogonal; let $P_-:\Hcal\to H_-$ be
the orthogonal projection. Observe that $H_-$ and $H_+$ are $S$-invariant;
assumption (a) means that the restriction $S_-\vert_{H_-}:H_-\to H_-$
is a negative isomorphism, i.e., its spectrum is contained
in $\left]-\infty,0\right[$; moreover, by (a), $W\cap H_+=\{0\}$,
and, by (b), $W^{\perp_B}\cap H_-=\{0\}$.
The thesis is equivalent to
the condition that the restriction $P_-\vert_{W}:W\to H_-$ be
an isomorphism, which is what we will prove now.

In first place, observe that $P_-\vert_W$ is injective, because
$\Ker\big(P_-\vert_W\big)=H_+\cap W=\{0\}$.
Next, $P_-\vert_W$ has dense image, because its adjoint
$P_W\vert_{H_-}$ is injective: $\Ker\big(P_W\vert_{H_-}\big)
=W^\perp\cap H_-=\{0\}$. For, if $x\in W^\perp\cap H_-$, since
the restriction $S\vert_{H_-}:H_-\to H_-$ is an isomorphism,
then $x=Sy$ for some $y\in H_-$, and for all $w\in W$ it is
$B(y,w)=\langle Sy,w\rangle=\langle x,w\rangle=0$, so that
$y\in W^{\perp_B}\cap H_-=\{0\}$, i.e., $x=0$.

Finally, we must prove that $P_-\vert_W$ has closed image;
it suffices to show that there exists $\beta>0$ such that
$\Vert P_-(x)\Vert^2\ge\beta\Vert x\Vert^2$ for all $x\in W$.
To this aim, define:
\[k=\inf_{y\in W, \Vert y\Vert=1}-B(y,y)>0,\]
and let $x\in W$ be fixed with $\Vert x\Vert=1$, $x=x_-+x_+$,
$x_-\in H_-$, $x_+\in H_+$. Then,
\[B(x_-,x_-)\le B(x_-,x_-)+B(x_+,x_+)=B(x,x)\le -k;\]
moreover, since $B$ is strongly negative definite on $H_-$, then
$-B$ is a Hilbert space inner product on $H_-$ equivalent to
$\langle\cdot,\cdot\rangle$, hence there exists $\alpha>0$ such that:
\[-B(y,y)\le\alpha\Vert y\Vert^2,\quad\forall\,y\in H_-.\]
From the last two inequalities we obtain:
\[\Vert P_-(x)\Vert^2=\Vert x_-\Vert^2\ge\frac k\alpha>0,\]
which concludes the proof.
\end{proof}

\begin{lem}\label{thm:RCP-comm}
Let $B$ be a Fredholm symmetric bilinear form on the Hilbert space $\Hcal$
and let  $W\subset \Hcal$ be a closed
subspace.
Then, the following are equivalent:
\begin{itemize}
\item[(a)]  $B\vert_W$ is RCPNI and $B\vert_{W^{\perp_B}}$ is RCPPI;
\item[(b)] $W$ is commensurable with $V^{\scriptscriptstyle-}(S)$, where
$S$ is the realization of $B$.
\end{itemize}
In particular, condition (b) is independent on the choice of an inner product on $\Hcal$.
\end{lem}
\begin{proof}
Assume that $W$ is commensurable with $V^{\scriptscriptstyle-}(S)$; then $W^\perp$  is commensurable with
$V^{\scriptscriptstyle-}(S)^\perp= V^{\scriptscriptstyle+}(S)\oplus\Ker(B)$. Moreover, since
$\Ker(B)$ is finite dimensional, then $W^\perp$ is also commensurable with
$V^{\scriptscriptstyle+}(S)$. By Proposition~\ref{thm:era1.6}, there exists an invertible self-adjoint
operator $T:\Hcal\to \Hcal$ such that $V^{\scriptscriptstyle-}(T)=W$, $V^{\scriptscriptstyle+}(T)=W^\perp$,
and with $S=T+K$, with $K$ compact. It follows easily that $B\vert_W$ is RCPNI
(namely, if $P$ denotes the orthogonal projection onto $W$, the realization
of $B\vert_W$ is $PS\vert_W=(PT+PK)\vert_W=(T+PK)\vert_W$), and $B\vert_{W^\perp}$ is RCPPI.
Observe in particular that $W\cap W^{\perp_B}=\Ker(B\vert_W)$ is finite dimensional.
To prove that $B\vert_{W^{\perp_B}}$ is RCPPI we argue as follows; denote by $P$
the orthogonal projection onto $W$ and by $P^{\scriptscriptstyle\perp}=1-P$ the orthogonal
projection onto $W^\perp$. As we
have observed, $W^{\perp_B}=S^{-1}(W^\perp)$; hence,
for all $x,y\in W^{\perp_B}$ we have:
\begin{equation}\label{eq:era6}\begin{split}B(x,y)&\,=\langle Sx,y\rangle=\langle
Sx,P^{\scriptscriptstyle\perp} y\rangle=\langle SPx,P^{\scriptscriptstyle\perp} y\rangle+\langle
SP^{\scriptscriptstyle\perp} x,P^{\scriptscriptstyle\perp} y\rangle=\\ &\,=\langle
P^{\scriptscriptstyle\perp}KPx,y\rangle +\langle
P^{\scriptscriptstyle\perp}TP^{\scriptscriptstyle\perp}x,y\rangle+
\langle P^{\scriptscriptstyle\perp}KP^{\scriptscriptstyle\perp}x,y\rangle.\end{split}
\end{equation}
In the above equality we have used the fact that $W$ and $W^\perp$ are $T$-invariant.
From \eqref{eq:era6} we deduce that $B\vert_{W^{\perp_B}}$ is represented by
a compact perturbation of the operator
$\widetilde T:W^{\perp_B}\to W^{\perp_B}$ given by $\widetilde
T=P^{\scriptscriptstyle\perp_B}P^{\scriptscriptstyle\perp}TP^{\scriptscriptstyle\perp}
\vert_{W^{\perp_B}}$ (where $P^{\scriptscriptstyle\perp_B}$ is the orthogonal
projection onto $W^{\perp_B} $) which is positive
semi-definite. The kernel of $\widetilde T$ is easily computed as the finite
dimensional space $W^{\perp_B}\cap T^{-1}\big(W\cap W^{\perp_B}\big)$; it follows that $\widetilde
T$ is a compact perturbation of a positive isomorphism of $W^{\perp_B}$,
which proves that (b) implies (a).

Conversely, assume that $B\vert_W$ is RCPNI and $B\vert_{W^{\perp_B}}$ is RCPPI.
Using functional calculus, write $W=W_0\oplus W_1$, with $W_0$ and $W_1$ orthogonal and
$B$-orthogonal, $B\vert_{W_0}$ strongly negative definite, $B\vert_{W_1}$
positive semi-definite, and $\Dim(W_1)<+\infty$; then,
\[\Hcal=W_0\oplus W_0^{\perp_B}.\]
Clearly,  $W^{\perp_B}\subset W_0^{\perp_B}$; we claim that
$W^{\perp_B}$ has finite codimension in $W_0^{\perp_B}$. Namely,
\[W_0^{\perp_B}\cap W_1^{\perp_B}\subset W^{\perp_B}\subset W_0^{\perp_B},\]
and since $W_1^{\perp_B}$ has finite codimension in $\Hcal$,
$W_0^{\perp_B}\cap W_1^{\perp_B}$ (and {\em a fortiori\/} $W^{\perp_B}$)
has finite codimension\footnote{If $X$ is a vector space, and
$S,Y\subset X$ subspaces. If $S$ has finite codimension in $X$, then
$S\cap Y$ has finite codimension in $Y$. Namely, $\mathrm{codim}_Y(S\cap Y)$
equals the dimension of the image of $Y$ by the projection $X\to X/S$
onto the finite dimensional space $X/S$.}
in $W_0^{\perp_B}$.
Since $B\vert_{W^{\perp_B}}$ is RCPP, it follows that also
$B\vert_{W_0^{\perp_B}}$ is RCPPI, and again we can write $W_0^{\perp_B}$
as a $B$-orthogonal direct sum $W_2\oplus W_3$ with $B\vert_{W_2}$
negative definite, $B\vert_{W_3}$ positive semi-definite, and $\Dim(W_2)<+\infty$.

Set $Z:=W_0\oplus W_2$; then $B\vert_Z$ is strongly negative definite;
moreover, $Z^{\perp_B}=W_3$, hence $B\vert_{Z^{\perp_B}}$ is positive
semi-definite; by Lemma~\ref{thm:WWperpB},
$Z$ is commensurable with $V^{\scriptscriptstyle-}(S)$
Since $W_0$ has finite codimension in $Z$, then $W_0$ and $Z$ are
commensurable, and since $W_0$ has finite codimension in $W$,
then $W_0$ and $W$ are commensurable. By transitivity, $W$ is commensurable
with $V^{\scriptscriptstyle-}(S)$, and the proof is concluded.
\end{proof}

Assume now that $B$ is a symmetric bilinear form, $S$ is its realization; if $W$ is closed
subspace of $\Hcal$ which is commensurable with $V^{\scriptscriptstyle-}(S)$, the one defines the {\em
relative index\/} of $B$ with respect to $W$, denoted by $\ind_W(B)$, the integer number:
\[\ind_W(B)=\mathrm{dim}_W\big(V^{\scriptscriptstyle-}(S)\big).\]

A subspace $Z$ of $\Hcal$ is said to be {\em isotropic\/} for the symmetric bilinear form $B$ if
$B\vert_Z\equiv0$.
\begin{lem}\label{thm:dimisotropic}
Let $B$ be a RCPPI nondegenerate symmetric bilinear form on $\Hcal$, and let $Z\subset \Hcal$ be an isotropic
subspace of
$B$. Then:
\[n^-(B)=n^-\big(B\vert_{Z^{\perp_B}}\big)+\mathrm{dim}(Z).\]
\end{lem}
\begin{proof}
Since $B$ is RCPPI, then the index $n^-(B)$ is finite, and so
$n^-\big(B\vert_{Z^{\perp_B}}\big)$ and $\Dim(Z)$ are finite.
Clearly, $Z\subset Z^{\perp_B}$; let $U\subset Z^{\perp_B}$ be a closed subspace
such that $Z^{\perp_B}=Z\oplus U$, so that $B\vert_U$ is nondegenerate and
$\Hcal=U\oplus U^{\perp_B}$. Moreover:
\[n^-(B)=n^-\big(B\vert_U\big)+n^-\big(B\vert_{U^{\perp_B}}\big).\]
Since $Z$ is isotropic, then $n^-\big(B\vert_U\big)=n^-\big(B\vert_{Z^{\perp_B}}\big)$;
to conclude the proof we need to show that $n^-\big(B\vert_{U^{\perp_B}}\big)=
\Dim(Z)$. To this aim, observe first that $\Dim(U^{\perp_B})=2\Dim(Z)$.
Namely, $\Dim(U^{\perp_B})=\Codim(U)$; moreover, $\Codim_{Z^{\perp_B}}(U)=
\Dim(Z)$, and $\Codim(Z^{\perp_B})=\Dim(Z)$.
Thus, keeping in mind that the dimension of an isotropic subspace of a nondegenerate
symmetric bilinear form is
less than or equal to the index and the coindex, we have:
\[n^-\big(B\vert_{U^{\perp_B}}\big)+n^+\big(B\vert_{U^{\perp_B}}\big)=\Dim(
U^{\perp_B})=2\,\Dim(Z)\le n^-\big(B\vert_{U^{\perp_B}}\big)+n^+\big(B\vert_{U^{\perp_B}}\big),
\]
which proves that $n^-\big(B\vert_{U^{\perp_B}}\big)=n^+\big(B\vert_{U^{\perp_B}}\big)=
\Dim(Z)$ and concludes the proof.
\end{proof}

\begin{lem}\label{thm:indtilde}
Let $B$ be a nondegenerate Fredholm symmetric bilinear form on $\Hcal$ and $W\subset \Hcal$ be a closed
subspace
such that $B\vert_{W^{\perp_B}}$ is RCPPI. Let $\widetilde W$ be any closed complement
of $W\cap W^{\perp_B}$ in $W$. Then the following identity holds:
\[n^-\big(B\vert_{\widetilde W^{\perp_B}}\big)=n^-\big(B\vert_{
W^{\perp_B}}\big)+\mathrm{dim}\big(W\cap W^{\perp_B}\big).\]
\end{lem}
\begin{proof}
We start with the observation that $\Ker\big(B\vert_{W}\big)=\Ker\big(B\vert_{W^{\perp_B}}\big)=W\cap
W^{\perp_B}$; this implies in particular that $B\vert_{\widetilde W}$  and $B\vert_{\widetilde
W^{\perp_B}}$ are nondegenerate.
Since $B\vert_{W^{\perp_B}}$ is RCPPI, then $n^-\big(B\vert_{W^{\perp_B}}\big)
$ and $\Dim(W\cap W^{\perp_B})=n$ are finite numbers.

Since $\Codim_{\widetilde W^{\perp_B}}\big(W^{\perp_B}\big)=n$, then:
\[n^-\big(B\vert_{\widetilde W^{\perp_B}}\big)\le n^-\big(B\vert_{W^{\perp_B}}\big)+n,\]
from which it follows that $ n^-\big(B\vert_{\widetilde W^{\perp_B}}\big)$ is finite;
moreover,
$B\vert_{\widetilde W^{\perp_B}}$ is RCPPI. The conclusion now follows
easily from
Lemma~\ref{thm:dimisotropic}, applied to the nondegenerate bilinear form $B\vert_{\widetilde
W^{\perp_B}}$ and the isotropic space $Z=W\cap W^{\perp_B}$.
\end{proof}

\begin{prop}\label{thm:indrelindcoind}
Let $B$ be a Fredholm symmetric bilinear form on $\Hcal$, $S$ its realization and let $W\subset \Hcal$ be
a closed subspace which is commensurable with $V^{\scriptscriptstyle-}(S)$. Then the relative index
$\ind_W(B)$ is given by:
\begin{equation}\label{eq:era7}
\ind_W(B)=n^-\big(B\vert_{W^{\perp_B}}\big)-n^+\big(B\vert_W\big).
\end{equation}
\end{prop}
\begin{proof}
Assume first that $B$ is nondegenerate on $W$; then we
have a direct sum decomposition $\Hcal=W\oplus W^{\perp_B}$. The relative
$\ind_{W}(B)$ does not change if we change the inner product of
$\Hcal$; we can therefore assume that $W$ and $W^{\perp_B}$ are orthogonal
subspaces of $\Hcal$. Then, $S=S^{\scriptscriptstyle-}\oplus S^{\scriptscriptstyle+}$,
where $S^{\scriptscriptstyle-}:W\to W$ is
the realization of $B\vert_{W}$
and $S^{\scriptscriptstyle+}:W^{\perp_B}\to W^{\perp_B}$ is the realization of
$B\vert_{W^{\perp_B}}$. Moreover,
$V^{\scriptscriptstyle-}(S)=V^{\scriptscriptstyle-}(S^{\scriptscriptstyle-})\oplus
V^{\scriptscriptstyle-}(S^{\scriptscriptstyle+})$. An immediate calculation  yields:
\[\begin{split} \ind_{W}(B)\,&=
\Dim\big(V^{\scriptscriptstyle-}(S)\cap W^{\perp_B}\big)-\Dim\big(
V^{\scriptscriptstyle-}(S)^\perp\cap W\big)
\\ &=\Dim\big(V^{\scriptscriptstyle-}(S)\cap
W^{\perp_B}\big)-\Codim_{W}\big(V^{\scriptscriptstyle-}(S^{\scriptscriptstyle-})\big) \\
&=\Dim\big(V^{\scriptscriptstyle-}(S^{\scriptscriptstyle+})\big)
-\Codim_{W}\big(V^{\scriptscriptstyle-}(S^{\scriptscriptstyle-})\big) \\
&=n^-\big(B\vert_{W^{\perp_B}}\big)-n^+\big(B\vert_{W}\big).\end{split}\]

Let us consider now the case that $B\vert_W$ is degenerate; by Lemma~\ref{thm:RCP-comm},
$B\vert_W$ is RCPNI, and so  $\Dim\big(W\cap W^{\perp_B}\big)=n<+\infty$.
Set $\widetilde W=\big(W\cap W^{\perp_B}\big)^\perp\cap W$, so that $B\vert_{\widetilde W}$
is nondegenerate; moreover, $V^{\scriptscriptstyle-}(S)$ is commensurable
with $\widetilde W$, because it has finite codimension in $W$.
We can then apply the first part of the proof, and we obtain:
\begin{equation}\label{eq:era8}\ind_{\widetilde W}(B)=n^-\big(B\vert_{\widetilde
W^{\perp_B}}\big)-n^+\big(B\vert_{
\widetilde W}\big).\end{equation}
Clearly,
\begin{equation}\label{eq:era9}
n^+\big(B\vert_{
\widetilde W}\big)=n^+\big(B\vert_{W}\big);\end{equation}
moreover, by definition of relative index:
\begin{equation}\label{eq:era10}
\ind_{\widetilde W}(B)=\ind_{W}(B)+n.\end{equation}
Finally, by Lemma~\ref{thm:RCP-comm}, $B\vert_{W^{\perp_B}}$ is RCPPI, and by Lemma~\ref{thm:indtilde}:
\begin{equation}\label{eq:era11}n^-\big(B\vert_{\widetilde
W^{\perp_B}}\big)=n^-\big(B\vert_{W^{\perp_B}}\big)+n.\end{equation}
Formulas \eqref{eq:era8}, \eqref{eq:era9}, \eqref{eq:era10} and \eqref{eq:era11}
yield \eqref{eq:era7} and conclude the proof.
\end{proof}

\end{section}

\begin{section}{Connectedness of the Grassmannian of $g$-negative subspaces.}
\label{sec:grasconn}
We will prove briefly in this appendix that, given a nondegenerate
bilinear form $g$ on $\R^n$ having index $k$
and denoting by $G_k^-(n)$ the open
subset  of the Grassmannian $G_k(n)$ consisting
of those $k$-dimensional planes in $\R^n$ on which $g$ is negative
definite, then $G_k^-(n)$ is (arc) connected.
The non trivial case is when $0<k<n$.
Clearly, it is not restrictive
to assume, as we will, that $g$ is the bilinear form whose matrix representation
in the canonical basis of $\R^n$ is given by:
\[g=\begin{pmatrix}-\mathrm I_k&0\\ 0&\mathrm{I}_{n-k}\end{pmatrix};\]
the group of isometries of $g$ will be denoted by  $\mathrm{O}(n,k)$:
\[\mathrm O(n,k)=\big\{A\in\mathrm{GL}(n,\R):A^*gA=g\big\}=
\big\{A\in\mathrm{GL}(n,\R):gA^*g=A^{-1}\big\},\]
and its Lie algebra by $\mathrm{so}(n,k)$:
\[\mathrm{so}(n,k)=\big\{H\in\mathrm{gl}(n,\R): gH+H^*g=0\big\}.\]

To prove our assertion we will show that the connected component of
the identity of $\mathrm{O}(n,k)$ acts transitively on $G^-_k(n)$.
Given $A\in\mathrm O(k,n)$, denote by $A_{\textrm{up}}$ and
$A_{\textrm{low}}$ respectively the upper left $k\times k$ block of $A$
and the lower right $(n-k)\times(n-k)$ block of $A$; it is easy to see
that both $A_{\textrm{up}}$ and $A_{\textrm{low}}$ are invertible,
and we will show that the sign of the determinant of the two blocks distinguish
the connected components of $\mathrm O(k,n)$:
\begin{lem}\label{thm:conncompOnk}
If $0<k<n$, $\mathrm O(k,n)$ has four connected components.
The connected component of the identity consists of those $A\in\mathrm O(k,n)$
such that $\mathrm{det}(A_{\textrm{up}})$ and $\mathrm{det}(A_{\textrm{low}})$
are positive.
\end{lem}
\begin{proof}
The four components are determined by the choice of the
signs of $\mathrm{det}(A_{\textrm{up}})$ and $\mathrm{det}(A_{\textrm{low}})$.
These sets are clearly open; let us prove that  they are connected.
Denote by by $\mathrm{Sym}(n)$ the vector space of $n\times n$
symmetric matrices, and by $\mathrm{Sym}_+(n)$ the subset of
the positive definite ones;
recall that the exponential map gives a diffeomorphism
$\exp:\mathrm{Sym}(n)\to\mathrm{Sym}_+(n)$.
Consider the diffeomorphism:
\begin{equation}\label{eq:diffeodecpol}
\mathrm{GL}(n,\R)\stackrel{\cong}{\longrightarrow}\mathrm O(n)\times\mathrm{Sym}_+(n)
\end{equation}
given by the polar decomposition $A\mapsto (U,P)$, $U=A\vert A\vert^{-1}$,
$P=\vert A\vert$.  It is easy to see that, given $A\in\mathrm{GL}(n,\R)$, $A$ belongs
to $\mathrm O(k,n)$ if and only if $U(A)$ and $P(A)$ do.
Namely, using the fact that $g^2=\mathrm I$,
if $gA^*g=A^{-1}$, then $gPU^*g=(gPg)(gU^*g)=P^{-1}U^{-1}$
and so $(gUg)(gP^{-1}g)=UP=A$. By the uniqueness of the polar decomposition,
$gUg=U$ and $gP^{-1}g=P$, i.e., $U,P\in\mathrm O(n,k)$.

Observe the following:
\begin{enumerate}
\item\label{itm:gruppi1} $\mathrm O(n,k)\cap\mathrm O(n)=\mathrm O(k)\times\mathrm{O}(n-k)$;
\item\label{itm:gruppi2} from the uniqueness of the square root, it follows easily that
if $A\in\mathrm O(n,k)\cap\mathrm{Sym}_+(n)$, then
$A^{\frac12}\in\mathrm O(n,k)$;
\item\label{itm:gruppi3} if $H\in\mathrm{Sym}(n)$, then $H\in\mathrm{so}(n,k)$ if and only if
$\exp(H)\in\mathrm O(n,k)$.\footnote{%
\eqref{itm:gruppi3} follows easily from \eqref{itm:gruppi2}.
Namely, if $\exp(H)\in\mathrm O(n,k)$, then $\exp(qH)\in\mathrm O(n,k)$ for
all $q\in\Q$, and since $\mathrm O(n,k)$ is closed, also $\exp(tH)\in\mathrm O(n,k)$
for all $t\in\R$. Differentiating $t\mapsto\exp(tH)$ at $t=0$ we get
that $H\in\mathrm{so}(n,k)$. The converse is trivial.}
\end{enumerate}
Thus, the restriction of the diffeomorphism \eqref{eq:diffeodecpol}
is a diffeomorphism:
\[\big(\mathrm O(k)\times\mathrm O(n-k)\big)\times\big(
\mathrm{so}(n,k)\cap\mathrm{Sym}(n)\big)\stackrel{\cong}{\longrightarrow}\mathrm O(n,k),\]
given by $(U,Z)\mapsto U\exp(Z)$.
The conclusion follows easily from the observation that
$\mathrm{so}(n,k)\cap\mathrm{Sym}(n)$ is a vector space, hence
contractible, and from the fact that, for $r>0$, $\mathrm O(r)$ has exactly two connected
components determined by the sign of the determinant.
\end{proof}

\begin{cor}\label{thm:Gnk-con}
The connected component of the identity of $\mathrm O(n,k)$ acts transitively
on $G_k^-(n)$; in particular, $G_k^-(n)$ is arc-connected.
\end{cor}
\begin{proof}
Let $S\in G_k^-(n)$ be fixed arbitrarily; let us show that there exists
$T$ in the connected component of the identity of $\mathrm O(n,k)$
such that $T\big(\R^k\oplus\{0\}\big)=S$. Set $S'=S^\perp_g$,
choose a $g$-orthonormal basis $b_1,\ldots,b_k$ of $S$,
a $g$-orthonormal basis $b_{k+1},\ldots,b_n$ if
$S'$ and, denoting by $e_1,\ldots,e_n$ the canonical basis of
$\R^n$, let
$\widetilde T\in\mathrm{GL}(n)$ be such that $\widetilde T(e_i)=b_i$
for all $i=1,\ldots,n$.
Clearly, $\widetilde T\in\mathrm O(n,k)$, because it send a $g$-orthonormal
basis into another $g$-orthonormal basis; moreover, $
\widetilde T\big(\R^k\oplus\{0\}\big)=S$.
Replacing $b_1$ with $-b_1$ in the choice of an orthonormal basis
of $S$ has the effect of changing the sign of the determinant of the
upper left block $k\times k$ of the matrix representation of $\widetilde T$
in the canonical basis, while replacing $b_{k+1}$ with $-b_{k+1}$
has the effect of changing the sign of the lower right $(n-k)\times (n-k)$
block of such matrix. In conclusion, the appropriate choice for the sign of
$b_1$ and $b_{k+1}$ can be made to ensure that the corresponding
operator $\widetilde T$ belongs to the connected component
of the identity of $\mathrm O(n,k)$, and this concludes the proof.
\end{proof}

\end{section}



\begin{thebibliography}{99}


\bibitem{Abbo1} A. Abbondandolo, {\em Morse Theory for Hamiltonian Systems}, Pitman Research Notes in
Mathematics, vol.\ 425, Chapman \& Hall, London, 2001.

\bibitem{abbomajer} A. Abbondandolo, P. Majer, {\em Morse Homology on Hilbert Spaces}, Commun.\ Pure
Appl.\ Math.\ {\bf64} (2001), 689--760.

\bibitem{APS} M.\ F.\ Atiyah, V.\ Patodi, I.\ M.\ Singer,
{\em Spectral Asymmetry and Riemannian
Geometry I, II, III}, Proc.\ Cambr.\ Phil.\ Soc.\
{\bf77}  (1975), 43--69,  {\bf78}  (1975), no.\ 3, 405--432,
and {\bf79} (1976), no.\ 1, 71--99.

\bibitem{BEE} J.\ K.\ Beem, P.\ E.\ Ehrlich, K.\ L.\ Easley,
{\em Global Lorentzian Geometry},
Marcel Dekker, Inc., New York and Basel, 1996.

\bibitem{BGM} V.\ Benci, F.\ Giannoni, A.\ Masiello, {\em
Some Properties of the Spectral Flow in Semiriemannian Geometry},
J.\ Geom.\ Phys.\ {\bf27} (1998), 267--280.

\bibitem{BoosFuru} B.\ Booss-Bavnek, K.\ Furutani,
{\em The Maslov Index: a Functional Analytical Definition
and the Spectral Flow Formula}, Tokio J.\ Math.\ {\bf21}
(1998), No.\ 1, 1--34.

\bibitem{BoosFurub} B.\ Booss-Bavnek, K.\ Furutani,
{\em Symplectic Functional Analysis and Spectral Invariants},
Contemporary Mathematics vol.\ {\bf242} (1999), 53--83.

\bibitem{FurOts1} B. Booss--Bavnek, K. Furutani, N. Otsuki, {\em Criss-cross Reduction of the
Maslov Index and a proof of Yoshida--Nicolaescu Theorem},
Tokio J.\ Math.\ {\bf 24} (2001), no.\ 1, 113--128.


\bibitem{CapLeeMil} S. E. Cappell, R. Lee, E. Y. Miller, {\em On the Maslov index},
 Comm.\ Pure Appl.\ Math.\  {\bf 47}  (1994),  no.\ 2, 121--186.

\bibitem{CapLeeMil2} S. E. Cappell, R. Lee, E. Y. Miller,
{\em Self-adjoint elliptic operators and manifold decompositions, Part I: Low
eigenmodes and stretching},
 Comm.\ Pure Appl.\ Math.\  {\bf 49}  (1996),  825--866.

\bibitem{CapLeeMil3} S. E. Cappell, R. Lee, E. Y. Miller,
{\em Self-adjoint elliptic operators and manifold decompositions, Part II:
Spectral Flow and Maslov index},
 Comm.\ Pure Appl.\ Math.\  {\bf 49}  (1996), 869--909.

\bibitem{ConZeh} C. C. Conley, E. Zehnder, {\em Morse-type index theory for
flows and periodic solutions of Hamiltonians}, Commun.\ Pure Appl.\
Math.\ {\bf37} (1984), 207--253.

\bibitem{CPPT} J.\ C.\ Corr\^ea Eidam, A.\ L.\ Pereira,
P.\ Piccione, D.\ V.\ Tausk, {\sl On the Equality between the Maslov Index and the Spectral
Index for the semi-Riemannian Jacobi Operator},
J.\  Math.\ Anal.\  Appl.\ {\bf268}, No.\ 2 (2002), 564--589.

\bibitem{Dan} M.\ Daniel, {\em An Extension of a Theorem of Nicolaescu
on Spectral Flow and Maslov Index}, Proc.\ Amer.\ Math.\ Soc.\ {\bf 128}, No.\ 2 (1999),
611--619.

\bibitem{DanKirk} M. Daniel, P. Kirk, {\em A general splitting formula
for the spectral flow}, Michigan Math.\ J.\ {\bf46} (1999), no.\ 3, 589--617.

\bibitem{deGos1} M.\ de Gosson, {\em La relation entre $\mathrm{Sp}_\infty$,
rev\^etement universel du groupe symplectique, et $\mathrm{Sp}\times\Z$},
C.\ R.\ Acad.\ Sci.\ Paris, t.\ 310, S\'erie I,
(1990), 245--248.

\bibitem{deGos2} M.\ de Gosson, {\em Le d\'efinition de l'indice de Maslov sans
hypoth\`ese de transversalit\'e}, C.\ R.\ Acad.\ Sci.\ Paris, t.\ 310, S\'erie I,
(1990), 279--282.

\bibitem{deGos3} M.\ de Gosson, {\em The structure of $q$-symplectic geometry},
J.\ Math.\ Pures Appl.\ {\bf71} (1992), 429--453.

\bibitem{deGos4} M.\ de Gosson, {\em Lagrangian Path Intersections and the Leray Index},
Contemporary Mathematics, Vol.\ 258 (2000), 177--184.

\bibitem{deGosdeGos} M.\ de Gosson, S. de Gosson,
{\em Symplectic Path Intersections and the Leray Index}, preprint.

\bibitem{Duis} J.\ J.\ Duistermaat, {\em On the Morse Index in Variational Calculus},
Adv.\ in Math.\ {\bf 21} (1976), 173--195.

\bibitem{DunSch} N.\ Dunford, J.\ T.\ Schwarz, {\em Linear Operators\/} vol.\ III,
Pure and Applied Mathematics vol.\ VII, Wiley Interscience, 1971.

\bibitem{FarLev} M.\ S.\ Farber, J.\ P.\ Levine, {\em Jumps of
the eta-invariant}, Math.\ Z.\ {\bf223} (1996), no.\ 2,
197--246.

\bibitem{Flo1} A. Floer, {\em A Relative Morse Index for the Symplectic Action},
Commun.\ Pure Appl.\ Math.\ {\bf 41} (1988), 393--407.

\bibitem{Flo2} A. Floer, {\em The unregularized gradient flow of the
symplectic action},
Commun.\ Pure Appl.\ Math.\ {\bf 41} (1988), 775--813.

\bibitem{Flo3} A. Floer, {\em Morse Theory for Lagrangian Intersections},
J.\ Diff.\   Geom.\   {\bf 28} (1988), 513--547.

\bibitem{FurOts2} K. Furutani, N. Otsuki, {\em Maslov index in the infinite dimension
and a splitting formula}, Japan J.\ Math.\ {\bf28} (2002), no.\ 2, 215--243.

\bibitem{FitzPejsaRecht} P.\ M.\ Fitzpatrick, J.\ Pejsachowicz, L.\ Recht, {\em Spectral
Flow and Bifurcation of Strongly Indefinite Functionals Part I. General Theory}, J.\ Funct.\ Anal.\
{\bf 162} (1) (1999), 52--95.

\bibitem{GiaPicPorCOMPTES} R.\ Giamb\`o,   P.\ Piccione, A.\ Portaluri, {\em Computation
of the Maslov index and the spectral flow  via partial
signatures},  C.\ R.\ Math.\ Acad.\ Sci.\ Paris \textbf{338} (2004), no.\ 5, 397--402.

\bibitem{asian} F.\ Giannoni, A.\ Masiello, P.\ Piccione,
D.\ Tausk, {\em A Generalized Index Theorem for Morse--Sturm Systems
and Applications to semi-Riemannian Geometry}, Asian Journal of
Mathematics Vol.\ 5, no.\ 3 (2001).

\bibitem{GohSig} I. C. Gohberg, E. I. Sigal, {\em An operator
generalization of the logarithmic residue theorem and the theorem
of Rouch\'e}, Math.\ USSR Sbornik, vol.\ 13, No.\ 4 (1971),
603--625.

\bibitem{Hel1} A.\ D.\ Helfer, {\em Conjugate Points on Spacelike Geodesics
or Pseudo-Self-Adjoint Morse-Sturm-Liouville Systems}, Pacific J.\ Math.\
{\bf164}, n.\ 2 (1994), 321--340.

\bibitem{Hor} L. H\"ormander, {\em Fourier integral operators}, Acta Math.\ {\bf127} (1971),
79--183.

\bibitem{kato} T.\ Kato, {\em Perturbation Theory for
Linear Operators}, Grundlehren der Mathematischen
Wissenschaften, vol.\ 132, Springer--Verlag, New York/Berlin,
1980.


\bibitem{KirKla} P.\ Kirk, E.\ Klassen, {\em The spectral flow of
the odd signature operator and higher Massey products}, Math.\
Proc.\ Cambridge Philos.\ Soc.\ {\bf 121} (1997), no.\ 2,
297--320.

\bibitem{LioVer} G. Lions, M. Vergne, {\em The Weil representation, Maslov index
and theta series}, Progress in Mathematics No.\ 6, Birk\"auser, Boston--Basel,
1980.

\bibitem{MasPic} A.\ Masiello, P.\ Piccione, {\em On the
Spectral Flow in Lorentzian Manifolds},
Annali di Matematica Pura e Applicata \textbf{182} (2003), 81--101.

\bibitem{pacific} F.\ Mercuri, P.\ Piccione, D.\ Tausk, {\em
Stability of the Conjugate Index, Degenerate Conjugate Points and the Maslov Index in semi-Riemannian
Geometry},  Pacific J.\ Math.\ {\bf206}  (2002),  no.\ 2,
375--400.

\bibitem{morvan} J.-M.\ Morvan, {\em Maslov, Duistermaat, Conley--Zehnder
invariants in Riemannian Geometry}, Geometry and Topology of Submanifold, V
(Leuven/Brussels, 1992), 174--200, World Sci.\ Publishing, River Edge, NJ, 1993.

\bibitem{ON} B.\ O'Neill, {\em Semi-Riemannian Geometry with Applications
to Relativity}, Academic Press, New York, 1983.

\bibitem{Nico} L. I. Nicolaescu, {\em The Maslov Index, the Spectral Flow, and
Decomposition of Manifolds}, Duke Math.\ J.\ {\bf80} (1995), 485--533.

\bibitem{Phillips} J.\ Phillips, {\em Self-adjoint Fredholm Operators and Spectral Flow},
Canad.\ Math.\ Bull.\ {\bf39} (4) (1996), 460--467.

\bibitem{geobif} P. Piccione, A. Portaluri, D. V. Tausk, {\em Spectral Flow,
Maslov Index and  Bifurcation of
semi-Riemannian Geodesics},  Ann.\ Global Anal.\ Geom.\ \textbf{25} (2004), no.\ 2, 121--149.

\bibitem{london} P.\ Piccione, D.\ V.\ Tausk, {\em An Index Theorem for Non Periodic
Solutions of Hamiltonian Systems}, Proceedings of the London
Mathematical Society  (3) 83 (2001), 351--389.

\bibitem{topology} P.\ Piccione, D.\ V.\ Tausk, {\em The Morse Index Theorem
in semi-Riemannian Geometry}, Topology  {\bf41}  (2002),  no.\
6, 1123--1159.

\bibitem{calcvar} P.\ Piccione, D.\ Tausk, {\em An Index Theory for Paths that
are Solutions of a Class of Strongly Indefinite
Variational Problems},  Calculus of Variations and PDE's {\bf15}  (2002),
 no.\ 4, 529--551.

\bibitem{fechado} P.\ Piccione, D.\ V.\ Tausk, {\em On the Distribution
of Conjugate Points along semi-Riemannian Geodesics},
Communications in Analysis and Geometry {\bf11} (2003), No.\ 1, 33--48.
({\tt LANL math.DG/0011038})

\bibitem{rabier} P.\ J.\ Rabier, {\em Generalized Jordan chains
and two bifurcation theorems of Krasnosel'skii},
Nonlinear Anal.\ {\bf13} (1989), 903--934.

\bibitem{RobSal} J.\ Robbin, D.\ Salamon, {\em The Maslov Index
for Paths}, Topology {\bf32}, No.\ 4 (1993), 827--844.

\bibitem{RobSal2} J. Robbin, D. Salamon, {\em The Spectral Flow and
the Maslov Index}, Bull.\ London Math.\ Soc.\ {\bf27} (1995),
1--33.

\bibitem{SalZeh} D. Salamon, E. Zehnder, {\em Morse theory for periodic
solutions of Hamiltonian systems and the Maslov index}, Commun.\ Pure Appl.\ Math.\
{\bf45} (1992), 1303--1360.

\bibitem{Yoshi} T.\ Yoshida, {\em Floer Homology and Splitting of Manifolds},
Annals of Mathematics (2) {\bf134}, no.\ 2 (1991), 277--323.

\bibitem{War} F. W. Warner, {\em The Conjugate Locus of a Riemannian Manifold},
Amer.\ J.\ Math.\ {\bf87} (1965), 575--604.


\end{thebibliography}
\end{document}